\documentclass[a4paper,10pt]{article}

\usepackage[english]{babel}
\usepackage{geometry}
\usepackage{fancyhdr}
\usepackage{graphicx}
\usepackage{subfig}
\usepackage{booktabs}
\usepackage{multirow}
\graphicspath{{figures/}{./}}
\DeclareGraphicsExtensions{{.png},{.pdf},{.jpg}}
\usepackage{caption}
\usepackage{float}
\usepackage{pb-diagram,pb-xy}
\usepackage{bold-extra}
\usepackage[all,cmtip]{xy}
\usepackage{latexsym,verbatim}
\usepackage{amsfonts,amssymb}
\usepackage{amsmath}
\usepackage{amsthm}
\usepackage{amsxtra}
\usepackage{enumerate}
\usepackage{bm}
\usepackage[utf8]{inputenc}
\usepackage[T1]{fontenc}
\usepackage{cite}
\usepackage{hyperref}
\hypersetup{
	colorlinks=true,
	linkcolor=blue,
	citecolor=green
}
\usepackage{mathtools}
\DeclarePairedDelimiter\ceil{\bigg\lceil}{\bigg\rceil}

\usepackage{authblk}
\theoremstyle{definition}
\newtheorem{remark}{Remark}
\numberwithin{remark}{section}
\begin{document}

\title{Nonlinear lumped-parameter models for blood flow simulations in networks of vessels} 
\author[1]{Beatrice Ghitti}
\author[2]{Eleuterio F. Toro}
\author[1]{Lucas O. M\"{u}ller}
\affil[1]{ \small{Department of Mathematics, University of Trento, Trento, Italy; e-mail: \tt{beatrice.ghitti@unitn.it \& lucas.muller@unitn.it}} }
\affil[2]{ \small{Laboratory of Applied Mathematics DICAM, University of Trento, Trento, Italy; e-mail: \tt{eleuterio.toro@unitn.it}} }
\date{}
\maketitle
\thispagestyle{empty}

\begin{abstract}
\noindent To address the issue of computational efficiency related to the modelling of blood flow in complex networks, we derive a family of \textit{nonlinear} lumped-parameter models for blood flow in compliant vessels departing from a well-established one-dimensional model. These 0D models must preserve important nonlinear properties of the original 1D model: the nonlinearity of the pressure-area relation and the pressure-dependent parameters characterizing the 0D models, the resistance $R$ and the inductance $L$, defined in terms of a time-dependent cross-sectional area subject to pressure changes. We introduce suitable coupling conditions to join 0D vessels through 0D junctions and construct 0D networks preserving the original 1D network topology. The newly derived nonlinear 0D models are then applied to several arterial networks and the predicted results are compared against (i) the reference 1D results, to validate the models and assess their ability to reproduce good approximations of pressure and flow waveforms in all vessels at a much lower computational cost, measured in terms of CPU time, and (ii) the linear 0D results, to evaluate the improvement gained by including certain nonlinearities in the 0D models, in terms of agreement with the 1D results.\\

\noindent \textit{Keywords.} Blood flow; lumped-parameter models; nonlinearity; coupling; arterial networks.\\

\noindent \textbf{Mathematics Subject Classification.} 76-10; 35L65; 65M08.
\end{abstract}

\section*{Introduction}\label{sec:intro}

One-dimensional (1D) models and zero-dimensional (0D) or lumped-parameter models of blood flow, based on simplified representations of the components of the cardiovascular system, can contribute strongly to the study and the deep understanding of the circulatory physiology and pathology. These two families of models can be derived from the three-dimensional (3D), time-dependent, incompressible Navier-Stokes equations by exploiting specific features of blood flow, such as the basically cylindrical morphology of the vessels. Usually, in fluid-structure interaction (FSI) problems, the Navier-Stokes equations are also coupled to the equations of solid mechanics for the 3D structure, the deformable vessel wall.
These 3D models offer great level of detail and potentially accurate description of relevant quantities, but their numerical discretization is very challenging and requires high computational resources. Even though they are highly simplified with respect to the local dynamics, 1D and 0D models can then provide reasonably good approximations of more complex models at a much lower computational effort. Under suitable assumptions on the velocity profile and the tube law relating pressure and area, 1D models can provide accurate descriptions of wave propagation phenomena at a lower computational cost; 0D models can not resolve wave propagation and the spatial variation of physical quantities of interest, but they are cheaper than 3D and 1D models, representing simple mathematical objects with physical parameters that can be usually quantified from clinical measurements. The selection of the appropriate dimensionality in a model representation (from 0D to 3D) depends on problem-specific characteristics (like flow conditions and spatial/temporal scales), the aim of the application in mind and thus the required level of accuracy by which the physical process has to be described.\\
0D cardiovascular system representation and analysis started with the modelling of arterial flow using the well-known Windkessel model. The first and simplest 0D mono-compartment description is the famous two-element Windkessel model, which was first proposed by Stephen Hales in 1733, and later formulated mathematically by Otto Frank in 1899 \cite{Sagawa:1990a}. The Windkessel model consists of two parallel elements, a capacitor $C$, describing the storage properties of large arteries, and a resistor $R$, that describes the dissipative nature of small peripheral vessels including arterioles and capillaries. Later on, this approach was expanded to cover the modelling of other cardiovascular components, such as the heart, heart valves and veins, with either mono-compartment or multi-compartment 0D models, to simulate the global haemodynamics in the whole circulation system, by assuming a uniform distribution of the fundamental variables (pressure, flow and volume) within any 0D component of the system. 
In the multi-compartment approach, suitable $RLC$ models for a single vascular segment were then derived, as building blocks for constructing the entire vessel network model. In this regard, we refer to the works of Formaggia and Veneziani \cite{Formaggia:2003a}, Mili\v{s}i\'{c} and Quarteroni \cite{Milisic:2004a} and Formaggia \textit{et al.} \cite{FormaggiaQuarteroni:2009a}, where four typical compartment model configurations suited to the description of a vessel segment were derived, assuming that mean flow rate and pressure over the whole vessel are equivalent to either the input or the output values. Lumped-parameter models can also provide boundary conditions for 1D and local 3D models. For instance, in \cite{Alastruey:2008b}, the authors studied the effect of the parameters of the 0D outflow models on the waveforms propagated in an arterial network, with the aim of providing appropriate outflow 0D models for patient-specific simulations.\\
More recently, in a multiscale approach for modelling the whole cardiovascular system, as the one proposed by M\"{u}ller and Toro \cite{Mueller:2014a,Mueller:2014b}, 0D models have been used to model the heart, the pulmonary circulation and blood flow in arterioles, capillaries and venules, while 1D models have been adopted for large vessels of both arterial and venous networks. Indeed, 1D models have been extensively used to study wave propagation phenomena in arteries and, more recently, this has been extended also the venous circulation. A deep understanding of pressure and flow pulse wave propagation in the cardiovascular system and the impact of disease and anatomical variations on these propagation patterns can provide valuable information for clinical diagnosis and treatment \cite{Matthys:2007a,Blacher:2005a,Cruickshank:2002a,Fujimoto:2004a}.\\
However, modelling blood flow in highly complex networks, such as the global, closed-loop multiscale model for the whole cardiovascular system developed in M\"{u}ller and Toro \cite{Mueller:2014a,Mueller:2014b} or the ADAN (Anatomical Detailed Arterial Network) model proposed by Blanco  \textit{et al.} in \cite{Blanco:2014a,Blanco:2015a}, can result in computationally expensive simulations. The very high computational cost and execution time increase significantly when long time scales are to be simulated, taking several minutes per cardiac cycle. The situation is even more severe when, for instance, we want to face the modelling of several mechanisms that need to be integrated in the cerebral microcirculation, that are brain perfusion, waste clearance mechanisms and exchange of solutes between blood and various tissue beds.\\
Several works concerning lumped-parameter models to simulate arterial blood flow are found in the literature, which address the issues of execution time and optimization of topological complexity. In \cite{Fossan:2018a}, the authors presented a method to optimize/reduce the number of arterial segments included in 1D blood flow models and to find the model with the fewest number of necessary arteries for a given clinical application, by lumping distributed 1D segments into 0D Windkessel models, while preserving key features of flow and pressure waveforms. Similarly in \cite{Epstein:2015a}, in the context of patient-specific 1D blood flow modelling, Epstein  \textit{et al.} investigated the effect of a reduction in the number of arterial segments in a given distributed 1D model on the shape of the simulated pressure and flow waveforms, by systematically lumping peripheral 1D model branches into Windkessel models that preserve the net resistance and total compliance of the original model. In \cite{Safaei:2018a}, to address the issue of execution time and the question of granularity in the context of the modelling of the cerebral circulation, the authors proposed a lumped-parameter mathematical model, which was constructed using a bond graph formulation to ensure mass and energy conservation. In this work, the topology of the original 1D network was fully preserved and the model included arterial vessels with geometric and anatomical data based on the ADAN model \cite{Blanco:2015a}. Furthermore in \cite{Mirramezani:2020a}, a distributed lumped-parameter (DLP) modelling framework was proposed to efficiently compute blood flow and pressure in vascular domains at a computational cost that is orders of magnitude lower than that of computational fluid dynamics (CFD) simulations. By developing an expression of the generalized resistance $\mathcal{R}$, various sources of energy dissipation, including viscous dissipation, unsteadiness, flow separation, vessel curvature and bifurcations, were taken into account.\\
In the present work, the issue of computational efficiency and execution time, related to the modelling of blood flow in highly complex networks, has been faced by constructing lumped-parameter models departing from 1D models for blood flow in deformable vessels. The approach we propose is similar to the one presented in \cite{Safaei:2018a}, but here the main difference and novelty is that, even if they are simpler than the 1D models due to the fact that the space dependence is completely lost, these 0D models must preserve important properties of the original 1D models, such as the nonlinearity in the tube law relating the pressure to the cross-sectional area, and the dependence of the 0D model parameters on the time-varying vessel cross-section.\\
Indeed, it is important to emphasize that the mechanics of the cardiovascular system can exhibit strong nonlinearities. These nonlinear effects include, for example, the contribution of the convective acceleration terms in the momentum equation, the nonlinear relationship between pressure and volume which is observed in real vessels, the pressure-dependent vessel compliance, the collapse of vessels due to environmental pressure, \textit{etc}. For a complete review on these nonlinear phenomena characterizing the cardiovascular system see \cite{Shi:2011a}, for instance. Among these nonlinearities, pressure-dependent constitutive equations and vessel properties represent an example of great interest. Indeed, in the 0D model description, the values of the different components $C$, $L$ and $R$ are generally taken to be constant. However, since they represent real physical parameters, they are subject to the same nonlinearities, such as nonlinear constitutive material relations, as any other description of vascular mechanics \cite{Shi:2011a}. As the vessel diameter changes under changes of pressure, its compliance will change as will its resistance to flow. Furthermore, the vessel wall exhibits a nonlinear stress-strain curve \cite{Fung:1993a}, meaning that the compliance $C$ is also a function of the luminal pressure. These effects are typically included in 1D models, but neglected in 0D models. In particular, since the diameter changes in the arterial system are relatively small and the range of arterial pressures over the cardiac cycle is such that the material tends to operate in a relatively linear region of the stress-strain curve, it is a common practice to neglect the pressure dependence of the arterial properties \cite{Shi:2011a}. This is not true for veins, at least when they enter a collapsed state.\\
Several models based on \textit{in vivo} measurements or theoretical derivations have been proposed to describe the relation between cross-sectional area and internal pressure exhibited by vessels. Ursino  \textit{et al.} \cite{Ursino:1996a,Ursino:1999a} used an exponential curve to describe the nonlinear pressure-volume relationship in the peripheral and venous compartments in the simulation of carotid baro-regulation of pressure pulsation. Fogliardi  \textit{et al.} \cite{Fogliardi:1996a} compared the linear (\textit{i.e.} incorporating a constant compliance) and nonlinear (\textit{i.e.} including a nonlinear pressure-dependent compliance) formulation of the $RCR$ Windkessel model, concluding that no additional physiological information was gained when a pressure-dependent compliance was incorporated. In contrast, Li et al. \cite{Li:1990a,Berger:1992a} examined the consequences of incorporating a pressure-dependent compliance with exponential variation in a modified arterial system model given by a three-element $RCR$ vessel model and concluded that a pressure-dependent compliance could more accurately reflect the behaviour of the arterial system. Cappello  \textit{et al.} \cite{Cappello:1995a} developed a one-step computational procedure for estimating the parameters of the nonlinear three-element Windkessel model of the arterial system incorporating a pressure-dependent compliance. In \cite{Ursino:1997a}, the authors proposed a simple mathematical model of intracranial pressure dynamics, where the resistance in the the arterial-arteriolar cerebrovascular compartment was defined to be inversely proportional to the fourth power of inner radius.\\
However, when dealing with 0D models suited to the discretization of a single vascular segment \cite{Formaggia:2003a,Milisic:2004a,Safaei:2018a}, the nonlinearities characterizing the original 1D models have never been taken into account, assuming a linear relationship between the pressure and volume of the vascular compartment and constant parameters $C$, $L$ and $R$, which are usually defined in terms of a reference state $(A_0, P_0)$.\\
The main goal of the present work is then to derive lumped-parameter models for blood flow in deformable vessels in a way that: (i) important nonlinear properties of the 1D models are preserved; and (ii) when these 0D models are applied to a network of vessels, each 1D vessel is replaced by a ``0D vessel’’ by a per-segment (one-to-one 1D-to-0D) mapping, so that the topology of the original 1D network is naturally preserved.
Concerning the first point, we firstly require the relation between mean pressure and cross-sectional area to be nonlinear, as the corresponding tube law in the original 1D model. In this way, the constant compliance $C$ is replaced by the nonlinear pressure-area relation in which the mechanical properties of the vessel wall are embedded. In addition, the other components $L$ and $R$ of each 0D vessel are no longer taken to be constant, but are defined in terms of a time-dependent average cross-sectional area, in order to account for vessel cross-section changes under changes of pressure. This feature is expected to be relevant when modelling large deviations from a baseline state, such as hypertension, an haemorrhage or a collapsed state. Actually, remarkable differences will be also observed when comparing results obtained with the standard linear and the newly derived nonlinear 0D models even in physiologically normal conditions. Then, these results and observations prove that, by including these nonlinearities, we obtain more realistic and consistent 0D models with respect to the original 1D model.
Furthermore, the convective terms are not excluded \textit{a-priori} from these nonlinear 0D models. The dimensional analysis of the 1D equations and an exhaustive investigation of the contribution and relative importance of the convective terms in both 1D and 0D blood flow models will be crucial in deciding whether these terms can be reasonably neglected in the 0D models.\\
In order to validate and test these newly derived 0D models we reproduce several benchmark test cases proposed in \cite{Boileau:2015a}. We compare the 0D results obtained with the nonlinear 0D models against the original 1D model for different arterial networks, to assess the ability of such 0D models to produce reasonably good approximations of pressure and flow waveforms in all vessels of a network with respect to the reference 1D results. Furthermore, the nonlinear 0D results are also compared against the linear 0D results from the 0D models with linear pressure-area relation and constant parameters, to evaluate the improvement we obtain in the 0D results when including certain nonlinearities in the lumped-parameter models.\\

The paper is organized as follows. In Section \ref{sec:1Dmodel} we briefly introduce the governing equations of 1D blood flow in compliant vessels and we perform a dimensional analysis of these equations. Then, in Section \ref{sec:0Dmodel} a family of nonlinear 0D models is derived departing from the 1D model. First, we describe the derivation procedure and present the resulting system of ordinary differential equations (ODEs); then, we focus on the main features characterizing these lumped-parameter models, especially on the nonlinearity preserved in the pressure-area relation; finally, we conclude this section by considering the different 0D representations for a vessel segment depending on the different data prescribed at the inlet and outlet of the vessel. Afterwards, in Section \ref{sec:stability} we restrict to the linear case, namely to the standard 0D models with constant parameter, linear pressure-volume relation and without convective terms, to carry out a stability analysis of the corresponding ODE systems. Thereafter, in Section \ref{sec:0Djunction} we describe how to couple 0D vessels converging to a shared node (bifurcations/junctions of vessels) and how to couple 0D vessels to terminal Windkessel models. Finally, in Section \ref{sec:numexp} we perform several benchmark test problems by applying the derived family of 0D models to different arterial networks of increasing complexity and discuss the obtained results. We conclude with Section \ref{sec:conclusion}, where final remarks are made and perspectives for future work are outlined.

\section{One-dimensional (1D) blood flow model}\label{sec:1Dmodel}

A well-established formulation of the one-dimensional (1D) blood flow equations in deformable vessels is given by the following system
\begin{equation}\label{sys1D}
\left\{
\begin{aligned}
& \partial_t A+\partial_x q=0,\\
& \partial_t q +\partial_x\left(\alpha\frac{q^2}{A}\right) +\frac{A}{\rho}\partial_x p= f,
\end{aligned}
\right.
\end{equation}
where $x\in\left[ 0, l \right]$, with $l$ being the vessel length, is the axial coordinate along the longitudinal axis of the vessel and $t>0$ is time; $A(x,t)$ is the cross-sectional area of the vessel; $q(x,t)$ is the flow rate; $p(x,t)$ is the average internal pressure over a cross-section; $\alpha$ is a momentum correction factor, also called Coriolis coefficient, and $f(x,t) = -k_R \frac{q}{A}$ is the friction force per unit length, where $k_R > 0$ is the viscous resistance coefficient. Both parameter $\alpha$ and $k_R$ depend on the assumed velocity profile. Here, the following axisymmetric velocity profile was prescribed
\begin{equation}\label{velocityprofile}
u_x (x, \xi, t) = u(x,t) s(\xi) = u(x,t) \frac{\zeta + 2}{\zeta} \left[ 1 - \left(\frac{\xi}{r}\right)^{\zeta} \right],
\end{equation}
where $u_x (x, \xi, t)$ is the axial component of the fluid velocity, $u(x,t)=\frac{q}{A}$ is the mean velocity on each cross-section, $s(\xi)$ is the assumed velocity profile, $\xi$ is the vessel radial coordinate, $r(x,t)$ is the lumen radius and $\zeta$ is the velocity profile order. The viscous resistance per unit length of tube $k_R$ is defined as a function of the velocity profile as 
\begin{equation}
k_R = -\frac{\mu}{\rho} \int_{\partial S} \frac{\partial s}{\partial \vec{n}_{\partial S}} d\gamma,
\end{equation}
where $\partial S$ is the boundary of the vessel cross-section $S$ and $\vec{n}_{\partial S}$ is the outward normal vector to $\partial S$, and which, for the velocity profile chosen in (\ref{velocityprofile}), becomes
\begin{equation}
k_R = 2(\zeta + 2)\pi \frac{\mu}{\rho},
\end{equation}
where $\rho$ and $\mu$ are the constant blood density and viscosity, respectively. The momentum correction coefficient $\alpha$ is well-defined for unidirectional flow, namely
\begin{equation}
\alpha = \frac{\int_{S} u_x^2\, d\sigma}{A u^2},
\end{equation}
from which we have that the Coriolis coefficient $\alpha$ and the velocity profile order $\zeta$ are related by $\alpha = \frac{\zeta + 2}{\zeta +1}$. The value $\zeta = 9$, for which $\alpha = 1.1$, defines a flat velocity profile, which is especially valid for large arteries \cite{Hunter:1972}. The choice $\alpha = 1$, which indicates a completely flat velocity profile, is also commonly used since it simplifies the analysis of the resulting 1D model. In contrast, for a Poiseuille flow, the parabolic velocity profile is obtained by setting $\zeta=2$, for which $\alpha = \frac{4}{3}$.\\
Pressure $p(x,t)$ is related to the cross-sectional area $A(x,t)$ by the following algebraic relation
\begin{equation}\label{tubelaw}
p(x,t) - p_{ext}(x,t) = \psi(A(x,t); A_0, K, m, n, P_0), 
\end{equation}
with
\begin{equation}\label{transmuralpressure}
\psi(A(x,t); A_0, K, m, n, P_0) = K \left[ \left(\frac{A(x,t)}{A_0}\right)^m - \left(\frac{A(x,t)}{A_0}\right)^n \right] + P_0,
\end{equation}
where $p_{ext}(x,t)$ is the external pressure acting on the vessel and $P_0$ is the reference pressure at which $A=A_0$. The above relation describes the elastic deformation of the vessel wall with variations of the transmural pressure, assuming that viscoelastic effects are negligible. Pressure $p(x,t)$ also depends on the reference cross-sectional area $A_0$ and on parameters $K$, $m$ and $n$, which take into account geometrical and mechanical properties of the vessel.\\
In particular, if we assume all these parameters to be independent of $x$ and we consider arterial vessels, then the factor $K>0$ in (\ref{tubelaw}) denotes the arterial stiffness and it is modelled as in \cite{Toro:2016a,FormaggiaQuarteroni:2009a} by
\begin{equation}\label{kappa}
K\equiv K_a =\frac{\sqrt{\pi}h_0 E}{(1-\nu^2)\sqrt{A_0}},
\end{equation}
where $h_0$ is the vessel wall thickness, $E$ is the Young's modulus and $\nu$ is the Poisson ratio. We adopt $\nu=0.5$, which implies that the vessel wall is assumed to be incompressible. The parameters $m$ and $n$ are obtained from higher-order models or simply computed from experimental measurements. Typical values for arteries are $m=0.5$ and $n=0$.\\
We note that system (\ref{sys1D}) can be rewritten under the classical form of balance laws, that is
\begin{equation}\label{conservationlaws}
\partial_t  \bm{Q} +\partial_x \bm{F}(\bm{Q}) = \bm{S}(\bm{Q}),
\end{equation}
with
\begin{equation}
\bm{Q} = \left[ \begin{array}{c}
A\\[1ex]
q
\end{array} \right], 
\quad \bm{F}(\bm{Q}) = \left[ \begin{array}{c}
q\\[1ex]
\alpha \dfrac{q^2}{A} +\dfrac{K A}{\rho} \left[ \dfrac{m}{m+1} \left(\dfrac{A}{A_0}\right)^m - \dfrac{n}{n+1} \left(\dfrac{A}{A_0}\right)^n \right] 
\end{array} \right], 
\quad \bm{S}(\bm{Q})= \left[ \begin{array}{c}
0\\[1ex]
-k_R \dfrac{q}{A}
\end{array} \right],
\end{equation}
where $\bm{Q}$ is the vector of conserved variables, $\bm{F}$ is the flux function and $\bm{S}$ is the source term.\\
We introduce here also the wave speed, denoted by $c$, as follows
\begin{equation}\label{c}
c = \sqrt{\frac{A}{\rho} \frac{\partial p}{\partial A}}.
\end{equation}

\subsection{Dimensional analysis}\label{sec:1Dmodel:dimanalysis}

In order to assess the relative importance of each term in the 1D blood flow model (\ref{sys1D}), especially of convective, pressure and friction terms in the momentum balance equation, we perform here a dimensional analysis, similar to the analysis performed in \cite{Ghigo:2017a,Saito:2011a}. For this purpose, we introduce the following nondimensional variables
\begin{equation}\label{nondimvar}
t = T_0 \overline{t}, \quad  x = L_0 \overline{x}, \quad A = A_0 \overline{A}, \quad q = Q_0 \overline{Q} = \left(  A_0 U_0 \right) \overline{Q},
\end{equation}
where the constants $T_0$, $L_0$, $A_0$ and $U_0$ are orders of magnitude of the dimensional variables, so that nondimensional variables $\overline{t}$, $\overline{x}$, $\overline{A}$ and $\overline{Q}$ are of order $1$. In particular, $T_0$ is the time scale, $L_0$ is the longitudinal spatial scale, $A_0$ is the reference cross-sectional area and $U_0$ is a reference flow velocity.\\
The nondimensional equation of conservation of mass reads
\begin{equation}\label{nondimmasseq}
\frac{\partial \overline{A}}{\partial \overline{t}} + \left[ \frac{T_0 U_0}{L_0} \right] \frac{\partial \overline{Q}}{\partial \overline{x}} = 0.
\end{equation}
By rewriting the pressure gradient $\partial_x p$ in terms of the nondimensional variables (\ref{nondimvar}), after straightforward calculations we get the following nondimensional momentum balance equation
\begin{equation}\label{nondimmomentumeq}
\frac{\partial \overline{Q}}{\partial \overline{t}} + \gamma_C \frac{\partial}{\partial \overline{x}} \left( \alpha \frac{\overline{Q}^2}{\overline{A}}\right) + \gamma_P \frac{\partial \overline{A}}{\partial \overline{x}} = - \gamma_F\frac{\overline{Q}}{\overline{A}},
\end{equation}
where the three coefficients for the convective, pressure and friction terms have been introduced, respectively given by
\begin{equation}\label{coefficients:dimanalysis}
\gamma_C := \frac{T_0 U_0}{L_0}, \quad \gamma_P := \frac{T_0}{L_0 U_0} c^2, \quad \gamma_F := \frac{k_R T_0}{A_0}.
\end{equation}
The above coefficients are nondimensional quantities and their magnitudes indicate the relative importance of each of these terms in the momentum balance equation. In Section \ref{sec:numexp} we will exploit this kind of analysis to decide whether the convective terms can be neglected or not in the family of nonlinear 0D models we are going to derive in the next section.

\section{Derivation of zero-dimensional (0D) models}\label{sec:0Dmodel}

Here we extend the traditional approach of deriving lumped-parameter models for blood flow in a vascular segment \cite{Formaggia:2003a,Milisic:2004a,FormaggiaQuarteroni:2009a}, in a way to preserve certain properties (nonlinear characteristics) of the original 1D blood flow models. The proposed strategy to do so will be extensively described in this section, where we first derive a family of 0D model for a simple vascular compartment formed by a single vessel and then, by application of appropriate matching conditions obtained from conservation principles, we couple different 0D models to build more complex networks of vessels.

\subsection{Governing ODE system}\label{sec:0Dmodel:derivation}

First of all, given a vessel with $x \in \left[ x_L, x_R \right]$ of length $l = \lvert x_R - x_L \rvert$, we introduce the integral averages of the physical quantities of interest over the vessel length, as follows
\begin{equation}\label{meanvar}
\left\{ \begin{aligned}
& \widehat{Q}(t) = \frac{1}{l} \int_{x_L}^{x_R} q(x,t)\, dx\, :\, \textrm{mean (volumetric) flow rate,}\\
& \widehat{A}(t) = \frac{1}{l} \int_{x_L}^{x_R} A(x,t)\, dx\, :\, \textrm{mean cross-sectional area,}\\
\end{aligned}
\right.
\end{equation}
and we define the volume $V(t)$ of the vessel compartment as 
\begin{equation}\label{volume-area}
V(t) := \widehat{A}(t)\, l.
\end{equation}
Integrating in space the continuity equation in (\ref{sys1D}) over the interval $[x_L, x_R]$ leads to the following ordinary differential equation (ODE) in time for the volume $V(t)$
\begin{equation}\label{ode1}
\frac{d}{dt} V(t) + Q_R(t) - Q_L(t) = 0,
\end{equation}
where we have used definition (\ref{volume-area}) to rewrite the mass conservation equation in terms of the volume $V(t)$ and we have set
\begin{equation}
Q_L(t) = q(x_L, t), \quad Q_R(t) = q(x_R, t),
\end{equation}
to denote the flow at the inlet and outlet of the vessel, respectively.\\

When considering the momentum balance equation in (\ref{sys1D}), two simplifying assumptions are added in the standard approach of deriving 0D models \cite{Formaggia:2003a,Milisic:2004a,FormaggiaQuarteroni:2009a,Safaei:2018a}: (i) the contribution of the convective term $\partial_x (\alpha q^2/A)$ is neglected, assuming this term to be small compared to the other terms; (ii) the variation of $A$ with respect to $x$ is small compared to that of $p$ and $q$, replacing $A$ in the momentum equation with a constant value for the area, generally assumed to be the area at rest $A_0$. Indeed, the first assumption is particularly suited to represent the peripheral circulation, where blood flow is in general quite slow, while the second assumption is reasonable when the axial average is carried out over short segments.\\
However, in order to preserve certain important properties of the original 1D models in deriving the 0D models, we start integrating in space the momentum equation in (\ref{sys1D}) over the interval $[x_L, x_R]$ without considering the above simplifying assumptions. By including the contribution of the convective term, straightforward calculations yield
\begin{equation}\label{integralA}
\frac{d}{dt} \left( \frac{1}{l} \int_{x_L}^{x_R} q\, dx \right) + \frac{\alpha}{l} \left( \frac{Q_R(t)^2}{A_R(t)} - \frac{Q_L(t)^2}{A_L(t)} \right) + \frac{1}{l} \int_{x_L}^{x_R} \frac{A}{\rho} \partial_x p\, dx = - \frac{k_R}{l} \int_{x_L}^{x_R} \frac{q}{A}\, dx.
\end{equation}
Observing that space integrals of the pressure gradient and the viscous force depend on the area $A(x,t)$, we approximate the variable $A$ by its spatial average $\widehat{A}(t)$, rather than by a constant value $A_0$ as done in the traditional approach, and since this quantity is no longer space-dependent we can bring it outside of these integrals, to get
\begin{equation}
\frac{d}{dt}\widehat{Q}(t) + \frac{\alpha}{l} \left( \frac{Q_R(t)^2}{A_R(t)} - \frac{Q_L(t)^2}{A_L(t)} \right) + \frac{\widehat{A}(t)}{\rho l} \int_{x_L}^{x_R} \partial_x p\, dx = - \frac{k_R}{\widehat{A}(t)} \left( \frac{1}{l}\int_{x_L}^{x_R} q\, dx \right),
\end{equation}
that is
\begin{equation}
\frac{d}{dt}\widehat{Q}(t) + \frac{k_R}{\widehat{A}(t)} \widehat{Q}(t)  + \frac{\widehat{A}(t)}{\rho l} \left( P_R(t) - P_L(t) \right) + \frac{\alpha}{l} \left( \frac{Q_R(t)^2}{A_R(t)} - \frac{Q_L(t)^2}{A_L(t)} \right) = 0,
\end{equation}
where, again, we have set:
\begin{equation}
P_L(t) = p(x_L, t), \quad A_L(t) = A(x_L, t), \quad P_R(t) = p(x_R, t), \quad A_R(t) = A(x_R, t),
\end{equation}
to denote the upstream and downstream pressures and cross-sectional areas. We then introduce the following parameters
\begin{equation}\label{LandR}
\left\{ \begin{aligned}
& L(\widehat{A}) := \frac{\rho l}{\widehat{A}}\, : \, \textrm{\textit{nonlinear} inductance,}\\
& R(\widehat{A}) := \frac{\rho k_R l}{\widehat{A}^2}\, : \, \textrm{\textit{nonlinear} resistance,}\\
\end{aligned} 
\right.
\end{equation}
to obtain the final form of the 0D momentum equation
\begin{equation}\label{ode2}
L(\widehat{A}) \frac{d}{dt}\widehat{Q}(t) + R(\widehat{A}) \widehat{Q}(t) + P_R(t) - P_L(t) + \frac{\alpha \rho}{\widehat{A}} \left( \frac{Q_R(t)^2}{A_R(t)} - \frac{Q_L(t)^2}{A_L(t)} \right) = 0.
\end{equation}
Equations (\ref{ode1}) and (\ref{ode2}) are then collected together in the following system of ODEs
\begin{equation}\label{sys0D}
\left\{
\begin{aligned}
& \frac{d}{dt} V(t) + Q_R(t) - Q_L(t) = 0,\\
& L(\widehat{A}) \frac{d}{dt}\widehat{Q}(t) + R(\widehat{A}) \widehat{Q}(t) + P_R(t) - P_L(t) + \frac{\alpha \rho}{\widehat{A}(t)} \left( \frac{Q_R(t)^2}{A_R(t)} - \frac{Q_L(t)^2}{A_L(t)} \right) = 0.
\end{aligned}
\right.
\end{equation}
The state variables of the above system are the volume $V(t)$ of the vessel compartment and the mean flow rate $\widehat{Q}(t)$ over the vascular segment; the parameters characterizing the 0D model, defined in (\ref{LandR}), are $R(\widehat{A})$, which represents the resistance induced to the flow by the blood viscosity and depends on the chosen velocity profile, and $L(\widehat{A})$, which represents the inertial term in the momentum equation and is called inductance of the flow. These parameters are said to be \textit{nonlinear} in the sense that they do no longer depend on a constant reference cross-sectional area $A_0$, but on the time-dependent average cross-section $\widehat{A}(t)$, which in turn, as we will see in Section \ref{sec:0Dmodel:nonlinearity}, will depend on the mean pressure $P(t)$ acting on the vessel in a nonlinear way. We point out that it is important to define these vessel properties in terms of a time-varying cross-sectional area $\widehat{A}$, in order to account also for possible deviations from the baseline state, such as hypertension, vessel collapse or postural changes. We also remark that if in definition (\ref{LandR}) the time-dependent area $\widehat{A}(t)$ is replaced by the reference value $A_0$, then these parameters become constant and coincide with the constant parameters found in standard 0D models, namely
\begin{equation}
L_0 = L(A_0) = \frac{\rho l}{\pi r_0^2}, \quad R_0 = R(A_0) = \frac{\rho k_R l}{A_0^2} = \frac{2(\zeta + 2)\mu l}{\pi r_0^4}.
\end{equation}
The ODE system (\ref{sys0D}) also involves input and output quantities, $Q_L$, $P_L$ and $Q_R$, $P_R$, respectively, that need to be defined along with initial conditions in order to close problem (\ref{sys0D}). As we will see in Section \ref{sec:0Dmodel:0Dconfig}, depending on the different possible assumptions about the data prescribed at the inlet and at the outlet of the vessel, we will obtain four different configurations, all of them describing flow and volume/pressure dynamic in a single vascular segment. Note that these are not boundary conditions, since the continuous space dependence has been lost in the axial average.

\subsection{Nonlinearity}\label{sec:0Dmodel:nonlinearity}

System (\ref{sys0D}) describes the temporal evolution of volume $V(t)$ and mean flow rate $\widehat{Q}(t)$. At this point, we are interested in relating the system state variables to another important physical quantity, the mean pressure $\widehat{P}(t)$. Indeed, the nonlinear 0D model (\ref{sys0D}) is derived without making any assumptions about the pressure law relating the mean pressure $\widehat{P}(t)$ to the average cross-sectional area $\widehat{A}(t)$, or, equivalently, to the volume $V(t)$ and the 0D mass conservation equation is obtained for the volume, not for the pressure. Therefore, in the following, we are going to characterize the relation to compute the pressure $\widehat{P}(t)$ from the area $\widehat{A}(t)$.\\
In the original 1D blood flow model, the pressure $p(x,t)$ is related to the cross-sectional area $A(x,t)$ by the elastic tube law (\ref{tubelaw}), which is a nonlinear relationship describing the behaviour of vessel walls in response to changes in the transmural pressure. On the one hand, in the traditional approach of deriving 0D models, where the convective terms are neglected and the model parameters $L_0$ and $R_0$ are constant, pressure $\widehat{P}$ and volume $V$ are linearly related via the constant compliance $C_0$, as follows
\begin{equation}\label{linearPV}
\widehat{P}(t) = P_0 + \frac{V(t) - V_0}{C_0} + P_{ext}(t),
\end{equation}
with 
\begin{equation}\label{constantC}
C_0 := l \left( \frac{\partial A}{\partial p} \right)\bigg|_{A=A_0} = l \frac{2A_0}{K} = \frac{3}{2} \frac{\pi r_0^3 l}{E h_0},
\end{equation}
where the last two expressions of $C_0$ are specific for arteries. Since coefficient $C_0$, which represents the mass storage capacity due to the compliance of the vessel, is constant, the nonlinearity of the 1D pressure-area relation is lost in the resulting 0D models. Indeed, in the traditional approach, the 0D mass conservation equation is derived for the pressure $\widehat{P}$, and not for the volume $V$, by integrating in space both the 1D continuity equation and the 1D tube law under suitable assumptions (for the complete derivation, see for instance \cite{FormaggiaQuarteroni:2009a}). This procedure leads to the linear pressure law stated in equation (\ref{linearPV}), where the nonlinearity of the original 1D pressure-area relation is completely lost.\\
Here we propose to directly compute the pressure $\widehat{P}$ from the average cross-sectional area $\widehat{A}$ via the nonlinear tube law (\ref{tubelaw}), given by
\begin{equation}
\widehat{P}(t) - P_{ext}(t) = K \left[ \left(\frac{\widehat{A}(t)}{A_0}\right)^m - \left(\frac{\widehat{A}(t)}{A_0}\right)^n \right] + P_0,
\end{equation} 
which in the case of arteries turns out to be
\begin{equation}\label{nonlinearPA}
\widehat{P}(t) - P_{ext}(t) = \dfrac{K}{\sqrt{A_0}} \left[ \sqrt{\widehat{A}(t)} - \sqrt{A_0} \right] + P_0,
\end{equation} 
so that the nonlinearity of the original 1D pressure-area relation of the vessel is fully preserved also in the family of derived 0D models. In particular, the 0D model (\ref{sys0D}) provides the temporal dynamic of volume $V$, from which the time-dependent average cross-sectional area $\widehat{A}$ can be easily computed as $\widehat{A}(t)=V(t)/l$, to get the mean pressure $\widehat{P}(t) \equiv \widehat{P}(\widehat{A}(t))$ via the nonlinear tube law (\ref{nonlinearPA}). The derived family of 0D models is then said to be \textit{nonlinear}, together with the fact that parameters $L$ and $R$ characterizing these 0D models depend on the time-dependent cross-section $\widehat{A}(t)$, rather than on the constant reference value $A_0$.\\
Rewriting the linear relation (\ref{linearPV}) in terms of the area $\widehat{A}$ and replacing the explicit expression for the constant compliance $C_0$ leads to
\begin{equation}\label{linearPA}
\widehat{A} = A_0 + \frac{C_0}{l} \left(\widehat{P}-P_0-P_{ext}\right) = A_0 \left[ 1 + \frac{2}{K}\left(\widehat{P}-P_0-P_{ext}\right) \right],
\end{equation}
namely, the average cross-sectional area $\widehat{A}$ depends linearly on pressure $\widehat{P}$. In contrast, from the nonlinear tube law (\ref{nonlinearPA}), this dependence is kept quadratic, as follows
\begin{equation}\label{nonlinearAP}
\widehat{A} = A_0 \left[ 1 +\frac{\widehat{P}-P_0-P_{ext}}{K} \right]^2 = A_0 \left[ 1 + \frac{2}{K} \left(\widehat{P}-P_0-P_{ext}\right) +\frac{\left(\widehat{P}-P_0-P_{ext}\right)^2}{K^2} \right].
\end{equation}
As a consequence, by adopting relation (\ref{linearPV}) or, equivalently, (\ref{linearPA}), we are neglecting the second-order term coming from the tube law (\ref{nonlinearPA}) and thus the nonlinearity of the original 1D pressure-area relation is lost. This nonlinearity is then fully preserved if, in the derived 0D models, the pressure $\widehat{P}$ is computed from the area $\widehat{A}$ via the nonlinear tube law (\ref{nonlinearPA}).\\

In conclusion, the ODE system (\ref{sys0D}) involves the mean flow rate $\widehat{Q}$, the volume $V$ and the mean pressure $\widehat{P}$ over the entire vessel, where the pressure $\widehat{P}$ is related to the average cross-sectional area $\widehat{A}$, and thus to $V$, via the nonlinear tube law (\ref{nonlinearPA}). Furthermore, system (\ref{sys0D}) also depends on the input and output quantities exchanged by the vessel with the rest of the systems, namely $P_{L}$, $Q_{L}$, which from now on will be denoted by $P_{in}$, $Q_{in}$, and $P_{R}$, $Q_{R}$, which will be instead replaced by $P_{out}$, $Q_{out}$. As we will illustrate in Section \ref{sec:0Dmodel:0Dconfig}, some input/output data, along with initial conditions, need to be prescribed in order to close system (\ref{sys0D}).\\
For the sake of simplicity, from now on we will denote the flow rate $\widehat{Q}$ and pressure $\widehat{P}$ just by $Q$ and $P$, respectively.

\subsection{0D vessel configurations}\label{sec:0Dmodel:0Dconfig}

The ODE system (\ref{sys0D}) defines a family of nonlinear 0D models. Indeed, four different 0D models are obtained depending on the different possible assumptions about the data prescribed at the inlet and outlet of the vessel. These models determine all the possible configurations of the same 0D vessel, which are the $(P_{in}, Q_{out})$, $(Q_{in}, P_{out})$, $(P_{in}, P_{out})$ and $(Q_{in}, Q_{out})$-type 0D vessels, all of them describing flow and volume/pressure dynamic in a compliant vessel.\\
In the following, we will discard the contribution of the convective component in the momentum balance equation, originally included in the general formulation of the family of nonlinear 0D models (\ref{sys0D}). This choice, commonly adopted in the literature, will be fully justified in Section \ref{sec:numexp}, where we will extensively discuss whether it is reasonable or not to incorporate the contribution of the convective terms in the 0D models, also observing that including these terms in the 0D models is not straightforward as one would expect. Table \ref{table:equations0Dmodels} summarizes the different 0D vessel configurations and the corresponding ODE systems, which will be described in detail throughout the remaining of this section.

\subsubsection{$(P_{in}, Q_{out})$-type 0D vessel}\label{sec:0Dmodel:0Dconfig:PinQout}

Suppose that the data prescribed at the inlet and outlet of the vessel are $P_{in}$ and $Q_{out}$, respectively. This first 0D vessel type, the $(P_{in}, Q_{out})$-type vessel, is displayed in the first row of Table \ref{table:equations0Dmodels}. Then, the temporal dynamic of the state variables $Q$ and $V$, which are the unknowns under time derivative, is governed by the following system of ODEs
\begin{equation}\label{PinQout}
\left\{ \begin{aligned}
& \frac{d V}{dt} = Q - Q_{out},\\
& \frac{d Q}{dt} = \frac{1}{L(\widehat{A})}\left[ P_{in} - R(\widehat{A}) Q - P \right],
\end{aligned}
\right.
\end{equation}
where, in the above momentum balance equation, the mean pressure $P$ depends on the time-dependent area $\widehat{A}$ via the nonlinear relation (\ref{nonlinearPA}). Clearly, for this 0D model, we have $Q_{in} = Q$. Given the nonlinear resistance $R_{tot}(\widehat{A})$ of the entire vessel according to formula (\ref{LandR}), this total resistance has been split and equally distributed into two resistances in series, $R$ and $R_d$, in order to add the distal resistance $R_d = R_{tot}(\widehat{A})/2$ at the outlet of the vessel, as shown in Table \ref{table:equations0Dmodels}. Then, the outlet pressure $P_{out}$ is directly computed as
\begin{equation}
P_{out}(t) = P(t) - R_d Q_{out}(t).
\end{equation}
The above value of the pressure $P_{out}$ at the outlet of the vessel, obtained by splitting the total vessel resistance and adding a distal resistance to the vessel, can then be used to enforce the continuity of pressure either at 0D junctions, or in the coupling with terminal elements, as it will be described in Section \ref{sec:0Djunction}. For the same motivation, a proximal resistance will be added to the $(Q_{in}, P_{out})$-type 0D vessel and both a proximal and a distal resistance will be appended to the $(Q_{in}, Q_{out})$-type 0D vessel, as illustrated in Sections \ref{sec:0Dmodel:0Dconfig:QinPout} and \ref{sec:0Dmodel:0Dconfig:QinQout}, respectively. We observe that this choice does not involve additional hypothesis on the flow since the total vessel resistance to flow is kept the same and it is just split into two (or more) resistances in series.\\
We observe that, by adding the usual simplifying assumptions considered in the standard approach of deriving 0D models, namely that the model parameters $L$ and $R$ are constant, and pressure $P$ and volume $V$ are linearly related via the constant compliance $C_0$, then the well-established formulation of the \textit{linear} 0D blood flow model is restored, as follows
\begin{equation}\label{PinQoutLinear}
\left\{ \begin{aligned}
& \frac{d V}{dt} = Q - Q_{out},\\
& \frac{d Q}{dt} = \frac{1}{L_0}\left[ P_{in} - R_0 Q - P \right],
\end{aligned}
\right.
\end{equation}
where, as discussed in Section \ref{sec:0Dmodel:nonlinearity}, the following linear relation between $P$ and $V$ holds
\begin{equation}\label{linearPVnew}
P - P_{ext} = P_0 + \frac{V-V_0}{C_0}.
\end{equation}
The $(P_{in},Q_{out})$-type vessel described so far is displayed in the first row of Table \ref{table:equations0Dmodels}. This representation is precisely valid for the linear system (\ref{PinQoutLinear}), where the model parameters are constant, the convective terms are neglected and the pressure $P$ is linearly related to the volume $V$ via the compliance $C_0$. However, this description can still be conveniently used also for the nonlinear 0D model (\ref{PinQout}): the model parameters $L$ and $R$ are nonlinear, while the compliance $C$ is now replaced by the nonlinear tube law (\ref{nonlinearPA}) relating the mean pressure $P$ and the average cross-sectional area $\widehat{A}$, in which the mechanical properties of the vessel wall are embedded. In general, the $C$-element represents the elastic component of the vessel regardless of how pressure and area are related.

\paragraph{\normalfont\textit{Formulation of the mass conservation equation.}}
 
By computing the time derivative of both sides of the linear relation (\ref{linearPVnew}), the mass conservation equation can be rewritten in terms of pressure $P$, as follows
\begin{equation}\label{dPdtLinear}
\frac{d P}{dt} = \frac{Q - Q_{out}}{C_0}.
\end{equation}
Clearly, this equivalence between the two formulations of the 0D mass conservation equation, the one in (\ref{PinQoutLinear}) describing the dynamic of $V$ and the other (\ref{dPdtLinear}) the dynamic of $P$, is no longer true for the nonlinear 0D model (\ref{PinQout}), since now pressure $P$ depends in a nonlinear fashion on the cross-sectional area $\widehat{A}$, and thus on volume $V$, via the tube law (\ref{nonlinearPA}). Indeed, on the one hand, the continuity equation in system (\ref{PinQout}), describing the time-variation of volume $V$, is an exact relation obtained by directly integrating the 1D equation $\partial_t A + \partial_x q = 0$ over the vessel length $l = |x_R - x_L|$. On the other hand, using the fact that
\begin{equation*}
\frac{\partial A}{\partial t} = \frac{\partial A}{\partial p} \frac{\partial p}{\partial t},
\end{equation*}
the mass conservation equation in system (\ref{sys1D}) integrated along the axial direction can be rewritten as follows
\begin{equation}
\frac{1}{l} \int_{x_L}^{x_R} \frac{\partial A}{\partial p} \frac{\partial p}{\partial t}\, dx + \frac{1}{l} \left( Q_R(t) - Q_L(t)\right) = 0.
\end{equation}
To compute the integral in the above equation, we assume $\partial A/\partial p$ to be evaluated at $A=\widehat{A}$, being $\widehat{A}(t)$ the time-dependent average cross-sectional area of the vessel, namely we introduce the following approximation
\begin{equation}\label{approxdAdp}
\frac{\partial A}{\partial p} \approx \left( \frac{\partial A}{\partial p}\right) \bigg|_{A=\widehat{A}},
\end{equation}
so that this quantity is no longer space-dependent and can be brought outside of the integral, to get
\begin{equation}
\frac{1}{l} \left(\frac{\partial A}{\partial p} \right)\bigg|_{A=\widehat{A}} \int_{x_L}^{x_R} \frac{\partial p}{\partial t}\, dx + \frac{1}{l} \left( Q_R(t) - Q_L(t)\right) = 0,
\end{equation}
which can be finally rewritten as
\begin{equation}\label{dPdt}
C(\widehat{A}) \frac{d P(t)}{dt} + Q_R(t) - Q_L(t) = 0,
\end{equation}
where we have introduced the following nonlinear parameter
\begin{equation}\label{C}
C(\widehat{A}) := l \left(\frac{\partial A}{\partial p} \right)\bigg|_{A=\widehat{A}}\, :\, \text{\textit{nonlinear} compliance}.
\end{equation}
The parameter $C(\widehat{A})$ represents the vessel wall compliance and, like the parameters $L$ and $R$ defined in (\ref{LandR}), is said to be \textit{nonlinear}, in the sense that it depends on the time-dependent area $\widehat{A}(t)$. In the case of the $(P_{in}, Q_{out})$-type 0D vessel under study, the approximate equation (\ref{dPdt}) becomes
\begin{equation}\label{dPdtNonlinear}
C(\widehat{A}) \frac{d P}{dt} = Q - Q_{out}.
\end{equation}
The shape of this equation strongly recalls that of equation (\ref{dPdtLinear}) obtained in the linear case, but now, because of the approximation introduced in (\ref{approxdAdp}), it is no longer equivalent to the first exact formulation of the continuity equation in (\ref{PinQout}). For this reason, in the nonlinear 0D model, we will restrict ourselves to consider the mass conservation equation in terms of the volume $V$, as given in (\ref{PinQout}), which is exact, while the pressure $P$ will be always computed from the nonlinear relation (\ref{nonlinearPA}) in order to fully preserve the nonlinearity of the original 1D tube law.

\subsubsection{$(Q_{in}, P_{out})$-type 0D vessel}\label{sec:0Dmodel:0Dconfig:QinPout}

Suppose now that $Q_{in}$ and $P_{out}$ are given data. The 0D vessel configuration corresponding to these input data is displayed in the second row of Table \ref{table:equations0Dmodels} and the governing equations of the nonlinear 0D model for the $(Q_{in}, P_{out})$-type vessel read
\begin{equation}\label{QinPout}
\left\{ \begin{aligned}
& \frac{d V}{dt} = Q_{in} - Q,\\
& \frac{d Q}{dt} = \frac{1}{L(\widehat{A})} \left[ P - R(\widehat{A}) Q - P_{out} \right].
\end{aligned}
\right.
\end{equation}
Here, we clearly have $Q_{out} = Q$. Furthermore, also for this 0D vessel type, the nonlinear resistance $R_{tot}(\widehat{A})$ of the entire vessel, given in formula (\ref{LandR}), has been split into two equal resistances in series, $R_p$ and $R$, in order to add the proximal resistance $R_p = R_{tot}(\widehat{A})/2$ at the inlet of the vessel and to explicitly compute the inlet pressure $P_{in}$ as 
\begin{equation}
P_{in}(t) = P(t) + R_p Q_{in}(t),
\end{equation}
which will be used in the coupling procedure between 0D vessels at 0D junctions and in the coupling of a 0D vessel to terminal elements, in order to enforce the continuity of pressure.

\subsubsection{$(P_{in}, P_{out})$-type 0D vessel}\label{sec:0Dmodel:0Dconfig:PinPout}

With the $(P_{in}, Q_{out})$ and $(Q_{in}, P_{out})$-type vessels at hand, the last two 0D vessel types can be easily constructed by connecting two basic 0D configurations described so far, as can be clearly seen from Table \ref{table:equations0Dmodels}. Indeed, if pressure is prescribed at both inlet and outlet of the vessel, $P_{in}$ and $P_{out}$, respectively, the corresponding system can be modelled by connecting a $(P_{in}, Q_{out})$-type 0D vessel to a $(Q_{in}, P_{out})$-type 0D vessel, yielding the configuration illustrated in the third row of Table \ref{table:equations0Dmodels}. Then, the nonlinear 0D model is governed by the following ODE system
\begin{equation}\label{PinPout}
\left\{ \begin{aligned}
& \frac{d V}{dt} = Q - Q_d,\\
& \frac{d Q}{dt} = \frac{1}{L(\widehat{A})} \left[ P_{in} - R(\widehat{A}) Q - P \right],\\
&  \frac{d Q_d}{dt} = \frac{1}{L_d(\widehat{A})} \left[ P - R_d(\widehat{A}) Q_d - P_{out} \right],
\end{aligned}
\right.
\end{equation}
where the quantities $Q$ and $Q_d$ define the flow rates through the first proximal and the second distal parts of the 0D vessel, respectively. By construction, the total resistance $R_{tot}(\widehat{A})$ and inductance $L_{tot}(\widehat{A})$, as defined in (\ref{LandR}), are equally distributed between these two proximal and distal portions of the segment, namely into $(R, L)$ and $(R_d, L_d)$, as shown in Table \ref{table:equations0Dmodels}.

\subsubsection{$(Q_{in}, Q_{out})$-type 0D vessel}\label{sec:0Dmodel:0Dconfig:QinQout}

Finally, assuming that both flow rates $Q_{in}$ and $Q_{out}$ are prescribed yields the last possible configuration, that is the $(Q_{in}, Q_{out})$-type 0D vessel, obtained by connecting a $(Q_{in}, P_{out})$-type 0D vessel to a $(P_{in}, Q_{out})$-type 0D vessel and displayed in the bottom row of Table \ref{table:equations0Dmodels}. In this case, the resulting system of ODEs reads
\begin{equation}\label{QinQout}
\left\{ \begin{aligned}
& \frac{d V}{dt} = Q_{in} - Q,\\
& \frac{d Q}{dt} = \frac{1}{L(\widehat{A})} \left[ P - R(\widehat{A}) Q - P_d \right],\\
& \frac{d V_d}{dt} = Q - Q_{out},
\end{aligned}
\right.
\end{equation}
where the quantities $V$, $P$ and $V_d$, $P_d$ define volumes and pressures in the first proximal and in the second distal compartments of the vascular segment, respectively.
The nonlinear resistance $R_{tot}(\widehat{A})$ over the entire vessel has been split into a proximal resistance $R_p$ at the inlet of the vessel, a resistance $R$ between the two capacitors and a distal resistance $R_d$ at the outlet of the vessel, as shown in Table \ref{table:equations0Dmodels}, so that the inlet and outlet pressures $P_{in}$ and $P_{out}$ are computed as
\begin{equation}
P_{in}(t) = P(t) + R_p Q_{in}(t), \qquad P_{out}(t) = P_d(t) - R_d Q_{out}(t).
\end{equation}
Admissible choices are also to set either $R_d=0$, so that we just have $P_{out} = P_d$, or $R_p=0$, that implies $P_{in} = P$. Also for this 0D vessel configuration, the inlet and outlet pressures $P_{in}$ and $P_{out}$ will be used, as described in Section \ref{sec:0Djunction}, to couple 0D vessels in a network by enforcing the conservation of mass and the continuity of pressure.\\
\begin{table}[h!]\small
	\centering
	\renewcommand\arraystretch{1.2}
	\begin{tabular}{l|c|c} 
		\toprule
		\textbf{0D vessel} & \textbf{ODE system} & \textbf{Representation} \\
		\midrule
		$(P_{in}, Q_{out})$ &
		$\left\{ \begin{aligned}
		& \dfrac{d V}{dt} = Q - Q_{out}\\
		& \dfrac{d Q}{dt} = \dfrac{1}{L(\widehat{A})}\left[ P_{in} - R(\widehat{A}) Q - P \right]
		\end{aligned}
		\right.$
		& \raisebox{-0.5\height}{\includegraphics[scale=0.65]{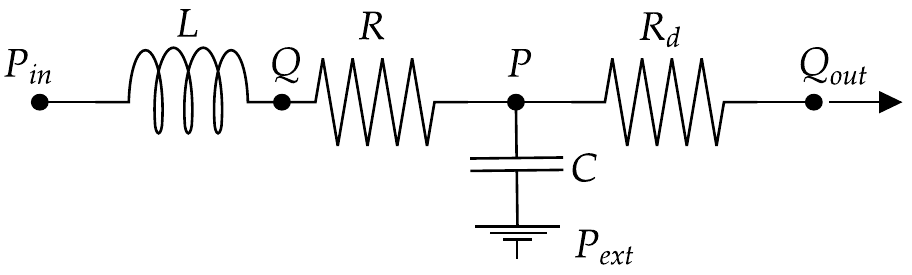}} \\
		\midrule
		$(Q_{in}, P_{out})$ &
		$\left\{ \begin{aligned}
		& \dfrac{d V}{dt} = Q_{in} - Q\\
		& \dfrac{d Q}{dt} = \dfrac{1}{L(\widehat{A})} \left[ P - R(\widehat{A}) Q - P_{out} \right]
		\end{aligned}
		\right.$
		&  \raisebox{-0.5\height}{\includegraphics[scale=0.65]{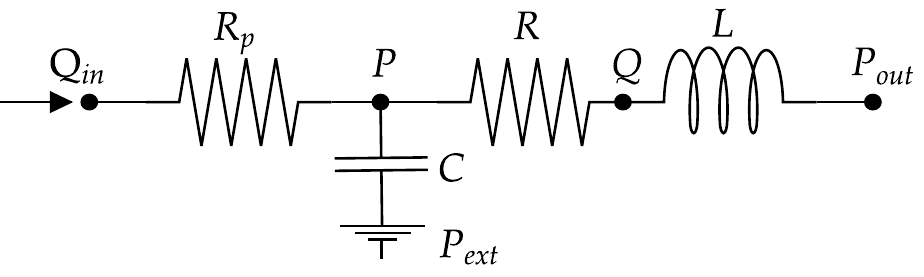}} \\
		\midrule
		$(P_{in}, P_{out})$ & 
		$\left\{ \begin{aligned}
		& \dfrac{d V}{dt} = Q - Q_d\\
		& \dfrac{d Q}{dt} = \dfrac{1}{L(\widehat{A})} \left[ P_{in} - R(\widehat{A}) Q - P \right]\\
		&  \dfrac{d Q_d}{dt} = \dfrac{1}{L_d(\widehat{A})} \left[ P - R_d(\widehat{A}) Q_d - P_{out} \right]
		\end{aligned}
		\right.$
		& \raisebox{-0.5\height}{\includegraphics[scale=0.65]{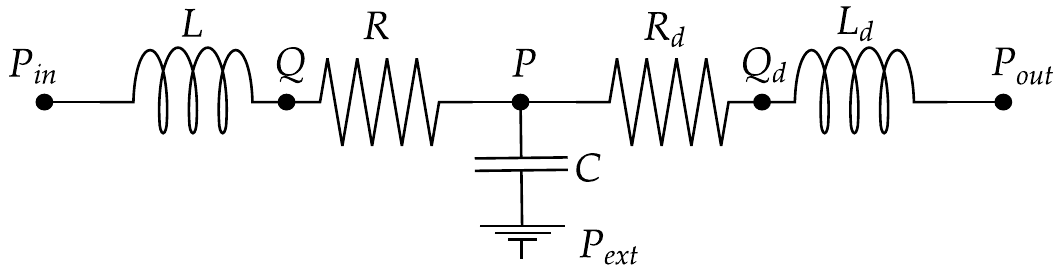}} \\
		\midrule
		$(Q_{in}, Q_{out})$ &
		$\left\{ \begin{aligned}
		& \dfrac{d V}{dt} = Q_{in} - Q\\
		& \dfrac{d Q}{dt} = \dfrac{1}{L(\widehat{A})} \left[ P - R(\widehat{A}) Q - P_d \right]\\
		& \dfrac{d V_d}{dt} = Q - Q_{out}
		\end{aligned}
		\right.$ 
		& \raisebox{-0.5\height}{\includegraphics[scale=0.65]{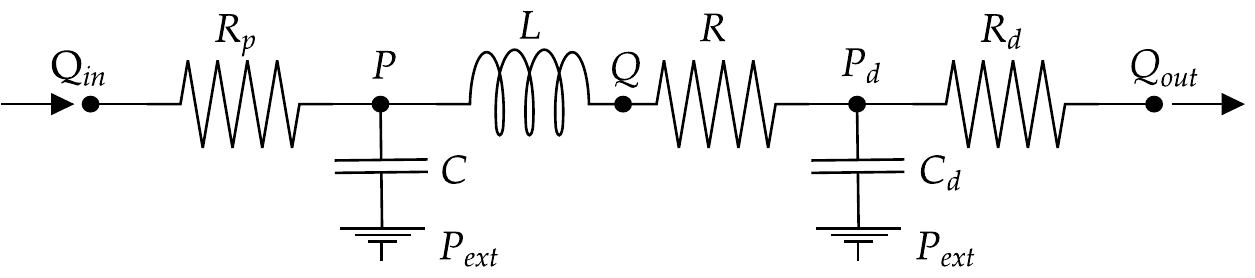}} \\
		\bottomrule
	\end{tabular}
	\caption{Summary of the four possible 0D vessel configurations depending on the different data prescribed at the inlet and outlet of the vessel: for each 0D vessel type, we report the associated governing ODE system (middle column) and display its representation in electric circuit notation (right column). In the linear case, the vessel wall compliance $C$ is a constant parameter, as defined in (\ref{constantC}); in the nonlinear case, $C$ represents the mechanical properties of the vessel wall embedded in the nonlinear pressure-area relation (\ref{nonlinearPA}).} \label{table:equations0Dmodels}
\end{table}

\section{Stability analysis} \label{sec:stability}

In this section, we are interested in studying the stability properties of the systems of ODEs governing the 0D blood flow models. In particular, for each of the four different 0D vessel configurations presented in Section \ref{sec:0Dmodel:0Dconfig}, we will perform the stability analysis of the corresponding linear ODE system with constant parameters $C_0$, $R_0$ and $L_0$, linear pressure-volume relation (\ref{linearPVnew}) and without convective terms. These systems of ODEs are linear and inhomogeneous, with periodic forcing terms. We point out that, to the best of our knowledge, such a stability analysis to investigate the behaviour of the exact solution of an ODE system arising from lumped-parameter models for blood flow has never been reported before, in the open literature.\\
When including the convective part into the 0D models, even if the model parameters are still constant and the pressure-volume relation is kept linear, the convective terms introduce a nonlinear component in the corresponding ODE systems. It is well-known that the theory on the stability of systems of ODEs is strictly valid only in the case in which the ODE system is linear. Indeed, in the nonlinear case, the eigenvalues of the Jacobian matrix associated to the ODE system can not be used to describe the behaviour of the exact solution of the original problem. The analytical study of the stability properties of the complete ODE systems including the convective terms turns out to be extremely complicated, if not impossible. Therefore, we limit our study to the linear case without convective terms. In general, we observed that results obtained for the linear case are valid also for the nonlinear case, as confirmed by numerical experiments presented in Section \ref{sec:numexp}. Moreover, in that section, we comment on numerical findings that suggest that the incorporation of convective terms has very strong implications on the stability of the resulting ODE global system, which in turn results in an extremely high computational cost and lack of robustness of its numerical treatment.\\

We first consider the $(P_{in}, Q_{out})$-type 0D vessel displayed in the top row of Table \ref{table:equations0Dmodels} and, by adding the following assumptions:
\begin{itemize}
	\item the model parameters $L_0$ and $R_0$ are constant,
	\item pressure $P$ is linearly related to volume $V$ via the constant compliance $C_0$ according to (\ref{linearPVnew}),
\end{itemize}
the resulting system of ODEs governing such a 0D vessel configuration reads
\begin{equation}\label{PinQoutLin1}
\left\{ \begin{aligned}
& \frac{d V}{dt} = Q - Q_{out},\\
& \frac{d Q}{dt} = \frac{1}{L_0} \left[ P_{in} - R_0 Q - P \right],
\end{aligned}
\right.
\end{equation}
where $R_0$ denotes the constant resistance between $Q$ and $P$. The data prescribed at the inlet and outlet of the vessel, $P_{in} \equiv P_{in}(t)$ and $Q_{out} \equiv Q_{out}(t)$, respectively, are given time-dependent functions, that we assume to be periodic of a certain period $T_0 > 0$. By using relation (\ref{linearPVnew}), the momentum equation in (\ref{PinQoutLin1}) can be reformulated in terms of the state variable $V$ as follows
\begin{equation}\label{PinQoutLin2}
\left\{ \begin{aligned}
& \frac{d V}{dt} = Q - Q_{out},\\
& \frac{d Q}{dt} = \frac{1}{L_0}\left[ P_{in} - R_0 Q - \frac{V}{C_0}  \right],
\end{aligned}
\right.
\end{equation}
where, for the sake of simplicity in the notation, in (\ref{linearPVnew}) we have set $V_0 = P_0 = P_{ext}=0$. Then, the above ODE system can be rewritten in matrix form as
\begin{equation}\label{matrixForm}
\frac{d \bm{x}(t)}{dt} =  \bm{A} \bm{x}(t) + \bm{b}(t),
\end{equation}
where we have set
\begin{equation}
\bm{x} (t) = \left[ \begin{matrix}
V(t)\\[2ex] Q(t)
\end{matrix}\right],
\quad
\bm{A} = \left[ \begin{matrix}
0 & 1\\[2ex] -\dfrac{1}{C_0 L_0} & -\dfrac{R_0}{L_0}
\end{matrix} \right],
\quad
\bm{b}(t) =  \left[ \begin{matrix}
-Q_{out}(t)\\[2ex] 
\dfrac{P_{in}(t)}{L_0}
\end{matrix}\right].
\end{equation}
Namely, $\bm{x}(t)$ is the vector of unknowns, the model state variables $V$ and $Q$, $\bm{A}$ is the constant coefficient matrix and $\bm{b}(t)$ is the time-dependent vector periodic forcing function, providing external data to the system. As the coefficient matrix $\bm{A}$ is constant, we have a non-homogeneous linear system of ODEs with constant coefficients.\\
The stability of the exact solution of the complete ODE system (\ref{matrixForm}) is determined by the real part of the eigenvalues of the coefficient matrix $\bm{A}$. In particular, we are going to use the following two results: 
\begin{itemize}
	\item[(i)]  Given a linear homogeneous system of ODEs with constant coefficients, that is an ODE system of the form (\ref{matrixForm}) with null forcing function $\bm{b}(t)\equiv 0$, a necessary and sufficient condition for this system to be \textit{asymptotically stable} is that all eigenvalues of $\bm{A}$ have strictly negative real part. 
	\item[(ii)] If the forcing function $\bm{b}(t)$ is periodic and if the homogeneous part of system (\ref{matrixForm}) is asymptotically stable, then the exact solution of the original inhomogeneous problem will converge to the periodic solution of system (\ref{matrixForm}) as $t$ increases, for any admissible choice of the initial condition.
\end{itemize}
These are well-known results and further details and proofs can be found in Appendix \ref{appendix:stabilityanalysis:part1}.\\
The eigenvalues $\lambda_{1,2}$ associated to matrix $\bm{A}$ are the roots of the following second-degree characteristic polynomial
\begin{equation}\label{polynomial}
\lambda^2 + \frac{R_0}{L_0} \lambda + \frac{1}{C_0 L_0} = 0,
\end{equation}
whose discriminant is 
\begin{equation}\label{Delta}
\Delta = \left( \frac{R_0}{L_0} \right)^2 - \frac{4}{C_0 L_0}.
\end{equation}
At this point, we will show that, regardless of the sign of the above discriminant $\Delta$, the eigenvalues associated to $\bm{A}$ have always strictly negative real part, condition that ensures the asymptotic stability of the homogeneous part of system (\ref{matrixForm}). We distinguish and analyze the following three cases:
\begin{itemize}
	\item[1.] If $\Delta < 0$, then the eigenvalues $\lambda_{1,2}$ of matrix $\bm{A}$ are complex and conjugate, given by
	\begin{equation}\label{eigenvalues1}
	\lambda_{1,2} = - \frac{R_0}{2 L_0} \pm i \frac{1}{2} \sqrt{-\Delta} =  - \frac{R_0}{2 L_0} \pm i \frac{1}{2} \sqrt{\frac{4}{C_0 L_0} - \left( \frac{R_0}{L_0} \right)^2},
	\end{equation}
	with strictly negative real part, that is
	\begin{equation}
	Re(\lambda_1) = Re(\lambda_2) = - \frac{R_0}{2 L_0} <0.
	\end{equation}
	\item[2.] If $\Delta = 0$, then the two eigenvalues associated to matrix $\bm{A}$ are equal, real and strictly negative, namely
	with strictly negative real part, that is
	\begin{equation}\label{eigenvalues2}
	\lambda_1 = \lambda_2 = - \frac{R_0}{2 L_0} <0.
	\end{equation}
	\item[3.] If $\Delta > 0$, then the eigenvalues $\lambda_{1,2}$ of matrix $\bm{A}$ are distinct, real and both negative, given by 
	\begin{equation}\label{eigenvalues3}
	\lambda_1 = - \frac{1}{2} \left( \frac{R_0}{L_0} + \sqrt{\Delta} \right), \quad \lambda_2 = - \frac{1}{2} \left( \frac{R_0}{L_0} - \sqrt{\Delta} \right).
	\end{equation}
	From (\ref{eigenvalues3}), the first eigenvalue $\lambda_1$ is clearly strictly negative, and it is straightforward to verify that the second eigenvalue $\lambda_2$ is also strictly negative. Indeed, the following chain of inequalities holds
	\begin{equation}
	0 < \Delta = \left( \frac{R_0}{L_0} \right)^2 - \frac{4}{C_0 L_0} < \left( \frac{R_0}{L_0} \right)^2 \quad \Longrightarrow \quad  \sqrt{\Delta} = \sqrt{\left( \frac{R_0}{L_0} \right)^2 - \frac{4}{C_0 L_0}} < \frac{R_0}{L_0},
	\end{equation}
	which implies $\lambda_2 < 0$.
\end{itemize}
In conclusion, the eigenvalues associated to the constant coefficient matrix $\bm{A}$ of system (\ref{matrixForm}) have always strictly negative real part. Therefore, the homogeneous part of system (\ref{matrixForm}) is asymptotically stable and, as a consequence, the complete ODE system (\ref{PinQoutLin1}) is stable, meaning that, for any choice of the initial condition, the exact solution will converge to the periodic solution.\\

Under the assumption that the forcing terms prescribed at the inlet and outlet of the vessel, $Q_{in}\equiv Q_{in}(t)$ and $P_{out}\equiv P_{out}(t)$, respectively, are periodic functions of $t$, the stability analysis of the $(Q_{in}, P_{out})$-type vessel is similar to that of the $(P_{in}, Q_{out})$-type vessel performed above. Indeed, it is straightforward to check that the eigenvalues of the coefficient matrix $\bm{A}$ associated to the linear ODE system corresponding to the $(Q_{in}, P_{out})$-type 0D vessel are the same to those of the $(P_{in}, Q_{out})$-type 0D vessel, thus leading to the same stability properties of the exact solution of the ODE system.\\

Next, we move to the $(P_{in}, P_{out})$-type 0D vessel displayed in the third row of Table \ref{table:equations0Dmodels}. The governing ODE system, of which we want to investigate the stability, reads
\begin{equation}\label{PinPoutLin}
\left\{ \begin{aligned}
& \frac{d V}{dt} = Q - Q_d,\\
& \frac{d Q}{dt} = \frac{1}{L_0/2} \left[ P_{in} - \frac{R_0}{2} Q - P \right],\\
&  \frac{d Q_d}{dt} = \frac{1}{L_0/2} \left[ P - \frac{R_0}{2} Q_d - P_{out} \right],
\end{aligned}
\right.
\end{equation}
where, as usual, the pressure data prescribed at the inlet and outlet of the vessel, $P_{in} \equiv P_{in}(t)$ and $P_{out} \equiv P_{out}(t)$, respectively, are assumed to be time-dependent periodic functions of a certain period $T_0 > 0$. In this case, the eigenvalues associated to the constant coefficient matrix $\bm{A}$ of the homogeneous part of system (\ref{PinPoutLin}) are given by
\begin{equation}
\lambda_1 = -\frac{R_0}{L_0}, \quad \lambda_{2,3} = \begin{cases}
- \dfrac{R_0}{2 L_0} \pm i \dfrac{1}{2} \sqrt{-\Delta}, & \text{ if } \Delta = \left( \dfrac{R_0}{L_0} \right)^2 - \dfrac{16}{C_0 L_0} <0,\\[2ex]
- \dfrac{R_0}{2 L_0},  & \text{ if } \Delta = 0,\\[2ex]
- \dfrac{R_0}{2 L_0} \pm \dfrac{1}{2} \sqrt{\Delta}, & \text{ if } \Delta > 0.
\end{cases}
\end{equation}
Therefore, we conclude that all the above eigenvalues have always strictly negative real part, meaning that the homogeneous part of system (\ref{PinPoutLin}) is asymptotically stable. As a consequence, for any choice of the initial conditions, any solution of the inhomogeneous system (\ref{PinPoutLin}) will converge to the periodic one.\\

We consider now the linear ODE system governing the last 0D vessel configuration, the $(Q_{in}, Q_{out})$-type 0D vessel depicted in the bottom row of Table \ref{table:equations0Dmodels}, that is
\begin{equation}\label{QinQoutLin}
\left\{ \begin{aligned}
& \frac{d V}{dt} = Q_{in} - Q,\\
& \frac{d Q}{dt} = \frac{1}{L_0} \left[ P - R_0 Q - P_d \right],\\
& \frac{d V_d}{dt} = Q - Q_{out},
\end{aligned}
\right.
\end{equation}
where $R_0$ denotes the resistance element between $P$ and $P_d$ only. Algebraic manipulations yield the following eigenvalues associated to the constant coefficient matrix $\bm{A}$ of the homogeneous part of system (\ref{QinQoutLin}) 
\begin{equation}\label{eigenvalQinQout}
\lambda_1 = 0, \quad \lambda_{2,3} = \begin{cases}
- \dfrac{R_0}{2 L_0} \pm i \dfrac{1}{2} \sqrt{-\Delta}, & \text{ if } \Delta = \left( \dfrac{R_0}{L_0} \right)^2 - \dfrac{16}{C_0 L_0} <0,\\[2ex]
- \dfrac{R_0}{2 L_0}, & \text{ if } \Delta = 0,\\[2ex]
- \dfrac{R_0}{2 L_0} \pm \dfrac{1}{2} \sqrt{\Delta}, & \text{ if } \Delta > 0.
\end{cases}
\end{equation}
The eigenvalues $\lambda_{2,3}$ have always strictly negative real part regardless the sign of the discriminant $\Delta$, while the first eigenvalues $\lambda_1$ turns out to be equal to zero.\\
In the following, we are going to study the asymptotic properties of system (\ref{QinQoutLin}) in order to find suitable assumptions on the periodic forcing functions $Q_{in}(t)$ and $Q_{out}(t)$ ensuring the stability of such ODE system, for both cases $\Delta<0$ and $\Delta>0$. However, by analyzing the orders of magnitude of typical physiological values of all geometrical and physical parameters defining the elements $R_0$, $L_0$ and $C_0$, and thus the expression of $\Delta$ in (\ref{eigenvalQinQout}), it is easy to check that we always have $\Delta<0$. The above expression of $\Delta$ can be reformulated as follows 
\begin{equation}\label{deltaQinQout}
\begin{split}
\Delta & = \left( \dfrac{R_0}{L_0} \right)^2 - \dfrac{16}{C_0 L_0} = \dfrac{1}{L_0} \left( \dfrac{R_0^2}{L_0} - \dfrac{16}{C_0} \right)\\
& = \dfrac{1}{L_0}\dfrac{4}{\pi r_0^3} \left[ f_1 - f_2 \right],
\end{split}
\end{equation}
with
\begin{equation}\label{f1f2}
f_1 = (\zeta+2)^2 \dfrac{\mu^2 l }{\rho r_0^3}, \quad f_2 = \dfrac{8}{3} \dfrac{E h_0}{l}.
\end{equation}
We consider the following order of magnitude ranges in arterial vessels for the variables defining the above factors $f_1$ and $f_2$
\begin{equation}
\begin{cases}
\rho \sim 10^0 & \text{[g/cm$^3$]}, \\
\mu \sim 10^{-2} & \text{[dyne$\cdot$ s/cm$^2$]}, \\
l \sim 10^0-10^1 & \text{[cm]}, \\
r_0 \sim 10^{-1}-10^0 & \text{[cm]}, \\
 h_0 \sim 10^{-2}-10^{-1} & \text{[cm]}, \\
 E \sim 10^6-10^7 & \text{[dyne/cm$^2$]}.
\end{cases}
\end{equation}
Then, for the first term $f_1$, its order of magnitude approximately ranges between $10^{-3}$ and $10^2$, while the order of magnitude of the second factor $f_2$ is estimated to vary between $10^4$ and $10^6$. Therefore, the second term $f_2$ is always the largest one, thus ensuring the negativity of $\Delta$. These findings were also confirmed by computing exact values of $f_1$, $f_2$ and $\Delta$ for all vessels of all arterial networks considered in this paper and described in Section \ref{sec:numexp}. Table \ref{table:deltaQinQout} displays maximum and minimum values, mean value and corresponding standard deviation of the two factors $f_1$ and $f_2$ defined in (\ref{f1f2}), for the aortic bifurcation model (Section \ref{sec:numexp:test1}), the 37-artery network (Section \ref{sec:numexp:test2}) and the reduced ADAN56 model (Section \ref{sec:numexp:test3}). Even if there is high variability in the values of these two factors, we observe that the second factor $f_2$ is always orders of magnitude larger with respect to the first factor $f_1$, thus implying that the discriminant $\Delta$ corresponding to the $(Q_{in}, Q_{out})$-type 0D vessel and given in (\ref{deltaQinQout}) is always strictly negative for all vessels of the three arterial networks considered. Then, our case of interest is the one corresponding to $\Delta<0$.\\
\begin{table}[h!]\footnotesize
	\centering
	\renewcommand\arraystretch{1.2}
	\begin{tabular}{l|cccc|cccc}
		\toprule
		\textbf{Network} & max($f_1$) & min($f_1$) & mean($f_1$) & std.dev($f_1$) & max($f_2$) & min($f_2$) & mean($f_2$) & std.dev($f_2$)\\
		\midrule
		\textbf{Aortic bifurcation} & 7.19 & 2.47 & 5.61 & 2.72 & 160000.00 & 158117.65 & 158745.10 & 1086.78 \\
		\textbf{37-artery network} & 315.26 & 7.30e-02 & 82.72 & 97.86 & 3076923.08 & 17777.78 & 227502.13 & 563787.17 \\
		\textbf{ADAN56 model} & 1837.19 & 8.96e-03 & 113.44 & 313.31 & 1777910.46 & 5184.12 & 187403.52 & 354118.73 \\
		\bottomrule
	\end{tabular}
	\caption{Quantitative assessment of the two factors $f_1$ and $f_2$ defining the discriminant $\Delta$ of the characteristic polynomial associated to the coefficient matrix $\bm{A}$ of the $(Q_{in}, Q_{out})$-type 0D vessel, for the vessels of the aortic bifurcation, the 37-artery network and ADAN56 model.}\label{table:deltaQinQout}
\end{table}

In general, if a linear homogeneous system of ODEs with constant coefficients has a null eigenvalue, then the coefficient matrix $\bm{A}$ is singular, with non-trivial null space, and any vector of the null space is an equilibrium point for the system. In other words, the homogeneous system (\ref{homogeneous}) does not have a unique equilibrium point, but a line of equilibria, which can be either stable (but not asymptotically stable) or unstable, depending on the sign of the other eigenvalues. Hence, an homogeneous system of the form (\ref{homogeneous}) with coefficient matrix $\bm{A}$ having a zero eigenvalue is stable if all other eigenvalues of  $\bm{A}$ have strictly negative real part, in the sense that it has an attractive line of equilibria and each equilibrium is stable, but not asymptotically stable.\\
However, the stability of the homogeneous system is not sufficient to ensure the stability of the corresponding inhomogeneous system (\ref{matrixForm}), but an additional assumption on the periodic forcing function is needed in order to preserve the stability of ODE system, namely that any solution of  (\ref{matrixForm}) for any admissible choice of the initial condition will converge to the periodic one as $t \to +\infty$.

We state here only the obtained condition, but the full derivation of this assumption on the periodic forcing function $\bm{b}(t)$ is extensively provided in Appendix \ref{appendix:stabilityanalysis:part2}, first in the scalar case of a single ODE, then for a system of ODEs, specifically focusing on system (\ref{QinQoutLin}) governing the $(Q_{in}, Q_{out})$-type 0D vessel.\\
Given the inhomogeneous linear ODE system (\ref{QinQoutLin}), whose coefficient matrix $\bm{A}$ has a null eigenvalue $\lambda_1$ and two eigenvalues $\lambda_{2,3}$ with strictly negative real part, the following condition on the periodic forcing function $\bm{b}(t)$ ensures the stability of the exact solution of the ODE system
\begin{equation}\label{condtionQinQoutFINAL}
\int_{0}^{T_0}  \left[ Q_{in}(s) - Q_{out} (s) \right]\, ds = 0
\quad \Longleftrightarrow \quad 
\int_{0}^{T_0}  Q_{in}(s)\, ds = \int_{0}^{T_0}  Q_{out} (s)\, ds.
\end{equation}
Namely, under this assumption, the complete inhomogeneous system (\ref{QinQout}) is stable, in the sense that for any admissible choice of the initial condition the exact solution will converge to the periodic one. This is also a physically consistent condition: in order for the volume $V(t)+V_d(t)$ not to constantly increase/decrease and asymptotically explode, the integral over a period $[0, T_0]$ of the inflow $Q_{in}(t)$ entering the vessel must equal the integral over 
the same period of the outflow $Q_{out}(t)$ leaving the vessel.

\begin{remark}
The present work focuses on arteries. However, we expect that the family of nonlinear 0D models for blood flow derived in Section \ref{sec:0Dmodel} can be applied not only to arteries, but also to veins, by appropriately changing the geometrical and mechanical properties of vessels and the tube law relating the mean internal pressure to the vessel cross-sectional area. Indeed, in the tube law (\ref{tubelaw}), typical values for parameters $m$ and $n$ for collapsible tubes, such as veins, are $m\approx 10$ and $n=-1.5$. A relation for the venous stiffness $K_v$ can also be derived from considerations made on the collapse of thin-walled elastic tubes, or, alternatively, $K_v$ can be estimated from pulse wave velocities, as described in \cite{Mynard:2011a}. Furthermore, the stability analysis presented in this section for arterial vessels can be straightforwardly repeated also for veins, in order to study the stability properties of the corresponding ODE systems.
\end{remark}

\section{0D junctions and networks}\label{sec:0Djunction}

Equipped with the family of nonlinear 0D models for blood flow derived in Section \ref{sec:0Dmodel} for a single vessel, we consider now the coupling of 0D vessels to construct more complex networks of vessels.\\
Two or more 0D vessels can be coupled through 0D junctions, which satisfy the conservation of mass and also impose a common pressure on all branches, to ensure the continuity of pressure throughout the 0D junction. We immediately note that, in contrast to 1D junctions between 1D vessels, where the total pressure continuity can be enforced, in junction between 0D vessels we enforce pressure continuity only. The choice of imposing pressure continuity at each 0D junction allows us to always solve a linear coupling problem, whereas it is well-known that the nonlinear problem to be solved at 1D junctions can become very computationally expensive. However, on the other hand, in order to arrange compatible segment types into a network, with inlets and outlets coupled appropriately, restrictions on admissible 0D vessel types are necessary for vessels converging at a 0D junction. In particular, input data to prescribe at the inlet and outlet of each vessel are defined by the state of their adjacent compartments, in order to ensure the conservation of mass and continuity of pressure. The coupling of 0D vessels without any restrictions on compatible 0D vessel configurations would be possible by enforcing total pressure continuity coming at the cost of solving nonlinear coupling problem at each junction node.\\
The simplest 0D junctions are two-vessel and three-vessel junctions, connecting two or three vessels, respectively, which can be generalized in a generic 0D junction attaching an arbitrary number of 0D vessels. In addition, in open-loop networks, terminal vessels can be coupled to single-resistance or $RCR$ Windkessel elements to model the cumulative effects of all vessels distal to the terminal segments of the vessel network.

\subsection{Two-vessel junction (J2)}\label{sec:0Djunction:J2}

\begin{figure}[h!]
	\begin{center}
		\subfloat[][J2\label{fig:J2}]{\includegraphics[scale=0.8]{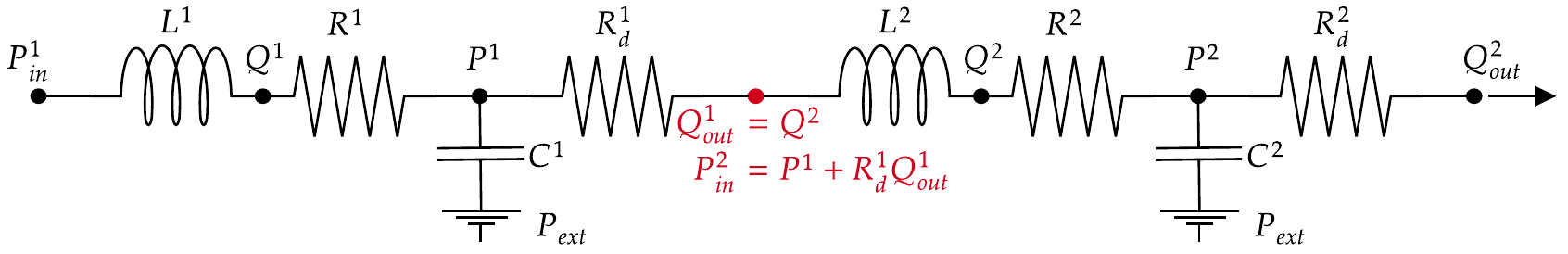}}\\
		\subfloat[][J3\label{fig:J3}]{\includegraphics[scale=0.7]{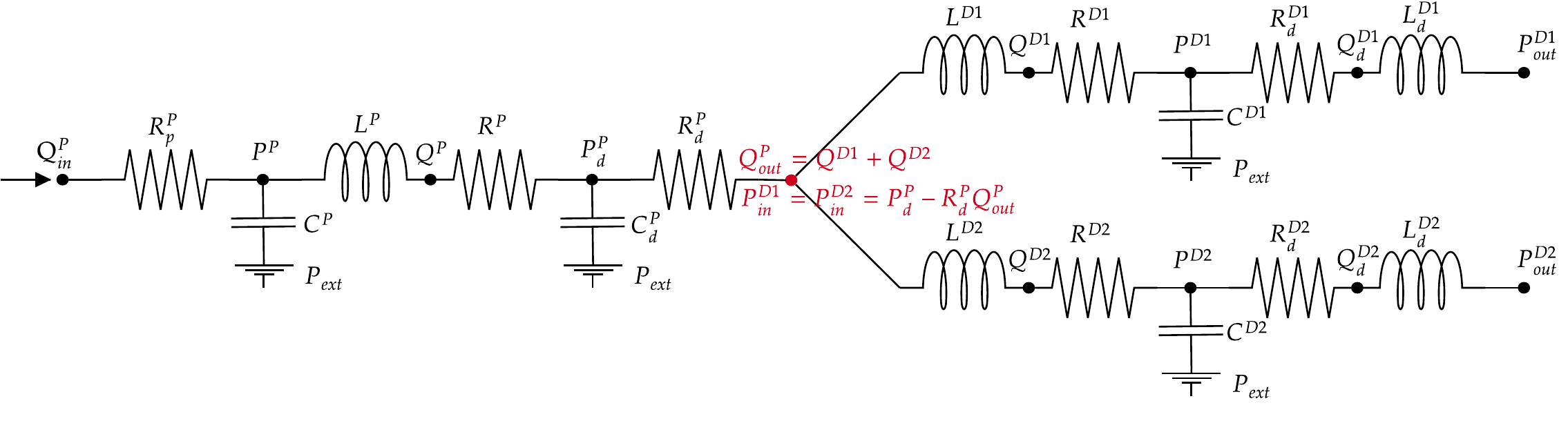}}\\
		\caption{(A) Example of a 0D representation for a two-vessel junction between two $(P_{in}, Q_{out})$-type 0D vessels. (B) Example of a 0D representation for a bifurcation branch (three-vessel splitting flow junction): $(Q_{in}, Q_{out})$-type 0D vessel for parent vessel and $(P_{in}, P_{out})$-type 0D vessels for daughter vessels.}
	\end{center}
\end{figure}

To describe the coupling procedure adopted, let us consider, for instance, two $(P_{in}, Q_{out})$-type 0D vessels connected in a simple two-vessel junction, as displayed in Figure \ref{fig:J2}. Coupling conditions are needed to determine, at the junction point, the flow rate $Q_{out}^1$ at the outlet of the first (parent) vessel and the pressure $P_{in}^2$ at the inlet of the second (daughter) vessel. By enforcing the conservation of mass, we obtain:
\begin{equation*}
Q_{out}^1 (t) = Q^2(t),
\end{equation*}
where $Q^2$ is the flow rate state variable in the second vessel; then, by imposing the continuity of pressure, we get:
\begin{equation*}
P_{in}^2(t) = P^1(t) - R_d^1 Q_{in}^1(t) =  P^1(t) + R_d^1 Q^2(t),
\end{equation*}
where $P^1$ is the pressure state variable in the parent vessel and $R_d^1$ is the distal resistance at the outlet of the same vessel, as illustrated in Figure \ref{fig:J2}. Note that in the case where $R_d^1=0$, we simply obtain $P_{in}^2 = P^1$.\\
From this test case, it is straightforward to conclude that all the possible pairs of vessel types that can be coupled to form a two-vessel 0D junction are such that the outlet of the parent vessel is of pressure type and the inlet of the daughter vessel is of flow type, or vice versa. All other configurations are not allowed, because no assumption would be made either on the flow rate $Q$ or on the pressure $P$ at the interface between vessels.\\
Moreover, we observe that the two-vessel junction may be also used to represent a single 0D vessel with two 0D compartments of the same type coupled in series, either $(P_{in}, Q_{out})$ or $(Q_{in}, P_{out})$. Indeed, the per-segment mapping replaces each 1D vessel of a network by a 0D vessel, which, in turn, can be composed of just one or more 0D compartments.\\

\subsection{Three-vessel junction (J3)}\label{sec:0Djunction:J3}

For the family of three-vessel junctions, we present here the coupling procedure adopted in the case of a \textit{splitting flow} junction, where the 0D junction represents the branching point at which the end of the parent vessel is connected to the inlets of the two daughter vessels, but the coupling conditions to be imposed in a \textit{merging flow} junctions, where the 0D junction represents the adjoining point at which the outlets of the two daughter vessels converge into the beginning of the parent vessel, can be easily derived in a similar way, as proposed in \cite{Safaei:2018a}.\\
We consider, for instance, a $(Q_{in}, Q_{out})$-type 0D vessel for the parent vessel and $(P_{in}, P_{out})$-type 0D vessels for both daughter vessels. By imposing the mass conservation, we get that the flow rate $Q_{out}^{P}$ at the outlet of the parent vessel must be equal to the sum of the two daughter branches' flows $Q^{D1}$ and $Q^{D2}$, that is
\begin{equation*}
Q_{out}^P(t) = Q^{D1}(t) + Q^{D2}(t).
\end{equation*}
Then, by enforcing the continuity of pressure, we have that the pressure at the inlet of both daughter vessels must be equal to the distal pressure in the parent vessel, namely
\begin{equation*}
P_{in}^{D1}(t) = P_{in}^{D2}(t) = P_d^P(t) - R_d^P Q_{out}^P(t) =  P_d^P(t) - R_d^P \left( Q^{D1}(t) + Q^{D2}(t) \right),
\end{equation*}
where the pressure $P_d^P$ is the distal pressure state variable in the parent vessel and $R_d^P$ is the distal resistance at the outlet of the same vessel, as depicted in Figure \ref{fig:J3}. We observe that in the particular case where $R_d^P=0$, the above condition of continuity of pressure  becomes $P_{in}^{D1} = P_{in}^{D2} = P_d^P$.\\
From this test case, we can then derive the following restrictions on compatible segment types for vessels converging in a three-vessel splitting flow junction: the outlet of the parent vessel must be of flow type, while the inlets of both daughter vessels must be of pressure type.

\subsection{Generic 0D junction} \label{sec:0Djunction:generic}

\begin{figure}[h!]
	\begin{center}
		\includegraphics[scale=0.9]{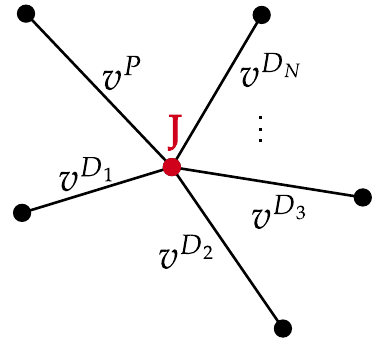}	
		\caption{Example of a generic 0D junction connecting $(N+1)$ 0D vessels: one vessel is chosen as parent vessel, while all other N segments are treated as daughter vessels.}\label{fig:genericJ}
	\end{center}
\end{figure}

In the most general situation, a 0D junction connects an arbitrary number of 0D vessels, sharing their inlets or outlets at the junction point, as displayed in Figure \ref{fig:genericJ}.\\
In order to couple appropriately all the 0D vessels converging at the junction, first of all, a single vessel has to be chosen as parent vessel (for instance, the vessel with the largest cross-sectional area), while all other vessels are classified as daughter vessels. As a consequence, the role of each vessel, either parent or daughter, will define the corresponding vessel type, depending on the inlet/outlet data to be prescribed at the junction point.\\
For the parent segment, at the vessel end shared at the junction a condition on the flow $Q^P$ is prescribed, which is computed by imposing the conservation of mass. Indeed, at the junction, the flow rate in each daughter vessel is known and denoted by $Q^{D_j}$. For instance, if we refer to the configuration illustrated in Figure \ref{fig:genericJ}, the flow rate $Q_{out}^P$ at the outlet of the chosen parent vessel is given by
\begin{equation}\label{consFlow}
Q_{out}^P (t) = \sum_{j=1}^{N} Q^{D_j} (t).
\end{equation}
In equation (\ref{consFlow}), the flow direction in each vessel, namely if the blood stream enters or leaves the 0D junction, is taken into account by the sign of $Q$. On the other hand, for each daughter, at the vessel end shared at the junction a condition on the pressure $P^{D_j}$ must be prescribed in order to enforce pressure continuity throughout the 0D junction. Then, this pressure must be equal to the pressure in the parent vessel, that is
\begin{equation}
P_{in/out}^{D_j} (t) = P^P (t) \quad j = 1,\ldots, N.
\end{equation}

\subsection{Terminal vessels}\label{sec:0Djunction:terminal}

\begin{figure}[h!]
	\begin{center}
		\includegraphics[scale=0.8]{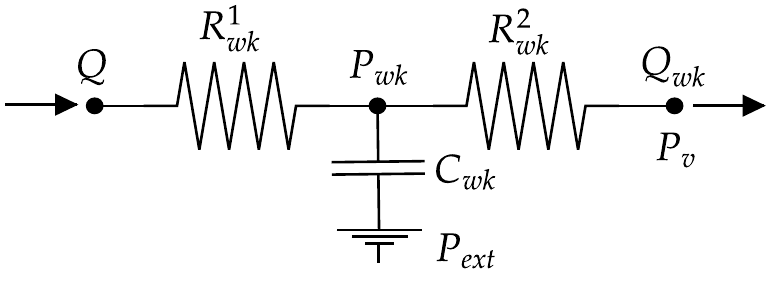}	
		\caption{$RCR$ Windkessel terminal element.}\label{fig:RCR}
	\end{center}
\end{figure}

In open-loop arterial networks, the cumulative effects of all distal vessels (small arteries, arterioles and capillaries) at ending locations of terminal arteries have to be taken into account. These effects can be modelled using either single-resistance or $RCR$ Windkessel elements coupled to the terminal arteries.\\
Each $RCR$ Windkessel element is composed of a proximal terminal resistance $R_{wk}^1$, a terminal capacitor $C_{wk}$ and a distal terminal resistance $R_{wk}^2$, as displayed in Figure \ref{fig:RCR}. Pressure $P_{wk}$ and flow rate $Q_{wk}$ in this terminal element are governed by
\begin{equation}
\left\{ \begin{aligned}
& \frac{d P_{wk}}{dt} = \frac{Q - Q_{wk}}{C_{wk}},\\
& Q_{wk} = \frac{P_{wk} - P_v}{R_{wk}^2},
\end{aligned}
\right.
\end{equation}
where $Q$ is the flow rate from the 0D vessel coupled to the $RCR$ element and $P_v$ is the constant outflow pressure.\\
Depending on the configuration of the 0D vessel coupled to the Windkessel model, either the flow rate $Q_{out}$ or the pressure $P_{out}$ must be prescribed at the outlet of the vessel. For instance, in the case of a $(P_{in}, Q_{out})$-type 0D vessel, the flow rate $Q_{out}$ has to be assigned at the outlet of the vascular segment, as illustrated in the top row of Table \ref{table:equations0Dmodels}, which is computed from the coupling to the Windkessel element, as follows
\begin{equation}
	Q_{out} = \frac{P - P_{wk}}{(R_d + R_{wk}^1)},
\end{equation}
where $P$ and $R_d$ are pressure and (if present) distal resistance in the 0D vessel, respectively. On the other hand, if we consider, for example, a $(Q_{in}, P_{out})$-type 0D vessel, the pressure $P_{out}$ to be enforced at the outlet of the vessel, as displayed in the second row of Table \ref{table:equations0Dmodels}, can be calculated as
\begin{equation}
	P_{out} = P_{wk} + R_{wk}^1 Q,
\end{equation}
where $Q$ is the flow rate in the 0D vessel coupled to the $RCR$ element.\\
Alternatively, if terminal vessels are coupled to single-resistance terminal elements, we simply get
\begin{equation}
Q_{out} = \frac{P - P_v}{(R_d + R_{wk})} \qquad \text{or} \qquad P_{out} = P_v + R_{wk} Q,
\end{equation}
where $R_{wk}$ is now the only peripheral resistance to the flow in the terminal element.

\section{Numerical experiments}\label{sec:numexp}

In this section, we validate the ability of the derived nonlinear 0D models to reproduce the essential blood flow distribution and the main features of pressure and flow waveforms in networks of deformable vessels with respect to the well-known and widely used 1D blood flow model (\ref{sys1D}), at a much lower computational cost.\\
We reproduce several benchmark problems proposed in \cite{Boileau:2015a}. First, we consider a simple test case: a model of blood flow in the aortic bifurcation (Section \ref{sec:numexp:test1}). Then, we assess the results in two different arterial networks: the 37-artery network of the aorta and the largest central system arteries constructed in \cite{Matthys:2007a}, for which \textit{in vitro} pressure and flow waveforms were acquired in \cite{Matthys:2007a} (Section \ref{sec:numexp:test2}), and a reduced version of the ADAN model developed by Blanco  \textit{et al.} \cite{Blanco:2014a,Blanco:2015a}, which contains the largest 56 systemic arteries of the human circulation (Section \ref{sec:numexp:test3}), referred to as ADAN56. For each test case, 0D results of pressure and flow rate in each vessel of the network are qualitatively and quantitatively compared to the 1D results obtained from the 1D model (\ref{sys1D}). Furthermore, in order to assess the properties of the newly derived nonlinear 0D models, we also compare the results from the nonlinear and the linear 0D models, where in the latter the convective terms are neglected \textit{a-priori}, the model parameters are constant and the pressure-volume relation is linear. For each test case, we provide graphical comparisons supported by error tables.\\
The numerical methods adopted to solve the 1D and 0D models are described in Section \ref{sec:numexp:methods}. The relative error metrics for pressure $P$ and flow rate $Q$ are introduced in Section \ref{sec:numexp:errors}. Before performing the 0D simulations and compare 0D results to 1D results, a detailed analysis of the contribution and relative importance of the convective component within the momentum equation is carried out in Section \ref{sec:numexp:convterms}, to decide whether it is reasonable or not to neglect the convective terms in the 0D models. 

\subsection{Numerical methods}\label{sec:numexp:methods}

\subsubsection{Second-order MUSCL-Hancock scheme for the 1D model}

System (\ref{sys1D}) governing 1D blood flow is solved using a second-order MUSCL-Hancock numerical scheme \cite{vanLeer:1984b}, where MUSCL stands for Monotonic Upstream-Centred Scheme for Conservation Laws, with ENO (Essentially Non-Oscillatory) reconstruction \cite{Harten:1987a,Harten:1987b} and numerical source computed following the ADER approach \cite{Toro:2001c,Toro:2009}.\\
In order to approximate the solutions of system (\ref{sys1D}), we first discretize the 1D space domain $\left[ 0, l \right]$ in $M$ cells of constant size $\Delta x$, where each cell $I_i = ( x_{i-\frac{1}{2}} , x_{i+\frac{1}{2}} )$ has centre located in $x_i$, with $x_{i\pm\frac{1}{2}} = x_i \pm \frac{\Delta x}{2}$, for $i = 1, \ldots, M$. The time domain $\left[ 0, t_{\text{end}} \right]$ is also discretized by assuming a constant time step $\Delta t$, which is restricted according to the usual Courant-Friedrichs-Lewy (CFL) stability condition. Then, if we consider system (\ref{sys1D}) written under the balance law form (\ref{conservationlaws}), given the approximate solution $\bm{Q}_i^n$ at time $t_n$ in each cell $I_i$, we can evolve the numerical solution to time $t_{n+1}=t_n + \Delta t$ by using a finite volume method of the form
\begin{equation}\label{FVS}
\bm{Q}_i^{n+1} = \bm{Q}_i^n - \frac{\Delta t}{\Delta x} \left[ \bm{F}_{i+\frac{1}{2}} - \bm{F}_{i-\frac{1}{2}} \right] + \Delta t \bm{S}_i,
\end{equation}
where $\bm{F}_{i+\frac{1}{2}}$ is the numerical flux that approximates the time-integral average over $[t_n, t_{n+1} ]$ of flux $\bm{F}$ at the cell interface $x=x_{i+\frac{1}{2}}$, while $\bm{S}_i$ is the numerical source in cell $I_i$ that approximates the volume-integral average over $V=[ x_{i-\frac{1}{2}} , x_{i+\frac{1}{2}} ] \times [t_n, t_{n+1} ]$ of the source term $\bm{S}$. Finite volume schemes (\ref{FVS}) may be interpreted as resulting from integrating the equations of system (\ref{conservationlaws}) on the control volume $V=[ x_{i-\frac{1}{2}} , x_{i+\frac{1}{2}} ] \times [t_n, t_{n+1} ]$, where suitable approximations of the integral averages have been introduced.\\
In this framework, the MUSCL–Hancock approach achieves a second-order extension of the well-known Godunov's first-order upwind method by computing the intercell flux $\bm{F}_{i+\frac{1}{2}}$ according to the following three steps:
\begin{itemize}
	\item[(I)] Data reconstruction and cell boundary values
	\begin{equation}
	\bm{Q}_i^L = \bm{P}_i(x_{i-\frac{1}{2}}), \quad \bm{Q}_i^R = \bm{P}_i(x_{i+\frac{1}{2}}),
	\end{equation}
	where $\bm{P}_i(x)$ is the first-degree reconstruction polynomial vector in cell $I_i = ( x_{i-\frac{1}{2}} , x_{i+\frac{1}{2}} )$, that is
	\begin{equation}\label{reconstructionpol}
	\bm{P}_i(x) = \bm{Q}_i + ( x- x_i) \Delta_i,
	\end{equation}
	with $\Delta_i$ being the slope vector associated to the reconstruction polynomial (\ref{reconstructionpol}), here computed by using the ENO criterion to preserve conservation and non-oscillatory properties.
	\item[(II)] Evolution of boundary extrapolated values by a time $\frac{1}{2}\Delta t$ accounting for source term
	\begin{equation}
	\begin{cases}
	\overline{\bm{Q}}_i^L = \bm{Q}_i^L - \dfrac{1}{2}\dfrac{\Delta t}{\Delta x} \left[ \bm{F}(\bm{Q}_i^R) - \bm{F}(\bm{Q}_i^L) \right] + \dfrac{1}{2} \Delta t \bm{S}(\bm{Q}_i^L),\\[2ex]
	\overline{\bm{Q}}_i^R = \bm{Q}_i^R - \dfrac{1}{2}\dfrac{\Delta t}{\Delta x} \left[ \bm{F}(\bm{Q}_i^R) - \bm{F}(\bm{Q}_i^L) \right] + \dfrac{1}{2} \Delta t \bm{S}(\bm{Q}_i^R).
	\end{cases}
	\end{equation}
	\item[(III)] Solution of a classical Riemann problem with data $\left( \overline{\bm{Q}}_i^R, \overline{\bm{Q}}_{i+1}^L \right)$ to obtain the similarity solution $\bm{Q}_{i+\frac{1}{2}}(x/t)$ to compute the intercell flux
	\begin{equation}
	\bm{F}_{i+\frac{1}{2}} = \bm{F}\left( \bm{Q}_{i+\frac{1}{2}}(0) \right). 
	\end{equation}
\end{itemize}
As last step, the numerical source is computed imitating the ADER approach \cite{Toro:2001c}, as follows
\begin{equation}
\bm{S}_i = \bm{S}\left( \bm{Q}_i^n + \frac{1}{2}\Delta t \left[ -\bm{A}(\bm{Q}_i^n) \Delta_i + \bm{S}(\bm{Q}_i^n) \right] \right),
\end{equation}
where $\bm{A}$ is the Jacobian matrix of system (\ref{conservationlaws}).\\

The coupling of several 1D vessels at junction points is treated following the methodology proposed in \cite{Sherwin:2003a} and extended in \cite{Spilimbergo:2021a} to achieve second-order accuracy also of the coupling procedure and preserve the global second-order accuracy in space and time over the entire 1D network. As in these cited papers, also here we will restrict to sub-critical flow conditions, \textit{i.e.} when $|u|<c$, which is a crucial assumption in ensuring the strictly hyperbolic nature of the PDE system and in determining the type of boundary conditions that can be applied to the 1D model.\\
In the case of $N_J$ vessels converging at a junction $J$ (for the networks considered in the present work, we will have $N_J=2$ for the two-vessel junction and $N_J=3$ for the three-vessel junction, only), the computational cells involved in the coupling of the $k$-th vessel, with $k = 1,\ldots, N_J$, provide the state $\bm{Q}_k^n = \left[ A_k^n, Q_k^n \right]^T$ at time $t^n$. Then, in  order to couple the $N_J$ vessels, we have to compute the unknown cross-sectional area $A_k^*$ and flow $Q_k^*$ for each vessel converging at node $J$, by imposing: (i) conservation of mass, (ii) continuity of total pressure; (iii) continuity of the generalized Riemann invariants. Therefore, to achieve the second-order coupling, we will set $\bm{Q}_k^n = \bm{Q}_k^{\text{MUSCL,}n}$, for $k = 1,\ldots, N_J$, where $\bm{Q}_k^{\text{MUSCL,}n}$ is the evolved boundary extrapolated value, given by either $\overline{\bm{Q}}_{k,1}^{L,n}$ if the first computational cell of vessel $k$ is converging to node $J$, or $\overline{\bm{Q}}_{k,M}^{R,n}$ if instead the last computational cell ($M$) of vessel $k$ converges to node $J$.\\
The same second-order reconstruction is also adopted for the coupling between terminal vessels and $RCR$ Windkessel/single-resistance terminal elements.\\

The number of computational cells used in the $j$-th vessel of each arterial network is defined as
\begin{equation}
M_j = \max\left\{ \ceil{ \frac{l_j}{\Delta x_{\text{max}}} }, 2 \right\},
\end{equation}
where $l_j$ is the length of vessel $j$. For all numerical experiments, before setting the maximum mesh size $\Delta x_{\text{max}}$, a mesh convergence study of the 1D solution was first carried out in order to select a reference sufficiently accurate 1D solution for the comparison with the 0D results. The values of $\Delta x_{\text{max}}$, which ensure a 1D mesh-independent solution, used in the different benchmark problems are displayed, together with all other computational parameters for the 1D/0D simulations, in Table \ref{table:param}.

\subsubsection{Numerical method for solving the 0D models}

In parallel with the fully 1D discretization considered so far, a vascular network can also be entirely modelled by using 0D vessels. In particular, as discussed in Section \ref{sec:0Djunction}, in order to arrange compatible segment types into a structure, the 0D configuration to be used for each vessel of the network has to be chosen so that inlets and outlets of vessels converging at the 0D junctions are all coupled appropriately. Hence, we end up with just one system of ODEs describing the dynamic of the entire network, where the single subsystems corresponding to each vessel are not isolated, but are connected to each other via the variables prescribed at the inlets and outlets of vessels. Indeed, pressures and/or flow rates imposed at the inlet and outlet of each vessel are defined by the state of the adjacent vessels, in order to ensure the conservation of mass and continuity of pressure.\\
The resulting ODE system can be written in compact form as follows:
\begin{equation}\label{globalODEsystem}
\frac{d \bm{y}(t)}{dt} = \bm{F} (t, \bm{y}(t)),
\end{equation}
where $\bm{y}(t)$ is the unknown vector containing the $2N$ state variables $(V_i(t), Q_i(t))$,  $i = 1, \ldots, N$, of the $N$ vessels of the network. Then, for all numerical tests, the global ODE system (\ref{globalODEsystem}) is solved using the four-step explicit fourth-order Runge-Kutta (RK4) method, with appropriate time step to guarantee the stability of the numerical scheme and the mesh independence of the solution (see Table \ref{table:param} for details).\\
It is worth also noting that numerical experiments have shown that the restrictions on the allowable time step ensuring the stability of RK4 method to solve the different arterial networks are the same when using linear and nonlinear 0D models. Indeed, as can be observed from Table \ref{table:param}, for each arterial network we adopt the same time step $\Delta t$ to solve both the linear and nonlinear fully 0D network configurations. Then, we conclude that the coupling between several nonlinear 0D models is not affecting the numerical stability of the discretized model.\\
\begin{table}[h!]\footnotesize
	\centering
	\begin{tabular}{lccc}
		\toprule
		\textbf{Parameter} & \textbf{Aortic bif.} & \textbf{37-artery network} & \textbf{ADAN56 model}\\
		\midrule
		$\Delta x_{\text{max}}$ & 2 mm & 1 mm & 1 mm \\
		CFL number & 0.9 & 0.9 & 0.9 \\
		RK4 time step, $\Delta t$ & 10$^{-3}$ s & 10$^{-4}$ s & 10$^{-4}$ s \\
		Cardiac cycle, $T_0$ & 1.1 s & 0.827 s & 1.0 s\\
		Final time, $t_{\text{end}}$ & 29.7 s & 24.81 s & 15.0 s \\
		\bottomrule
	\end{tabular}
	\caption{Computational parameters adopted in the 1D/0D numerical simulations for the three arterial networks.}\label{table:param}
\end{table}

\subsection{Error calculations}\label{sec:numexp:errors}

To provide a quantitative assessment of the predicted waveforms compared with the reference 1D solution and to measure the benefit, if any, that we get by preserving certain nonlinearities of the original 1D model in the newly derived 0D models, we introduce the following relative error metrics for pressure $P$ and flow rate $Q$:
\begin{equation}\label{metrics}
\begin{aligned}
\varepsilon_P^{\text{RMS}} = \sqrt{\dfrac{1}{n} \sum_{i=1}^{n} \left( \dfrac{P_i^{0D} - P_i^{1D}}{P_i^{1D}}\right)^2}, & \qquad \varepsilon_Q^{\text{RMS}} = \sqrt{\dfrac{1}{n} \sum_{i=1}^{n} \left( \dfrac{Q_i^{0D} - Q_i^{1D}}{\max_j \left( Q_j^{1D} \right)}\right)^2},\\
\varepsilon_P^{\text{SYS}} = \dfrac{\max \left( P^{0D} \right) - \max \left( P^{1D} \right)}{\max \left( P^{1D} \right)}, & \qquad \varepsilon_Q^{\text{SYS}} = \dfrac{\max \left( Q^{0D} \right) - \max \left( Q^{1D} \right)}{\max \left( Q^{1D} \right)},\\
\varepsilon_P^{\text{DIAS}} = \dfrac{\min \left( P^{0D} \right) - \min \left( P^{1D} \right)}{\min \left( P^{1D} \right)}, & \qquad \varepsilon_Q^{\text{DIAS}} = \dfrac{\min \left( Q^{0D} \right) - \min \left( Q^{1D} \right)}{\max \left( Q^{1D} \right)},
\end{aligned}
\end{equation}
where $i = 1, \ldots, n$ are time points over the cardiac cycle at which the solution is sampled, $P^{0D}$ and $Q^{0D}$ are 0D pressure and flow, either from the nonlinear or the linear 0D models, and $P^{1D}$ and $Q^{1D}$ are 1D pressure and flow at the midpoint of the vessel. We compare the solution obtained using 1D models sampled at this location since this is a commonly observed variable in this research field. Other choices are possible (like, for example, averaged quantities over the 1D domain) and would not affect the conclusions of this work (results not reported here).\\
All error metrics are calculated over a single cardiac cycle, once the numerical results are in the periodic regime. Periodicity is defined as the distance in $L^{\infty}$-norm between the normalized solutions over two consecutive cardiac cycles to be smaller than a threshold of $10^{-3}$ (pressure and cross-sectional area are normalized by the mean pressure and mean cross-sectional area, respectively, over the cardiac cycle; flow rate is normalized by the maximum flow over the cardiac cycle).

\subsection{Convective terms}\label{sec:numexp:convterms}

We note that we have made no assumptions yet about the contribution of the convective terms in the family of nonlinear 0D models governed by the system of ODEs (\ref{sys0D}). 
A first insight into the role of the convective term in the momentum balance equation of the 1D blood flow model (\ref{sys1D}) and its relative importance especially with respect to the pressure term is given by the dimensional analysis of the 1D equations carried out in Section \ref{sec:1Dmodel:dimanalysis} and, in the following, applied to the different arterial networks considered. The coefficients $\gamma_C$ and $\gamma_P$ characterizing the convective and pressure terms, respectively, in the nondimensional momentum balance equation (\ref{nondimmomentumeq}) are defined in terms of the average flow velocity $U_0$. Since for all the arterial networks of interest 1D simulations of blood flow have been performed to obtain 1D mesh-independent solutions, the average flow velocity $U_{1D}^{\text{mean}}$ and the maximum flow velocity $U_{1D}^{\text{max}}$ can be computed from the 1D results for each vessel of each network. Then, from the estimated velocities, we are able to quantify the nondimensional coefficients (\ref{coefficients:dimanalysis}) and to assess the contribution and relative importance of the convective, pressure and friction terms within the momentum equation. Values of the ratios $\gamma_C/\gamma_P$ and $\gamma_F/\gamma_P$ are displayed for the aortic bifurcation, some vessels of both the 37-artery network and ADAN56 model in Table \ref{table:dimanalysis}. Furthermore, for the 37-artery network maximum and mean values of the ratio $\gamma_C/\gamma_P$ are
\begin{equation}
\left( \gamma_C/\gamma_P \right)_{\text{max}}^{\text{net37}} = \begin{cases}
0.00259 & \text{if } U_0 = U_{1D}^{\text{mean}},\\
0.00947 & \text{if } U_0 = U_{1D}^{\text{max}},
\end{cases} 
\qquad
\left( \gamma_C/\gamma_P \right)_{\text{mean}}^{\text{net37}}  = \begin{cases}
4.801\text{e-04} & \text{if } U_0 = U_{1D}^{\text{mean}},\\ 
0.00223 & \text{if } U_0 = U_{1D}^{\text{max}},
\end{cases}
\end{equation} 
where the two maximum values of $\gamma_C/\gamma_P$ are found in the left anterior tibial and in the left iliac-femoral III arteries, respectively, while for the reduced ADAN56 model we have
\begin{equation}
\left( \gamma_C/\gamma_P \right)_{\text{max}}^{\text{adan56}}  = \begin{cases}
0.00791 & \text{if } U_0 = U_{1D}^{\text{mean}},\\
0.09804 & \text{if } U_0 = U_{1D}^{\text{max}},
\end{cases} 
\qquad
\left( \gamma_C/\gamma_P \right)_{\text{mean}}^{\text{adan56}}  = \begin{cases}
0.00195 & \text{if } U_0 = U_{1D}^{\text{mean}},\\
0.02233 & \text{if } U_0 = U_{1D}^{\text{max}}.
\end{cases}
\end{equation}  
where the two maximum values of $\gamma_C/\gamma_P$ are both achieved in the thoracic aorta VI.\\
The magnitude of these coefficient ratios clearly suggests that the pressure gradient is the dominating term in the momentum balance equation in (\ref{sys1D}), with respect to the convective and the friction terms. In particular, from Table \ref{table:dimanalysis} we observe that, on the one hand, the frictional losses become more and more important as we consider vessels of consecutive generations of bifurcation further from the aortic trunk, while, on the other hand, the pressure term is always significantly dominating over the convective component. As expected, the pressure gradient represents the main term in the momentum balance equation, while the contribution of the convective term turns out to be consistently smaller in all vessels of the three arterial networks.\\
In addition, this analysis shows that overall the ratio $\gamma_C/\gamma_P$ is larger in ADAN56 model than in the 37-artery network, suggesting that in ADAN56 model the contribution of the convective term is of greater importance.\\
\begin{table}[h!]\footnotesize
	\centering
	\renewcommand\arraystretch{1.2}
	\begin{tabular}{ll|cc|cc}
		\toprule
		\multicolumn{2}{c|}{} & \multicolumn{2}{c|}{$U_0 = U_{1D}^{\text{mean}}$} & \multicolumn{2}{c}{$U_0 = U_{1D}^{\text{max}}$} \\
		\midrule
		\textbf{Test case} & \textbf{Vessel name} & $\gamma_C/\gamma_P$ & $\gamma_F/\gamma_P$ & $\gamma_C/\gamma_P$ & $\gamma_F/\gamma_P$ \\
		\toprule
		\multirow{2}{*}{\textbf{Aortic bifurcation}} & Aorta & 2.537e-05 & 7.915e-05 & 0.00233 & 7.584e-04\\
		& Iliac artery & 2.027e-05 & 1.214e-04 & 9.479e-04 & 8.303e-04\\
		\midrule
		\multirow{8}{*}{\textbf{37-artery network}} & Aortic arch II & 1.047e-04 & 2.085e-05 & 0.00266 & 1.050e-04\\
		& Thoracic aorta II & 1.071e-04 & 7.530e-05 & 0.00273 & 3.802e-04 \\
		& L subclavian I & 2.199e-04 & 0.00174 & 0.00195 & 0.00518\\
		& R iliac-femoral II & 3.344e-04 & 0.00393 & 0.00327 & 0.01227\\
		& L ulnar & 5.963e-04 & 0.00734 & 0.00120 & 0.01041\\
		& R anterior tibial & 7.988e-04 & 0.01186 & 0.00200 & 0.01878\\
		& R ulnar & 6.389e-04 & 0.00603 & 0.00105 & 0.00775 \\
		& Splenic & 0.00121 & 0.01436 & 0.00223 & 0.01948 \\
		\midrule
		\multirow{12}{*}{\textbf{ADAN56 model}} & Aortic arch I & 0.00106 & 8.842e-05& 0.02354 & 4.168e-04\\
		& Thoracic aorta III & 0.00249 & 5.660e-05 & 0.03793 & 2.209e-04\\
		& Abdominal aorta V & 0.00228 & 3.677e-04 & 0.04167 & 0.00157\\
		& R common carotid & 7.868e-04 & 9.420e-04 & 0.01311 & 0.00384\\
		& R renal & 0.00453 & 0.00157 & 0.01482 & 0.00284 \\
		& R common iliac & 0.00224 & 0.00126 & 0.03577 & 0.00505 \\
		& R internal carotid &  0.00106 & 0.00324 &0.00994 & 0.00993\\
		& R radial & 0.00256 & 0.03949 & 0.01074 & 0.08010\\
		& R internal iliac & 0.00189 & 0.00230 & 0.01100 & 0.00555\\
		& R posterior interosseous & 0.00135 & 0.08253 & 0.00329 & 0.12892\\
		& R femoral II & 1.271e-04 & 0.00247 & 0.03294 & 0.03983\\
		& R anterior tibial & 0.00102 & 0.04334 & 0.01833 & 0.18414\\
		\bottomrule
	\end{tabular}
	\caption{Quantitative assessment of the relative importance of the convective, pressure and friction terms within the momentum balance equation, by computing the ratio between coefficients $\gamma_C$ and $\gamma_P$ and between coefficients $\gamma_F$ and $\gamma_P$, for the vessels of the aortic bifurcation, the 37-artery network and ADAN56 model.} \label{table:dimanalysis}
\end{table}

Equipped with these findings, we conclude that the convective terms can be neglected in the family of nonlinear 0D models derived in Section \ref{sec:0Dmodel} according to the following two main motivations:
\begin{itemize}
	\item the dimensional analysis performed so far shows that the pressure gradient is significantly larger with respect to the convective term in the 1D momentum balance equation;
	\item the ultimate goal of our work is to apply this family of nonlinear 0D models to larger and more complex networks of vessels, such as the global, closed-loop, multiscale model of M\"{u}ller and Toro \cite{Mueller:2014a,Mueller:2014b} and the complete ADAN model developed by Blanco  \textit{et al.} \cite{Blanco:2014a,Blanco:2015a}, in order to construct hybrid 1D-0D networks in the attempt of facing the issues of computational efficiency and execution time. These newly derived 0D models would then be applied not to all vessels of the network, but to small vessels where it is well-known that the convective terms are negligible.
\end{itemize}
Finally, it is worthy to note that numerical experiments (not reported) have shown that including the convective terms into the 0D models is not a straightforward operation. Indeed, numerical difficulties arise in solving the resulting 0D models even by using implicit methods and ODE solvers for stiff problems. We claim that these numerical issues are not related to the instability or stiffness of the ODE systems to be solved, but that the source of these problems lies in the fact that, according to the coupling approach described in Section \ref{sec:0Djunction}, the input data to be prescribed at the inlet/outlet of the vessels converging at the 0D junction are defined by the state of their adjacent compartments, but in each 0D vessel no interaction between the input data and internal vessel state is enforced. This coupling procedure is indeed different from the approach usually adopted for 1D junctions, where in the Riemann problem to be solved at the 1D junction the unknown boundary state vectors are connected not only among themselves, but also to the vessel initial condition states via non-linear waves. As a consequence, this produces ambiguity in determining which are the correct flow rates $Q_L$, $Q_R$ and cross-sections $A_L$, $A_R$ to be used in the convective terms difference originally included in (\ref{sys0D}).\\
In conclusion, dealing with convective terms in 0D blood flow models is clearly an open problem which, to the best of our knowledge, has never been addressed in previous scientific works.  Indeed, in the standard derivation of 0D blood flow models, convective terms are commonly neglected under the assumption that the contribution of the convective terms difference is small compared to the other terms in the momentum balance equation and can thus be discarded. No further discussion is found in the literature about this topic, which remains to be further investigated.

\subsection{Aortic bifurcation model}\label{sec:numexp:test1}

\begin{figure}[h!]
	\begin{center}
		\subfloat[][Abdominal aorta\label{fig:test1:aorta}] {\includegraphics[scale=0.7]{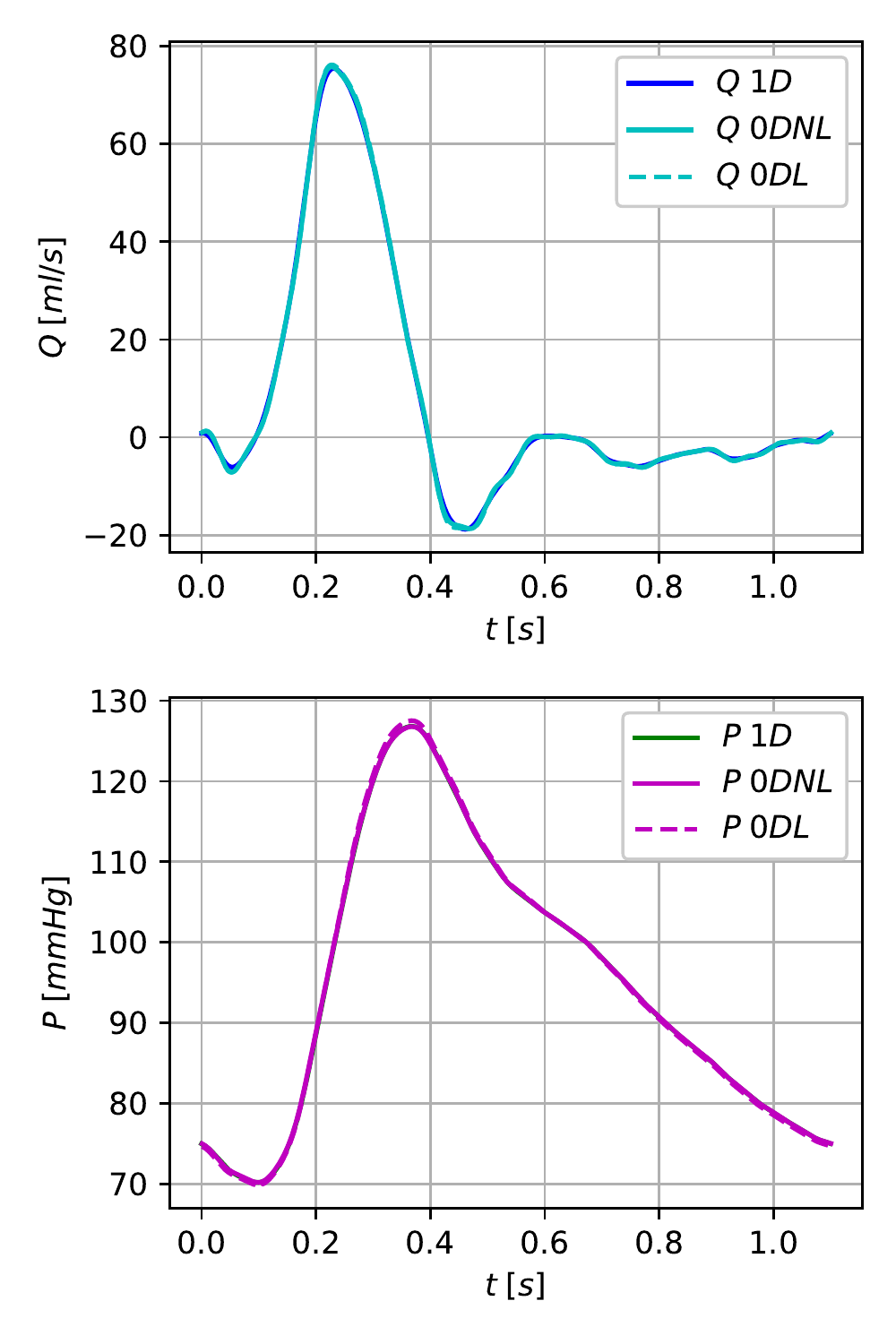}}\qquad
		\subfloat[][Iliac arteries\label{fig:test1:iliac}] {\includegraphics[scale=0.7]{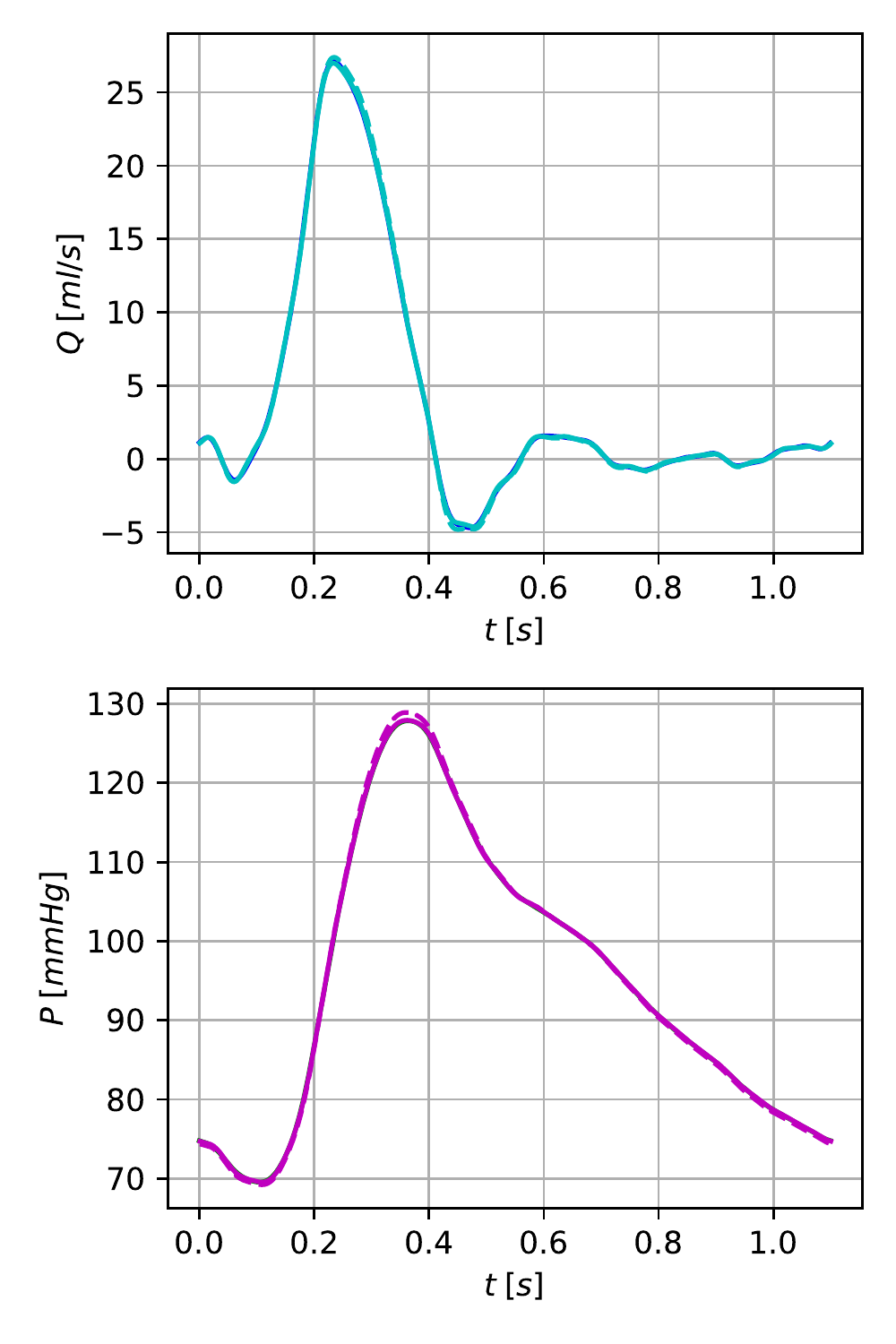}}\\
		\caption{Aortic bifurcation. Comparison between 1D, nonlinear 0D and linear 0D results in the abdominal aorta (A) and in the iliac arteries (B). 1D/Q 1D: 1D numerical solution at the midpoint of the vessel; P 0DNL/Q 0DNL: 0D numerical solution from the nonlinear 0D models; P 0DL/Q 0DL: 0D numerical solution from the linear 0D models.}\label{fig:plot:test1}
	\end{center}
\end{figure}

We simulate the abdominal aorta branching into the two iliac arteries using a single-bifurcation model, consisting of a three-vessel junction \cite{Boileau:2015a}. Both iliac arteries are coupled to a $RCR$ Windkessel terminal element of the rest of the systemic circulation. The geometrical and mechanical properties of this model are summarized in Table \ref{table:param:test1}. 1D/0D initial areas $A(x,0)/\widehat{A}(0)$ are computed using the tube law (\ref{tubelaw}) with $P_0 = P_d$ and $A_0 = A_d$. The inflow boundary condition $Q_{in}(t)$ is an \textit{in vivo} signal taken from \cite{Xiao:2014a} and available in the Supporting Information of \cite{Boileau:2015a}.\\
For the 0D simulations, the abdominal aorta (parent vessel) is discretized using a $(Q_{in}, Q_{out})$-type 0D vessel, while both iliac arteries (daughter vessels) are represented as $(P_{in}, P_{out})$-type 0D vessels. Flow rate $Q_{in}$ at the inlet of the parent vessel is given by the periodic inflow boundary condition; flow rate $Q_{out}$ at the outlet of the aorta and pressure $P_{in}$ at the inlet of both iliac arteries are determined by imposing the coupling conditions at a 0D junction; pressure $P_{out}$ at the outlet of both parent vessels is defined by the adjacent $RCR$ Windkessel terminal element. The computational parameters of the 1D/0D simulations are chosen according to Table \ref{table:param}.\\
Figure \ref{fig:plot:test1} shows a qualitative comparison of pressure and flow waveforms for the aortic bifurcation model. We can observe that, compared to the reference 1D results, the amplitude and the shape of these waveforms are well-captured by both the nonlinear and linear 0D results, with differences between the two 0D models that can not be appreciated in the scale of the figures. However, the 0D solution becomes less accurate when using the linear pressure-volume relation and constant parameters $C_0$, $R_0$ and $L_0$ in the 0D simulations. This is evident in the systolic phase of the pressure waveform, where the systolic peak corresponding to the linear 0D results is higher with respect to the reference 1D results, while the pressure peak is perfectly captured by the nonlinear 0D model. These observations are confirmed also by the quantitative assessment. Relative errors were determined with respect to the 1D solution and are displayed in Table \ref{table:errors:test1}. Relative errors for pressure and flow rate are consistently very small, relative root mean square errors are all smaller than 0.5\% for the pressure and smaller than 1\% for the flow rate. Furthermore, relative errors in the nonlinear 0D results are in general smaller than the corresponding relative errors in the linear 0D results.\\
\begin{table}[h!]\footnotesize
	\centering
	\begin{tabular}{lcc}
		\toprule
		\textbf{Property} & \multicolumn{2}{c}{\textbf{Value}}\\
		\midrule
		Blood density, $\rho$ & \multicolumn{2}{c}{1.060 g/cm$^3$} \\
		Blood viscosity, $\mu$ & \multicolumn{2}{c}{0.04 dyne$\cdot$s/cm$^2$} \\
		Velocity profile order, $\zeta$ & \multicolumn{2}{c}{9} \\
		Diastolic pressure, $P_d$ & \multicolumn{2}{c}{9.4$\hat{6}$ $\cdot$10$^4$ dyne/cm$^2$} \\
		External pressure, $p_{ext}$ & \multicolumn{2}{c}{0}\\
		\midrule
		& \textbf{Aorta} & \textbf{Iliac}\\
		\midrule
		Length, $l$ & 8.6 cm & 8.5 cm\\
		Radius at diastolic pressure, $r_d$ & 0.86 cm & 0.60 cm\\
		Area at diastolic pressure, $A_d$ & 2.3235 cm$^2$ & 1.1310 cm$^2$\\
		Wall thickness, $h$ & 1.032 mm & 0.72 mm\\
		Young's modulus, $E$ & 5.0$\cdot$10$^6$ dyne/cm$^2$ & 7.0$\cdot$10$^6$ dyne/cm$^2$\\
		WK resistance, $R_1$ & - & 6.8123$\cdot$10$^2$ dyne$\cdot$s/cm$^5$\\
		WK compliance, $C$ & - & 3.6664$\cdot$10$^{-5}$ cm$^5$/dyne\\
		WK resistance, $R_2$ & - & 3.1013$\cdot$10$^4$ dyne$\cdot$s/cm$^5$\\
		Outflow pressure, $P_{out}$ & - & 0\\
		\midrule
		Initial cross-sectional area, $A(x, 0)/\widehat{A}(0)$ & 1.8062 cm$^2$ & 0.94789 cm$^2$\\
		Initial velocity, $u(x, 0)/U(0)$ & 0 & 0\\
		Initial pressure, $p(x, 0)/P(0)$ & 0 & 0\\
		\bottomrule
	\end{tabular}
	\caption{Geometrical and mechanical properties of the aortic bifurcation model.}\label{table:param:test1}
\end{table}

\begin{table}[h!]\footnotesize
	\centering
	\renewcommand\arraystretch{1.2}
	\begin{tabular}{l | c | cccccc}
	\toprule
	\textbf{Vessel} & \textbf{0D model} & $\varepsilon_P^{\text{RMS}} (\%)$ & $\varepsilon_Q^{\text{RMS}} (\%)$ & $\varepsilon_P^{\text{SYS}} (\%)$ & $\varepsilon_Q^{\text{SYS}} (\%)$ & $\varepsilon_P^{\text{DIAS}} (\%)$ & $\varepsilon_Q^{\text{DIAS}} (\%)$ \\
	\midrule
	\multirow{2}{*}{Abdominal aorta} & 0D-NL & 0.231 & 0.620 & -0.027 & 0.424 & -0.013 & 0.127\\
	& 0D-L & 0.437 & 0.680 & 0.551 & 0.970 & -0.482 & 0.092\\
	\midrule
	\multirow{2}{*}{Iliac arteries} & 0D-NL & 0.149 & 0.406 & 0.051 & -0.484 & 0.010 & 0.212\\
	& 0D-L & 0.462 & 0.695 & 0.823 & 1.010 & -0.458 & -0.445\\
	\bottomrule
	\end{tabular}
	\caption{Aortic bifurcation. Relative errors (in \%) for pressure and flow between for both nonlinear (0D-NL) and linear (0D-L) 0D results with respect to the 1D results (1D) at the midpoint of the vessel, computed according to the relative error metrics (\ref{metrics}).}\label{table:errors:test1}
\end{table}

\subsection{37-artery network}\label{sec:numexp:test2}

We consider the arterial network presented in \cite{Matthys:2007a}, for which \textit{in vitro} pressure and flow measurements were acquired at multiple locations. This arterial tree is made up of 37 silicone vessels representing the largest central systemic arteries of the human vascular system. For the detailed topology of the network, the reader is referred to \cite{Matthys:2007a}. At the inlet of the ascending aorta, the flow rate measured \textit{in vitro} is prescribed as the inflow boundary condition $Q_{in}(t)$. Terminal vessels are coupled to single-resistance terminal models. In all vessels, initial 1D/0D cross-sectional areas $A(x,0)/\widehat{A}(0)$, corresponding to the zero initial pressure prescribed, are computed using the tube law (\ref{tubelaw}) with reference pressure $P_0 = P_d = 0$ and cross-sectional area $A_0 = A_d$. General parameters of this arterial network are given in Table \ref{table:test2}. For a complete set of parameters, we refer the reader to the Supplementary Information of \cite{Boileau:2015a}. We note that, differently from \cite{Matthys:2007a,Boileau:2015a}, the 1D results we present in this work do not include vessel tapering. For each vessel, the reference constant radius $r_0$, from which the reference cross-sectional area $A_0$ is computed, is obtained as the mean value between the cross-sectional radii at the inlet and outlet of the vessel, $r_{in}$ and $r_{out}$ respectively, given in \cite{Matthys:2007a,Boileau:2015a}. We point out that this is a limitation of our study that will be addressed in future works.\\

\begin{figure}[h!]
	\begin{center}
		\subfloat[][Aortic arch II] {\includegraphics[scale=0.7]{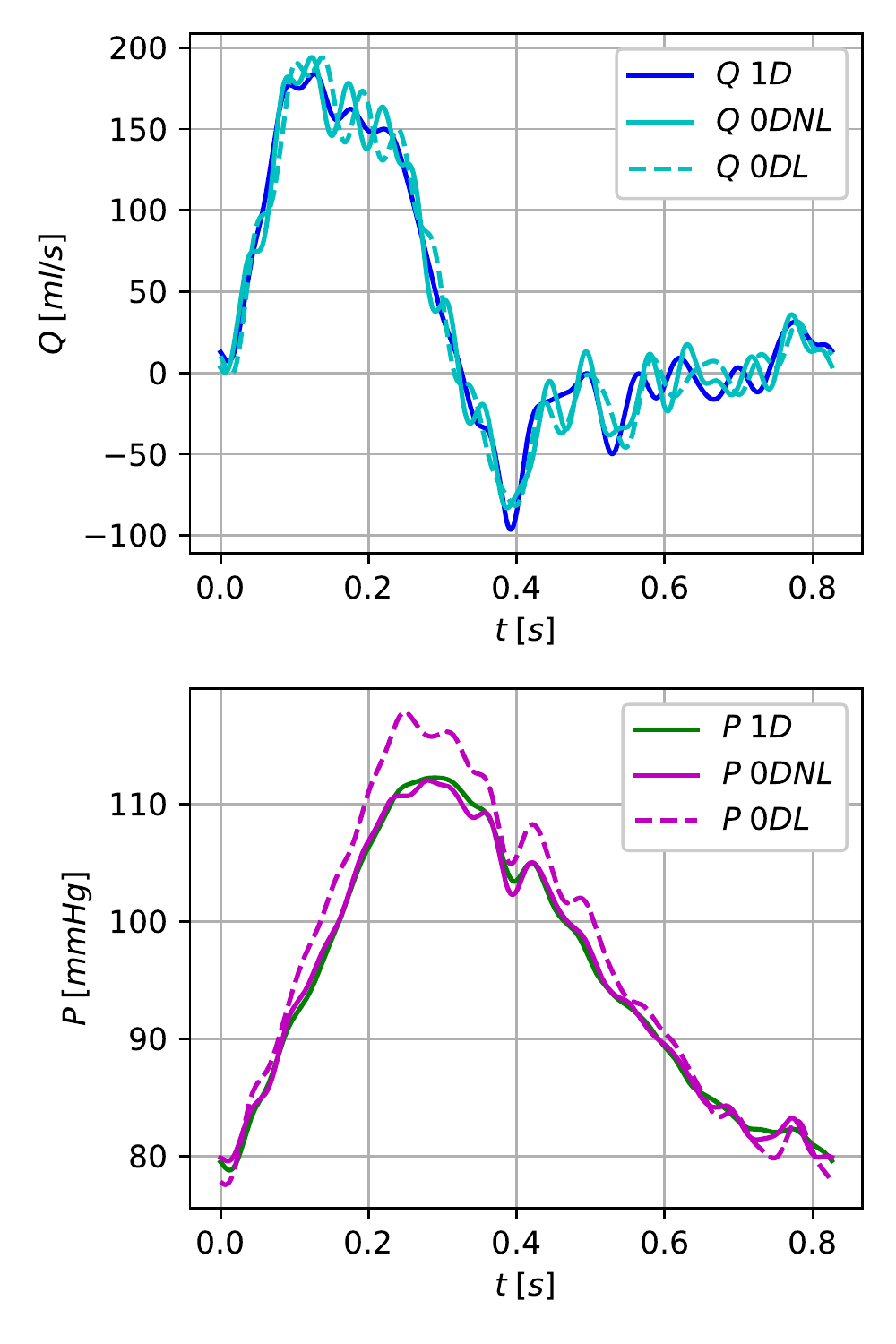}} 
		\subfloat[][Thoracic aorta II] {\includegraphics[scale=0.7]{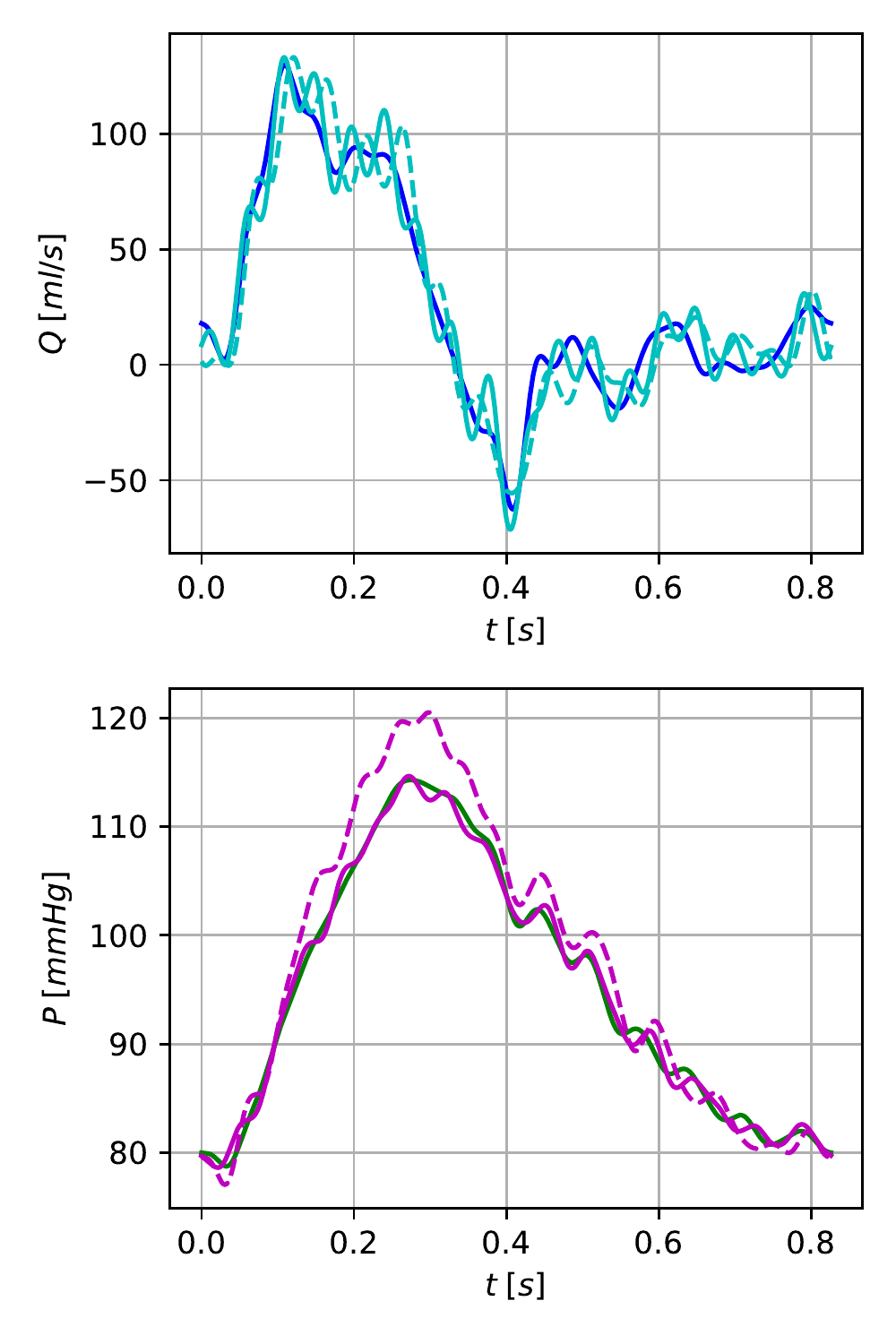}}\\
		\caption{37-artery network. Comparison between 1D (at vessel midpoint), nonlinear 0D and linear 0D results in two aortic segments.}\label{fig:plot:test2:aorta}
	\end{center}
\end{figure}

\begin{figure}[h!]
	\begin{center}
		\subfloat[][Left subclavian I] {\includegraphics[scale=0.7]{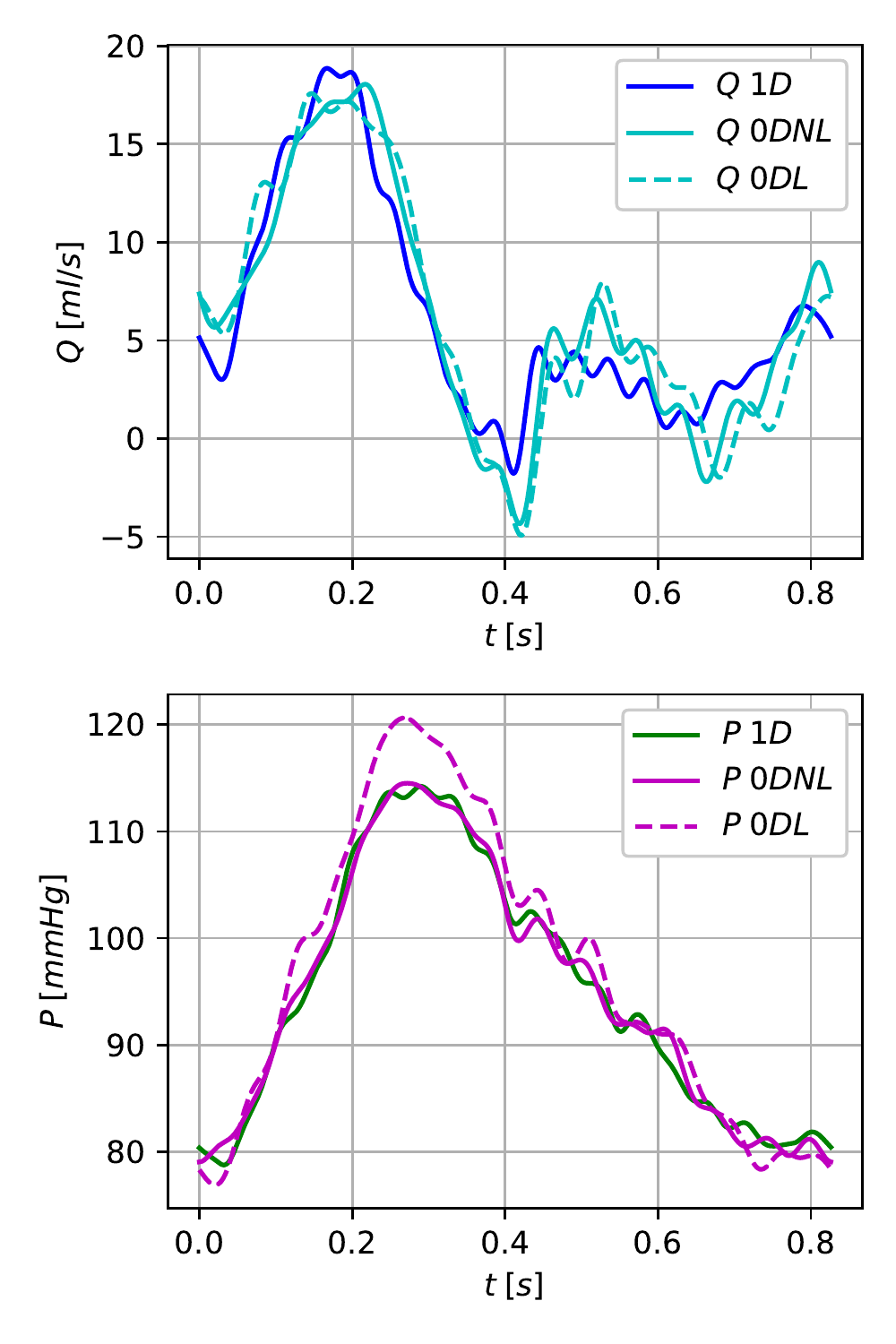}} 
		\subfloat[][Right iliac-femoral II] {\includegraphics[scale=0.7]{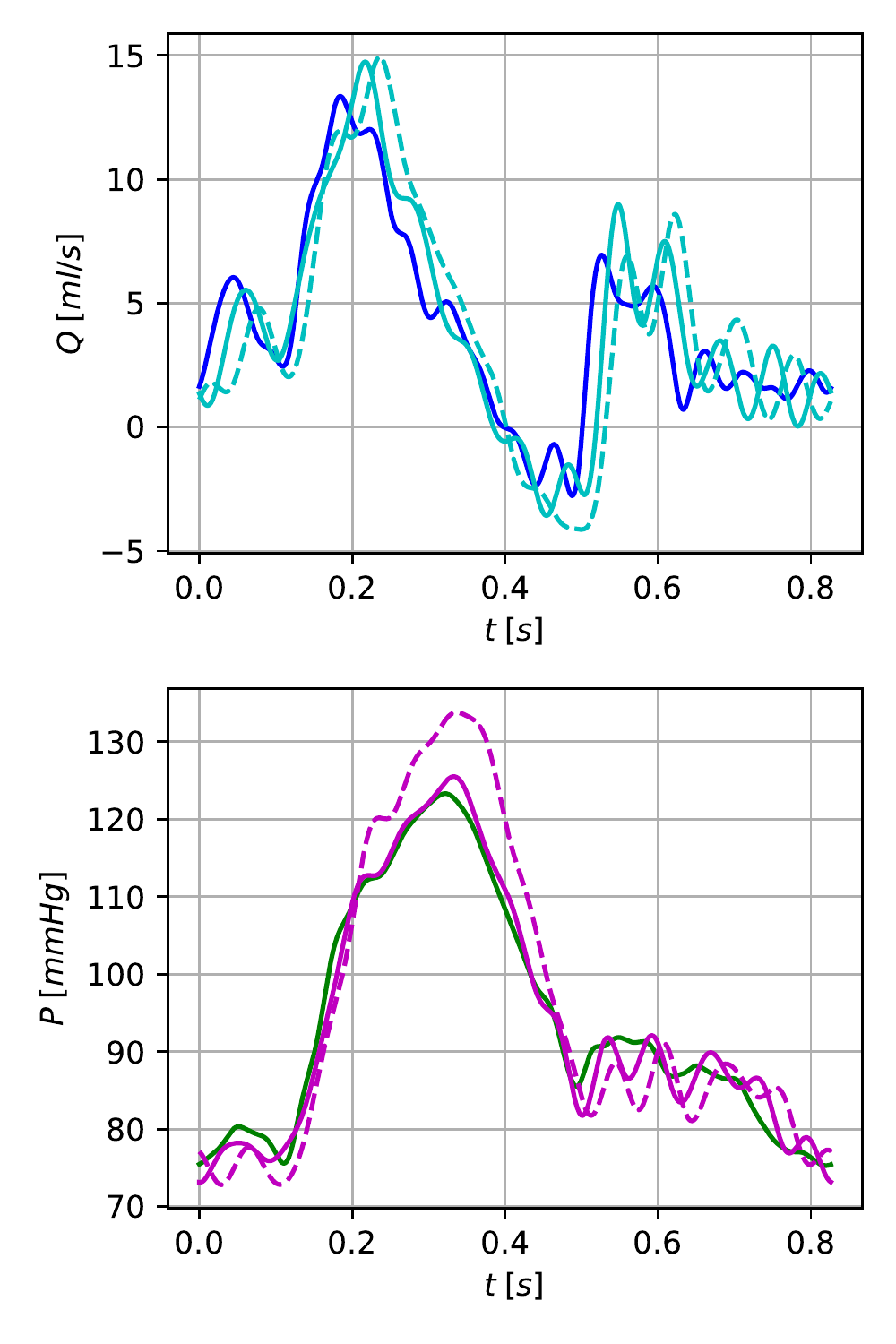}}\\
		\caption{37-artery network. Comparison between 1D (at vessel midpoint), nonlinear 0D and linear 0D results in two first-generation vessels.}\label{fig:plot:test2:gen1}
	\end{center}
\end{figure}

\begin{figure}[h!]
	\begin{center}
		\subfloat[][Left ulnar] {\includegraphics[scale=0.7]{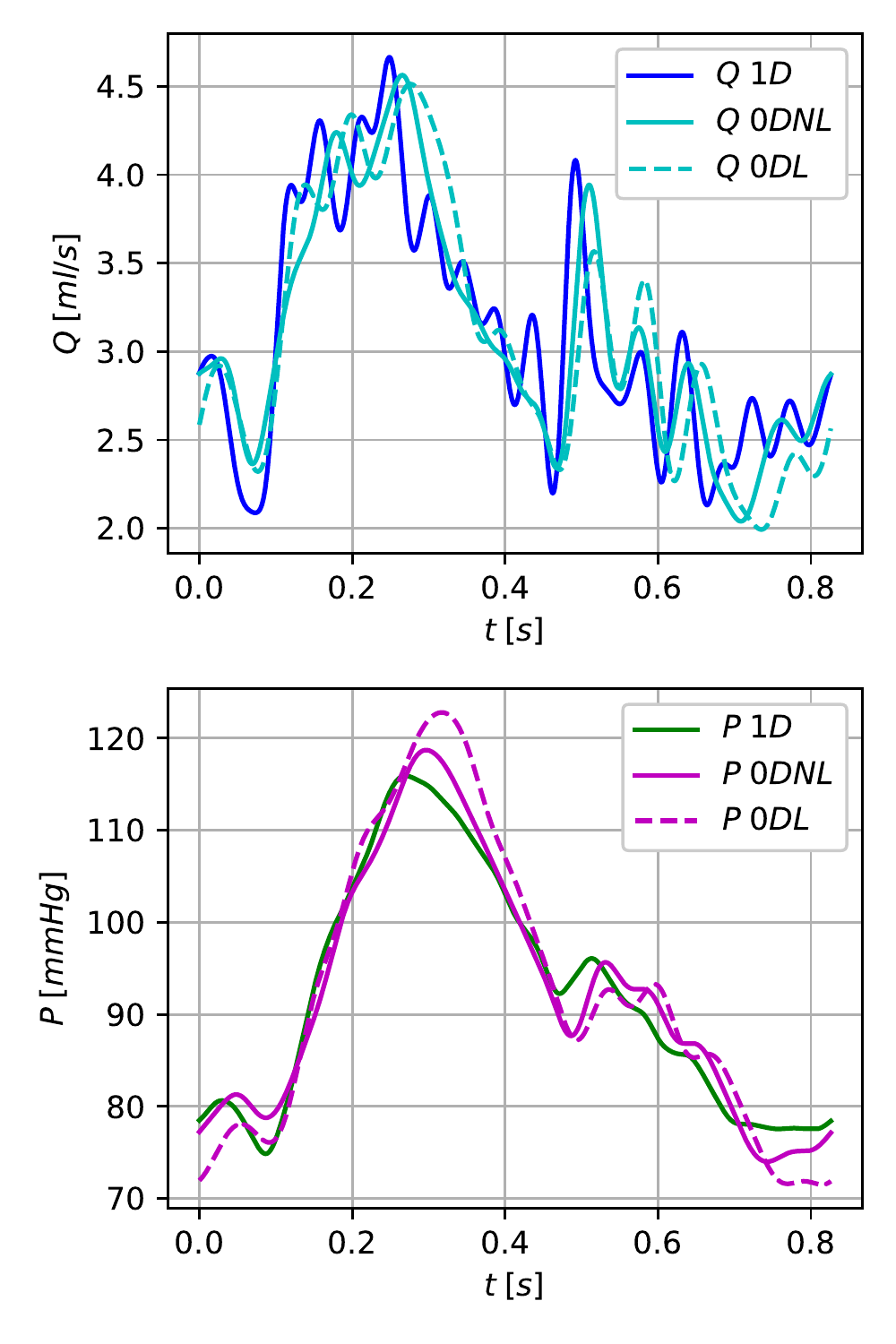}} 
		\subfloat[][Right anterior tibial] {\includegraphics[scale=0.7]{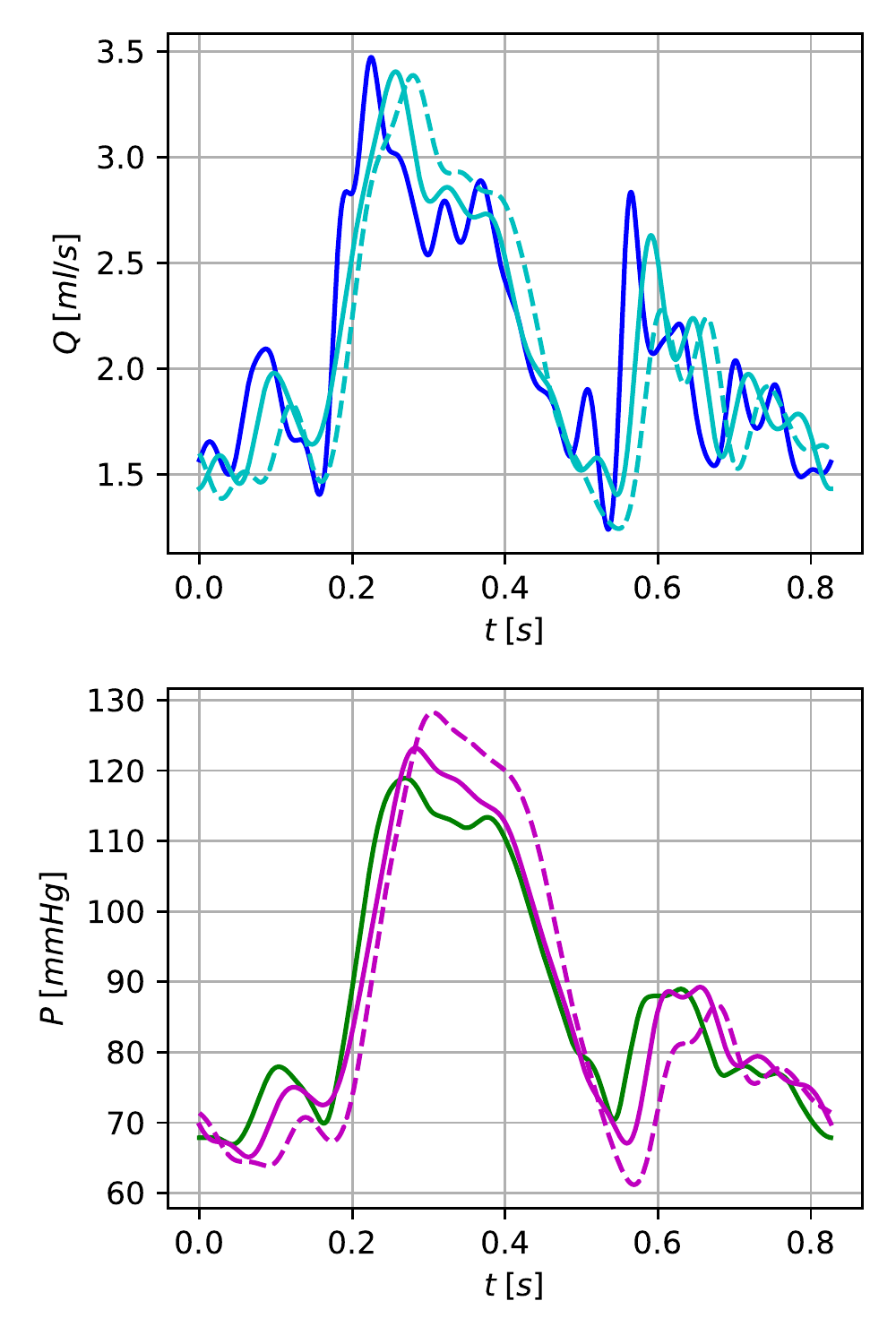}}\\
		\caption{37-artery network. Comparison between 1D (at vessel midpoint), nonlinear 0D and linear 0D results in two second-generation vessels.}\label{fig:plot:test2:gen2}
	\end{center}
\end{figure}

\begin{figure}[h!]
	\begin{center}
		\subfloat[][Right ulnar] {\includegraphics[scale=0.7]{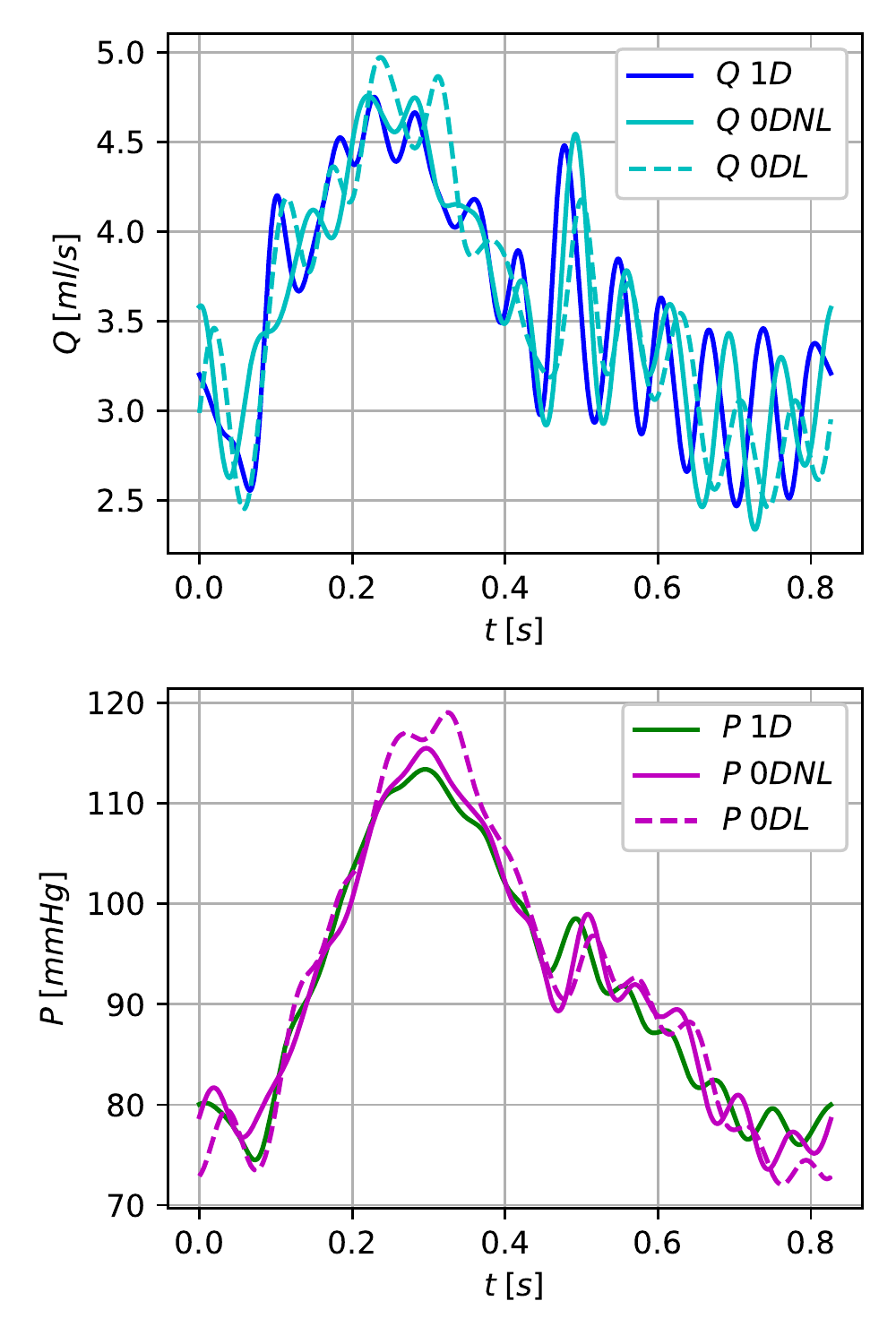}} 
		\subfloat[][Splenic] {\includegraphics[scale=0.7]{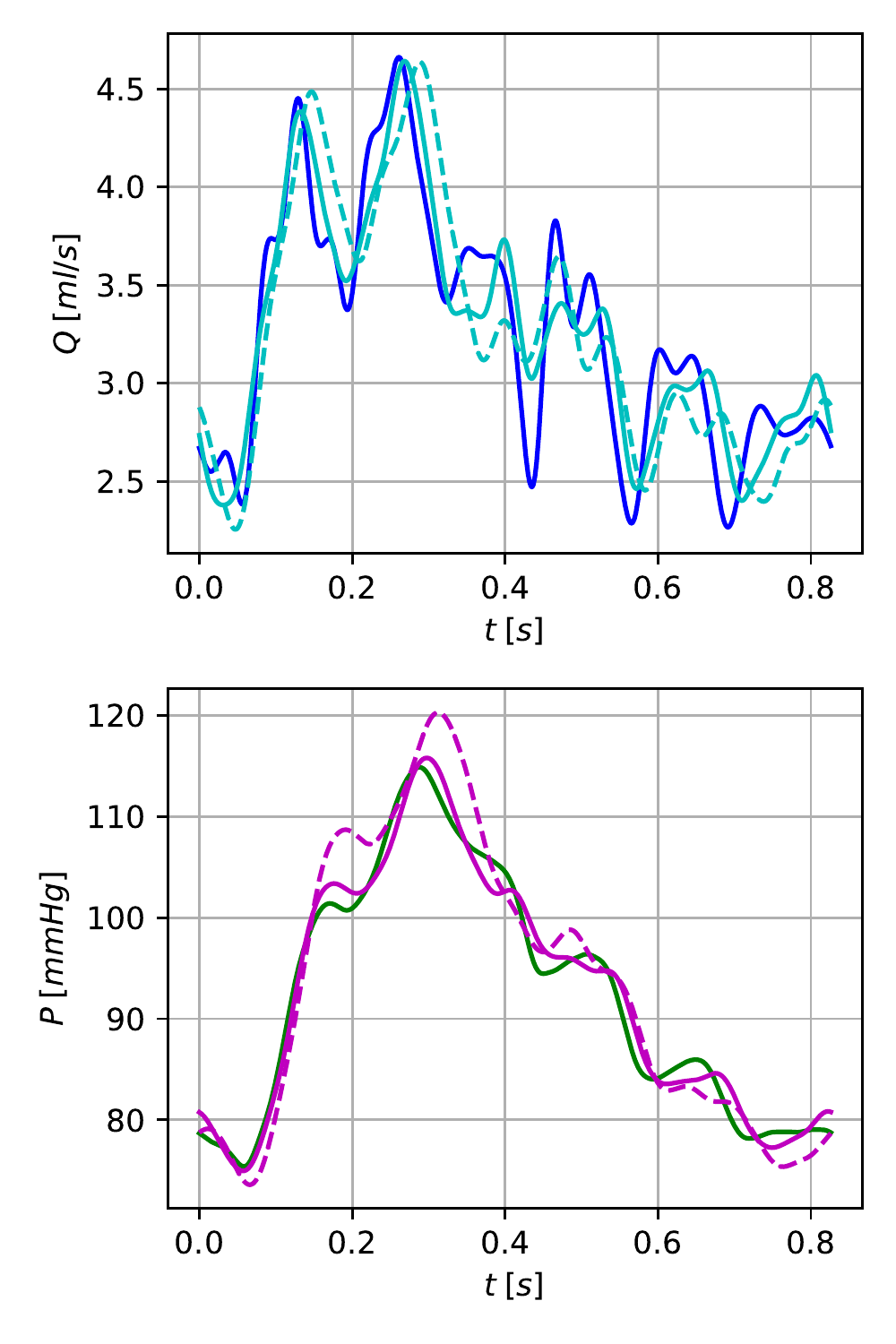}}\\
		\caption{37-artery network. Comparison between 1D (at vessel midpoint), nonlinear 0D and linear 0D results in two third-generation vessels.}\label{fig:plot:test2:gen3}
	\end{center}
\end{figure}

For the 0D simulations, the first vessel of the arterial network, the ascending aorta, is discretized using a $(Q_{in}, Q_{out})$-type 0D vessel, terminal vessels, which are coupled to single-resistance terminal elements, are represented as $(P_{in}, P_{out})$-type 0D vessels, while all other vessels, which are not at the extremities of the network, are modelled as two-split $(P_{in}, Q_{out})$-type 0D vessels. Results are shown for two aortic segments (aortic arch II and thoracic aorta II), two vessels of the first generation of bifurcations (left subclavian I and right iliac-femoral II), two of the second generation (left ulnar and right anterior tibial) and two of the third generation (right ulnar and splenic). Qualitative comparisons between nonlinear 0D, linear 0D and reference 1D solutions are shown in Figure \ref{fig:plot:test2:aorta} for the aortic segments and in Figures \ref{fig:plot:test2:gen1}-\ref{fig:plot:test2:gen3} for the vessels of first, second and third generations, respectively. Table \ref{table:errors:test2} displays the relative errors computed for both nonlinear and linear 0D results with respect to 1D results.\\
As expected, when we move to a more complex network, the differences between 0D and 1D results become more significant, even when adopting the nonlinear 0D models. In Figures \ref{fig:plot:test2:aorta}-\ref{fig:plot:test2:gen3} we observe that, especially in the flow waveforms, some oscillations are amplified, while other oscillations are not captured by the 0D models with respect to the reference 1D solution. However, if we focus on the nonlinear 0D results, we can conclude that they are quite satisfactory: overall, the predicted pressure and flow waveforms are in good agreement with the 1D results and, given the complexity of the flow to be simulated, the essential features, shape and amplitude, of these waves are well-captured. Furthermore, from this benchmark problem we clearly see how the nonlinearity included in the 0D models, through the nonlinear pressure-area relation and the nonlinear parameters, strongly improves the 0D results, especially the pressure waveforms: in general, the linear 0D solution overestimates the systolic peak in pressure, while the pressure peak reproduced by the nonlinear 0D results is much more in agreement with the reference 1D results.\\
All these observations are confirmed also by the quantitative assessment presented in Table \ref{table:errors:test2}. Relative errors in the nonlinear 0D results are overall smaller than the corresponding relative errors in the linear 0D results. Indeed, when moving from the linear to the nonlinear 0D models, the RMS relative errors decrease and, in particular, the systolic relative errors in pressure are significantly reduced, suggesting that overall the proposed family of nonlinear 0D models is able to better reproduce pressure and flow waveforms and to capture their essential features, such as the systolic peak of the pressure wave. Indeed, since the nonlinearity mostly enters in the evaluation of the pressure via the nonlinear pressure-area relation, this improvement is mainly reflected on the pressure waveforms, as expected. For the nonlinear 0D results, RMS relative errors are all smaller than 6.0\% for the pressure and smaller than 15.0\% for the flow rate; systolic relative errors are all smaller than 4.0\% for the pressure and smaller than 6.0\% for the flow rate, with an exception in the right iliac-femoral II artery, where the systolic peak in the flow rate seems to be slightly overestimated.\\
\begin{table}[h!]\footnotesize
	\centering
	\begin{tabular}{lc}
		\toprule
		\textbf{Property} & \textbf{Value} \\
		\midrule
		Blood density, $\rho$ & 1.050 g/cm$^3$ \\
		Blood viscosity, $\mu$ & 0.025 dyne$\cdot$s/cm$^2$ \\
		Velocity profile order, $\zeta$ & 9 \\
		Young's modulus, $E$ & 1.2$\cdot$10$^7$ dyne/cm$^2$ \\
		Diastolic pressure, $P_d$ & 0 \\
		External pressure, $p_{ext}$ & 0 \\
		Outflow pressure, $P_{out}$ & 4.2663$\cdot$10$^3$ dyne/cm$^2$ (3.2 mmHg)\\
		Initial velocity, $u(x, 0)/U(0)$ & 0\\
		Initial pressure, $p(x, 0)/P(0)$ & 0\\
		\bottomrule
	\end{tabular}
	\caption{General model parameters of the 37-artery network.}\label{table:test2}
\end{table}

\begin{table}[h!]\footnotesize
	\centering
	\renewcommand\arraystretch{1.2}
	\begin{tabular}{l | c | cccccc}
		\toprule
		\textbf{Vessel} & \textbf{0D model} & $\varepsilon_P^{\text{RMS}} (\%)$ & $\varepsilon_Q^{\text{RMS}} (\%)$ & $\varepsilon_P^{\text{SYS}} (\%)$ & $\varepsilon_Q^{\text{SYS}} (\%)$ & $\varepsilon_P^{\text{DIAS}} (\%)$ & $\varepsilon_Q^{\text{DIAS}} (\%)$ \\
		\midrule
		\multirow{2}{*}{Aortic arch II} & 0D-NL & 0.638 & 6.218 & -0.205 & 5.532 & 0.983 & 7.135\\
		& 0D-L & 2.780 & 7.641 & 4.958 & 5.443 & -1.568 & 7.967\\
		\midrule
		\multirow{2}{*}{Thoracic aorta II} & 0D-NL & 0.868 & 7.803 & 0.308 & 2.543 & -0.161 & -6.772 \\
		& 0D-L & 3.028 & 10.085 & 5.438 & 2.576 & -2.150 & 5.367 \\
		\midrule
		\multirow{2}{*}{L subclavian I} & 0D-NL & 1.193 & 10.638 & 0.251 & -4.304 & 0.369 & -13.605\\
		& 0D-L & 3.114 & 12.584 & 5.610 & -6.796 & -2.372 & -16.745\\
		\midrule
		\multirow{2}{*}{R iliac-femoral II} & 0D-NL & 2.560 & 13.841 & 1.786 & 10.405 & -2.846 & -6.015\\
		& 0D-L & 6.286 & 20.628 & 8.489 & 11.780 & -3.264 & -10.194\\
		\midrule
		\multirow{2}{*}{L ulnar} & 0D-NL & 2.839 & 7.199 & 2.401 & -2.167 & -1.141 & -1.025\\
		& 0D-L & 4.653 & 9.115 & 5.924 & -3.247 & -4.483 & -2.019\\
		\midrule
		\multirow{2}{*}{R anterior tibial} & 0D-NL & 5.291 & 7.785 & 3.600 & -1.957 & -2.675 & 4.733\\
		& 0D-L & 10.635 & 11.526 & 7.854 & -2.466 & -8.565 & 0.104\\
		\midrule
		\multirow{2}{*}{R ulnar} & 0D-NL & 2.766 & 8.322 & 1.855 & 0.160 & -1.248 & -2.766\\
		& 0D-L & 3.829 & 8.849 & 4.992 & 4.683 & -3.301 & -0.366\\
		\midrule
		\multirow{2}{*}{Splenic} & 0D-NL & 1.748 & 4.970 & 0.805 & -0.422 &-0.584 & 2.417 \\
		& 0D-L & 3.571 & 7.517 & 4.717 & -0.519 & -2.388 & -0.251 \\
		\bottomrule
	\end{tabular}
	\caption{37-artery network. Relative errors (in \%) for pressure and flow between for both nonlinear (0DNL) and linear (0DL) 0D results with respect to the 1D results (1D) at the midpoint of the vessel, computed according to the relative error metrics (\ref{metrics}).} \label{table:errors:test2}
\end{table}

\subsection{Reduced ADAN56 model}\label{sec:numexp:test3}

\begin{figure}[h!]
	\begin{center}
		\subfloat[][Aortic arch I] {\includegraphics[scale=0.55]{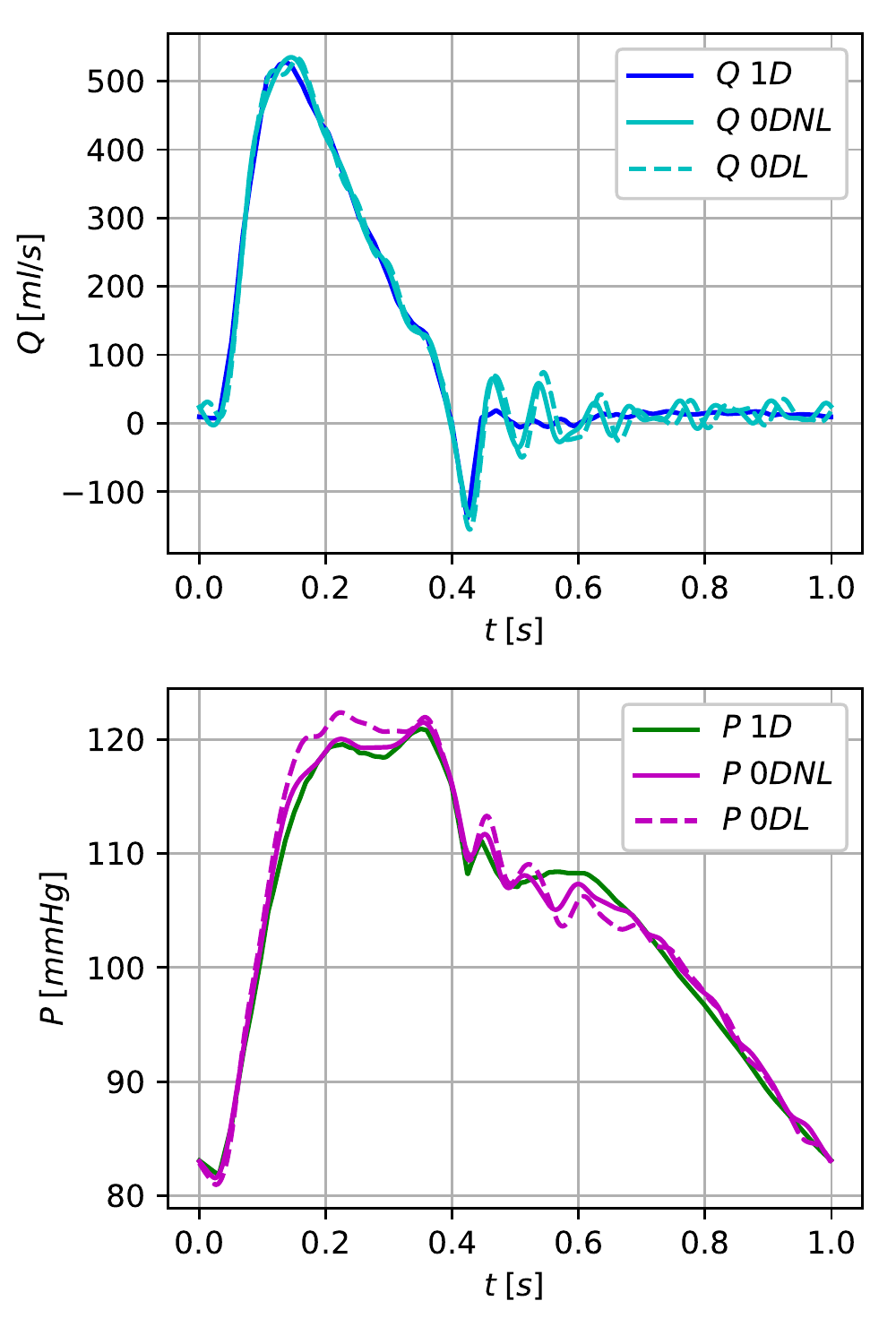}} 
		\subfloat[][Thoracic aorta III] {\includegraphics[scale=0.55]{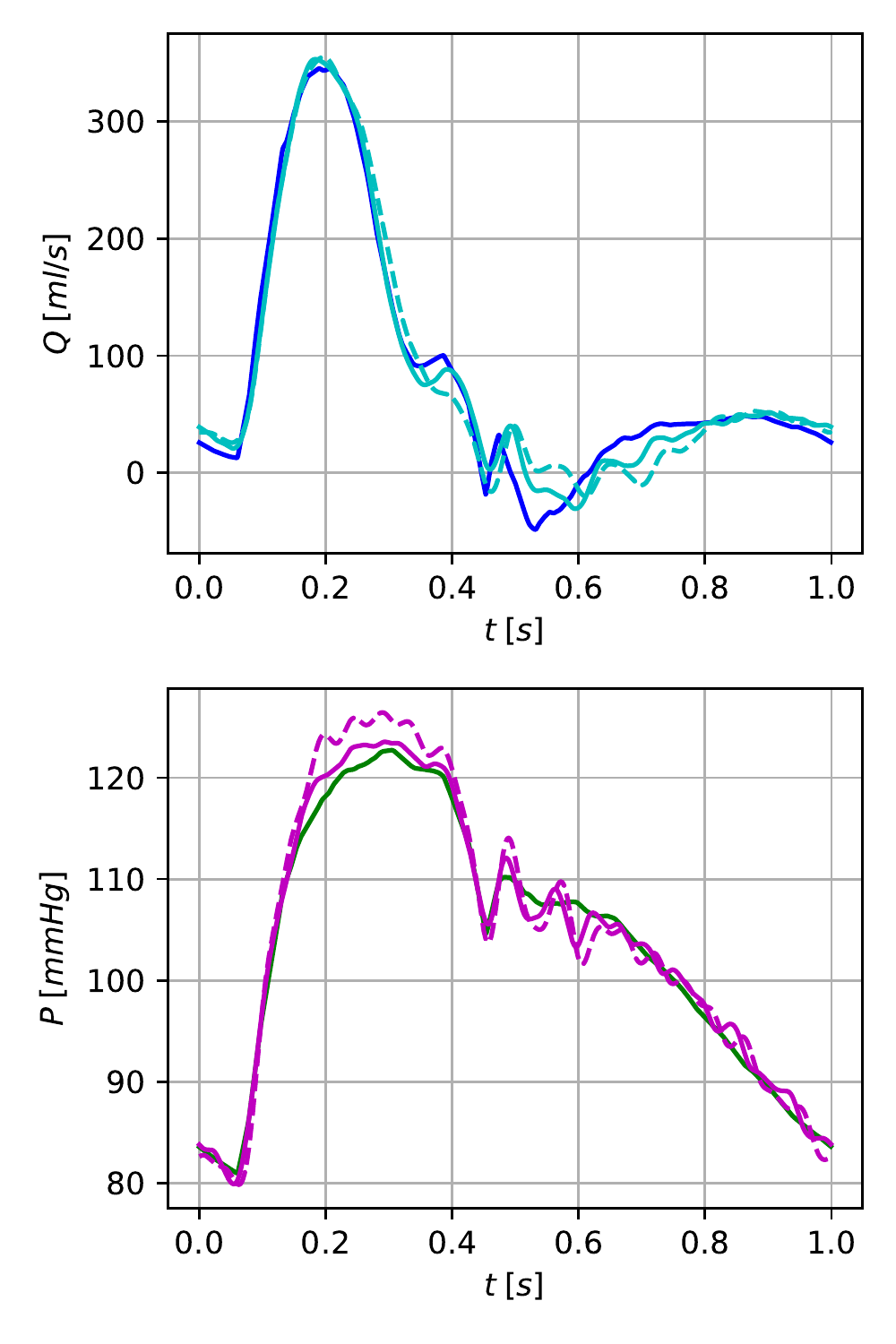}} 
	    \subfloat[][Abdominal aorta V] {\includegraphics[scale=0.55]{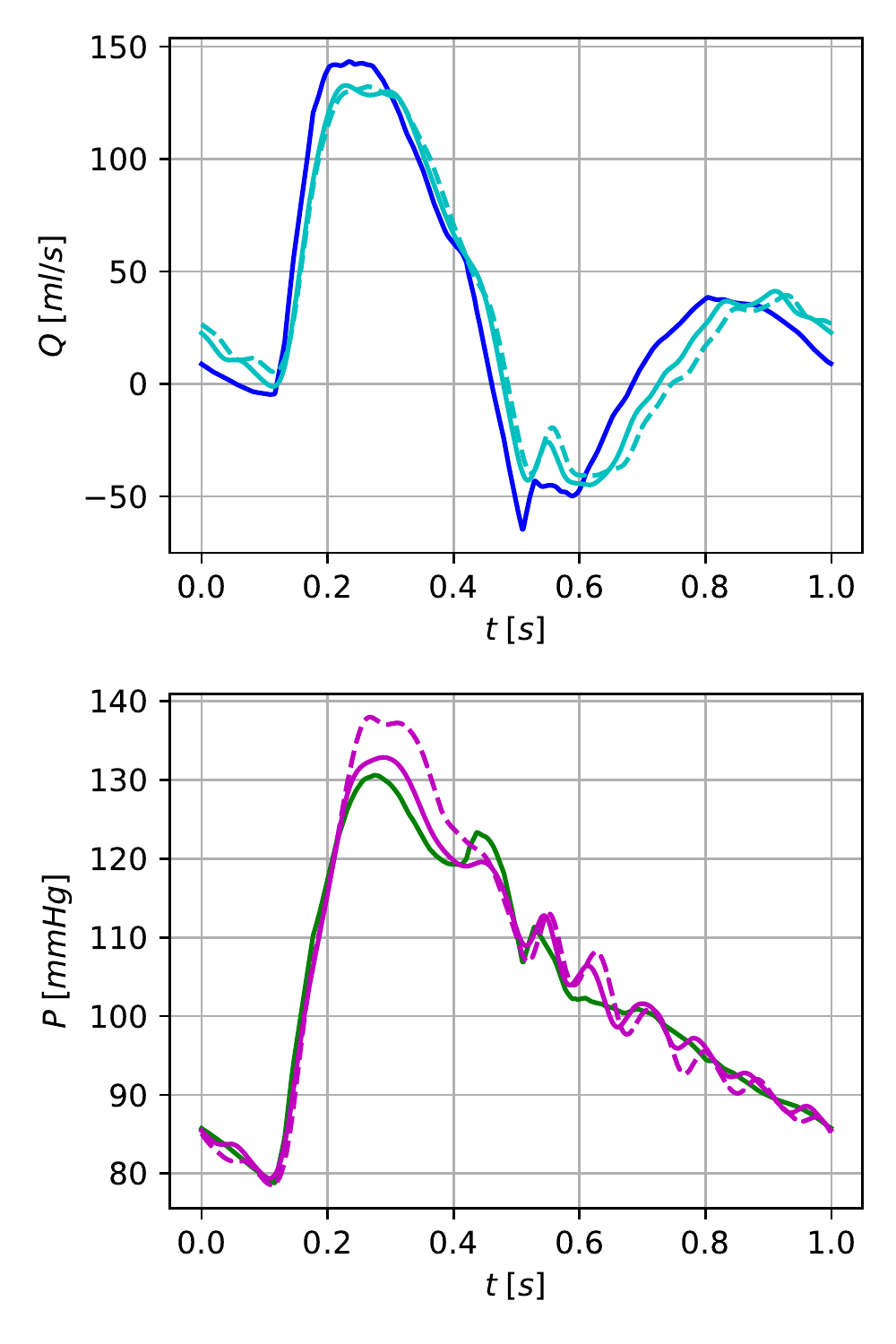}}	
		\caption{ADAN56 model. Comparison between 1D (at vessel midpoint), nonlinear 0D and linear 0D results in three aortic segments.}\label{fig:plot:test3:aorta}
	\end{center}
\end{figure}

\begin{figure}[h!]
	\begin{center}
		\subfloat[][Right common carotid] {\includegraphics[scale=0.55]{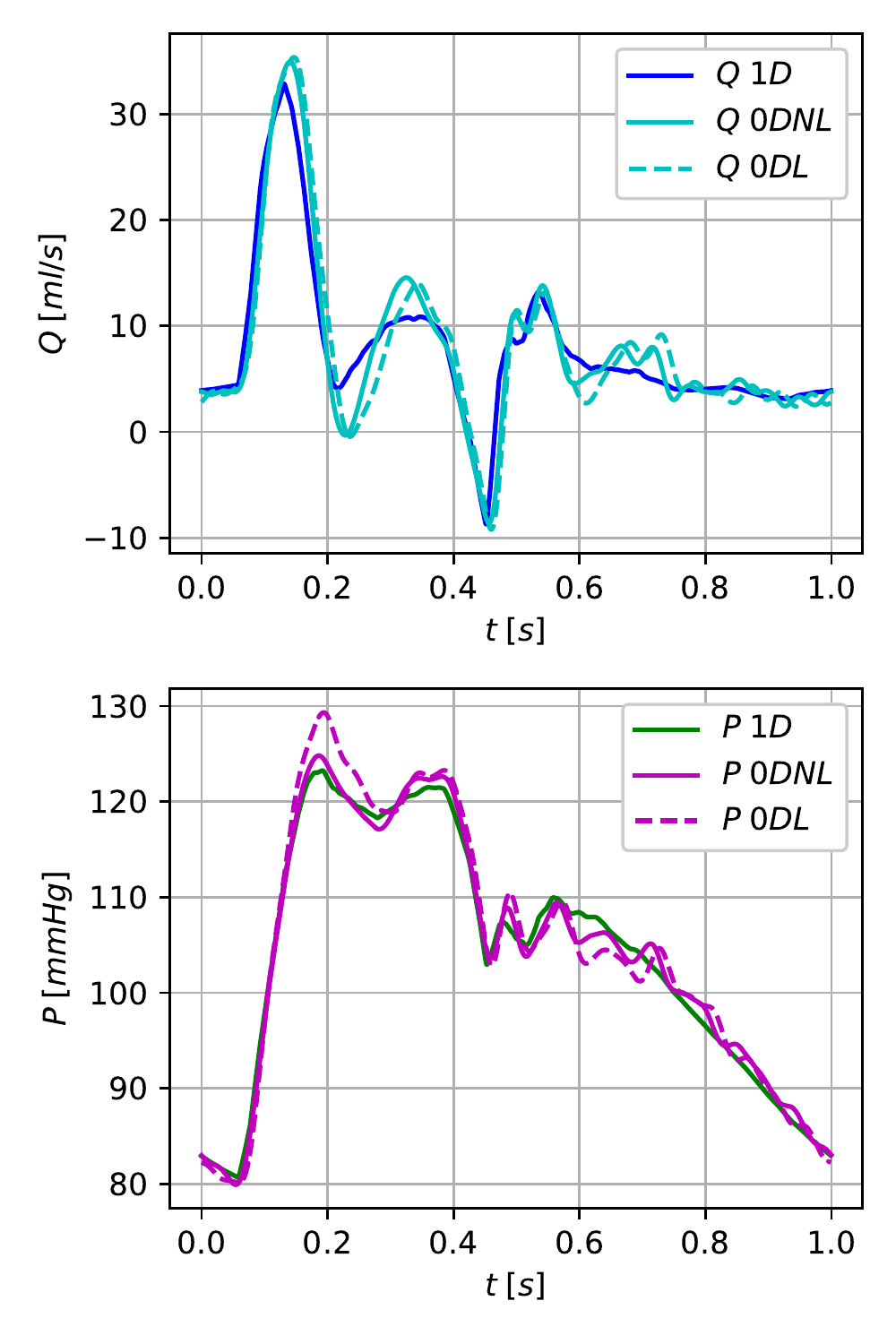}} 
		\subfloat[][Right renal] {\includegraphics[scale=0.55]{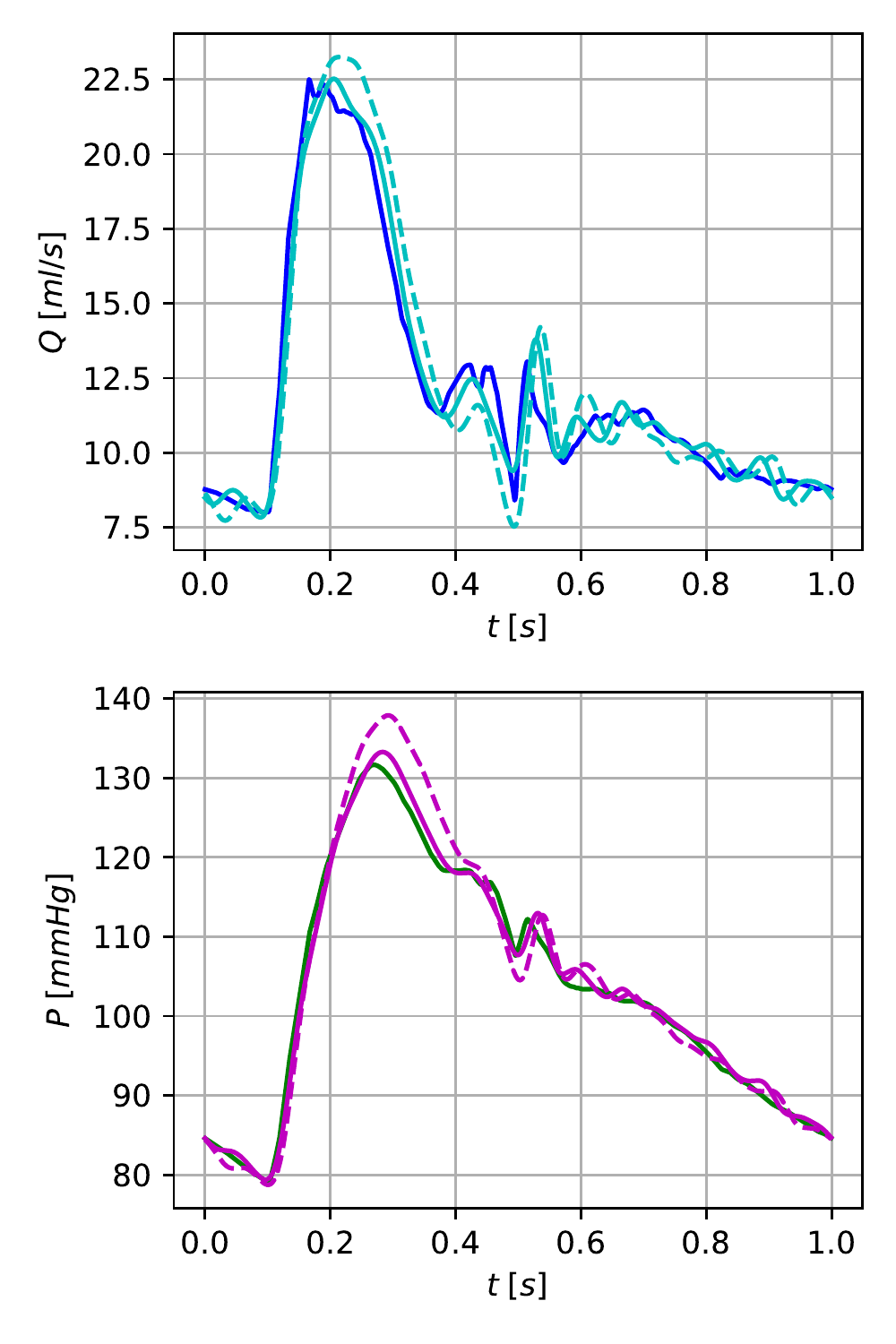}} 
		\subfloat[][Right common iliac] {\includegraphics[scale=0.55]{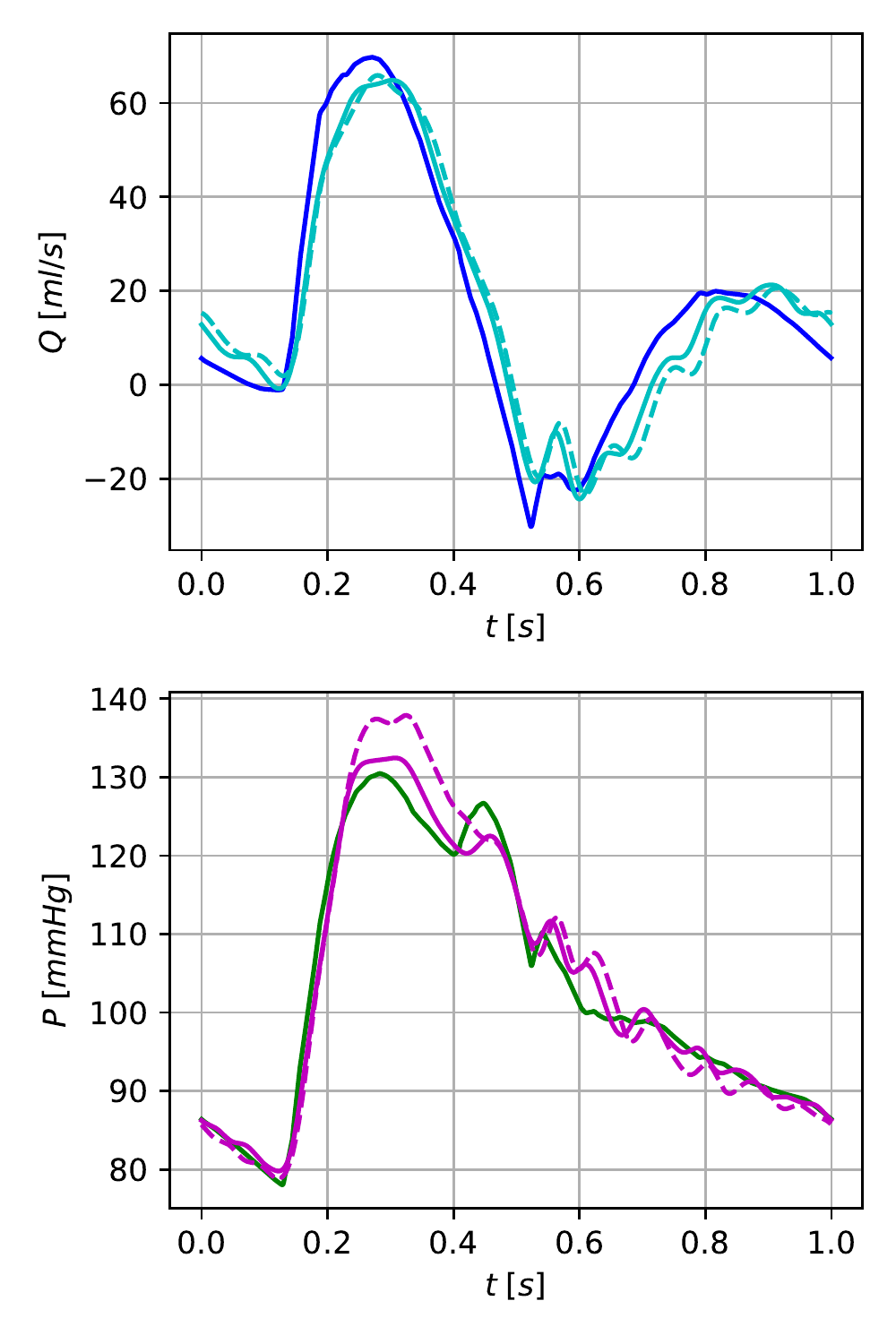}}	
		\caption{ADAN56 model. Comparison between 1D (at vessel midpoint), nonlinear 0D and linear 0D results in three first-generation vessels.}\label{fig:plot:test3:gen1}
	\end{center}
\end{figure}

\begin{figure}[h!]
	\begin{center}
		\subfloat[][Right internal carotid] {\includegraphics[scale=0.55]{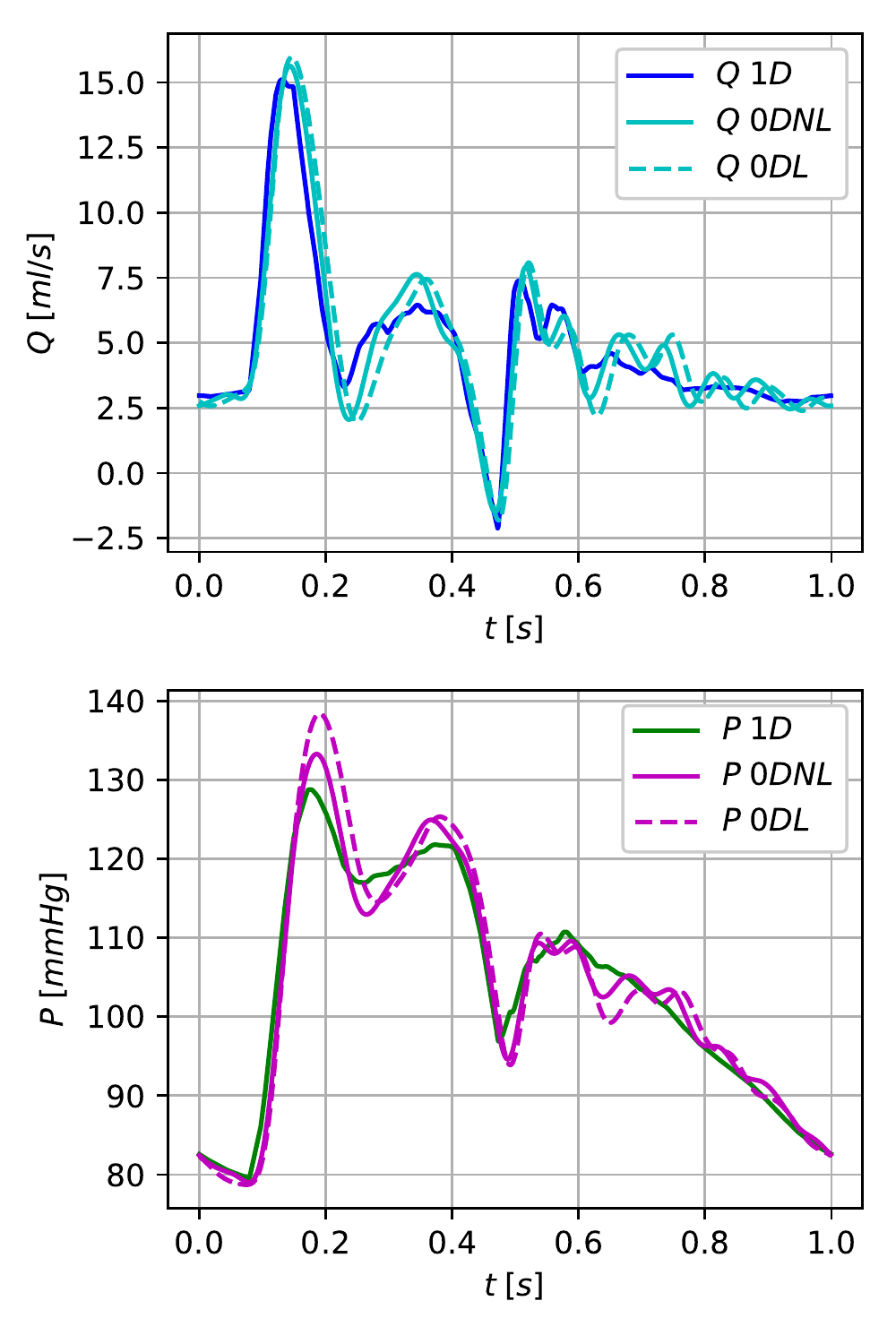}} 
		\subfloat[][Right radial] {\includegraphics[scale=0.55]{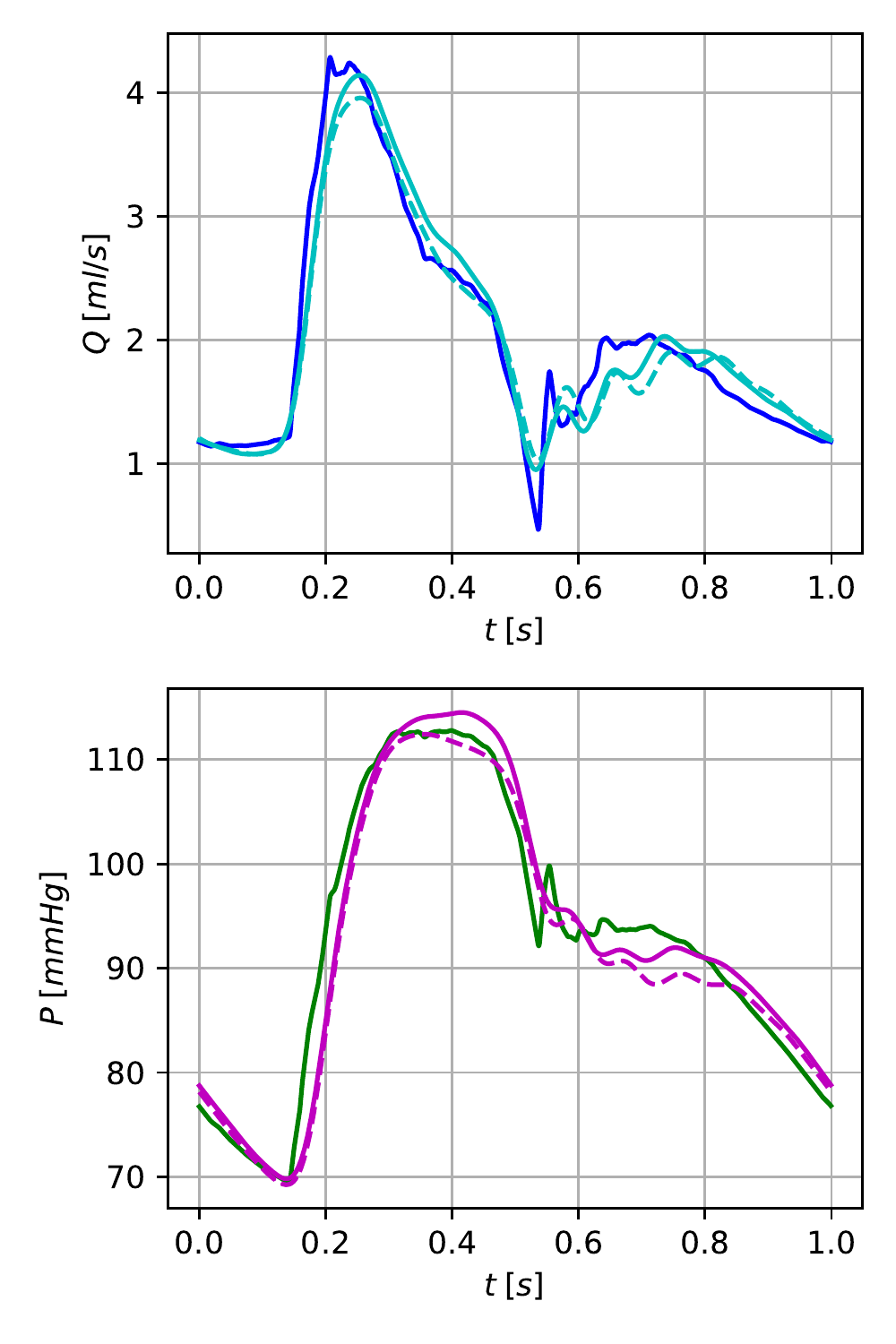}\label{fig:plot:test3:rightradial}} 
		\subfloat[][Right internal iliac] {\includegraphics[scale=0.55]{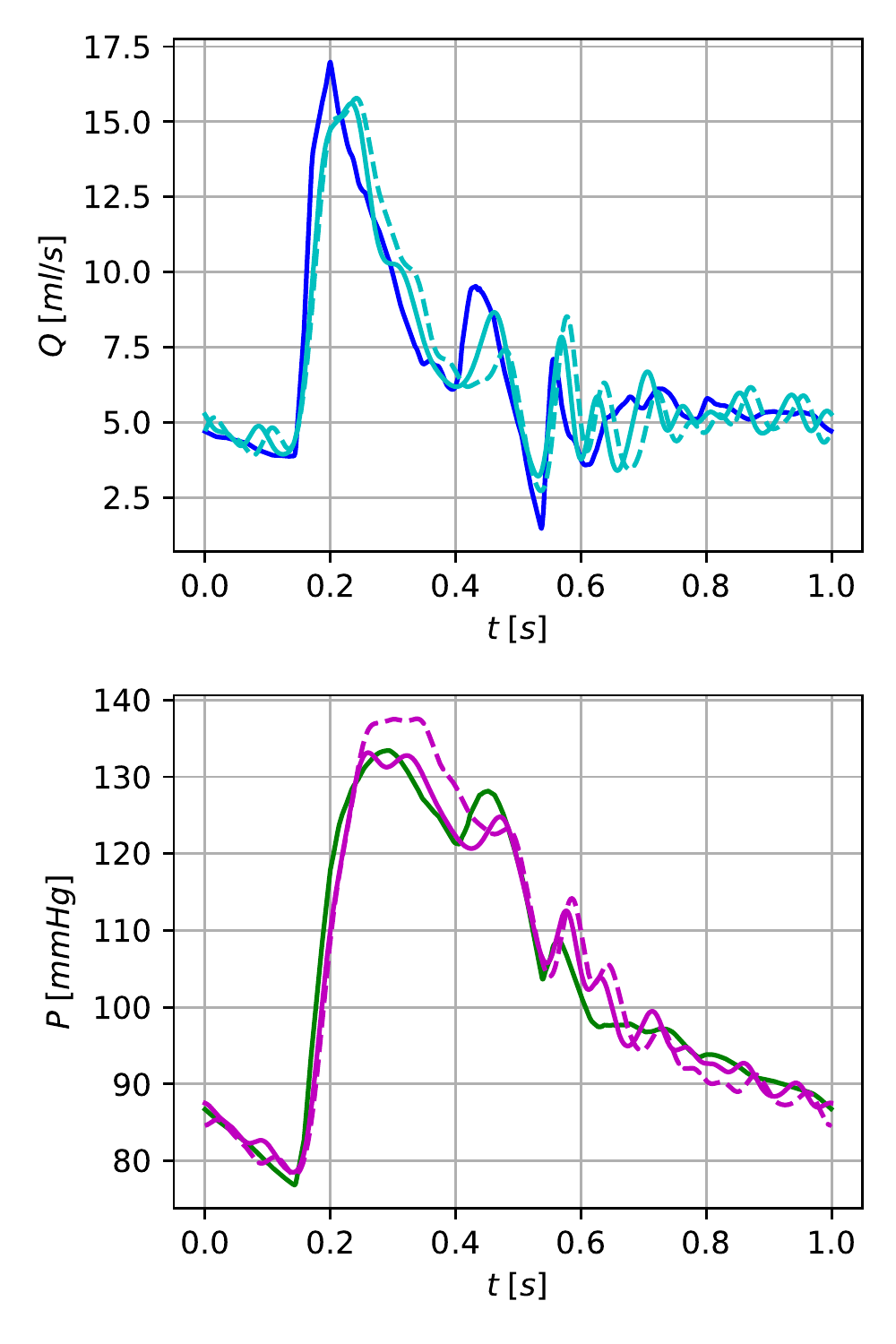}}	
		\caption{ADAN56 model. Comparison between 1D (at vessel midpoint), nonlinear 0D and linear 0D results in three second-generation vessels.}\label{fig:plot:test3:gen2}
	\end{center}
\end{figure}

\begin{figure}[h!]
	\begin{center}
		\subfloat[][Right posterior interosseous] {\includegraphics[scale=0.55]{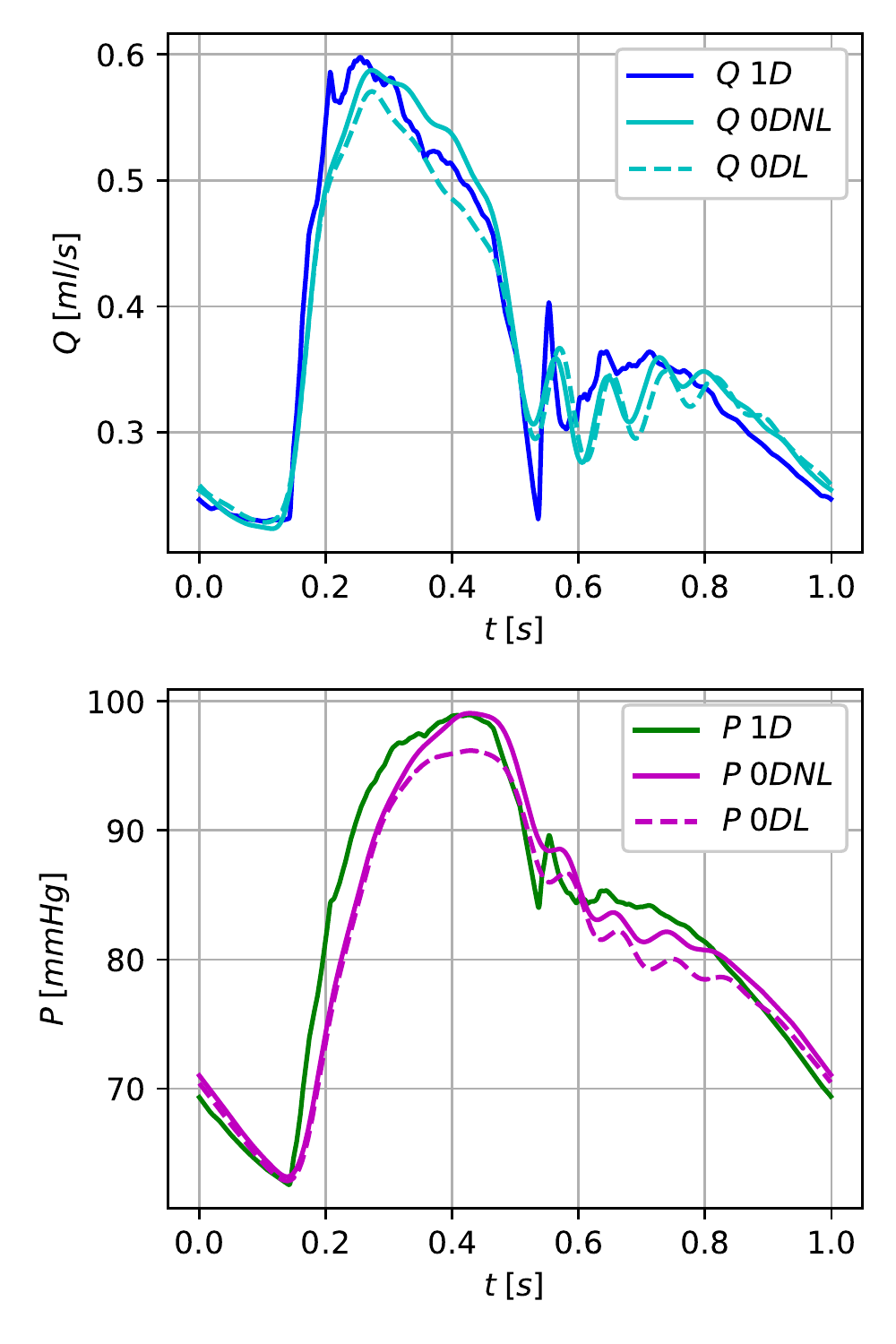}} 
		\subfloat[][Right femoral II] {\includegraphics[scale=0.55]{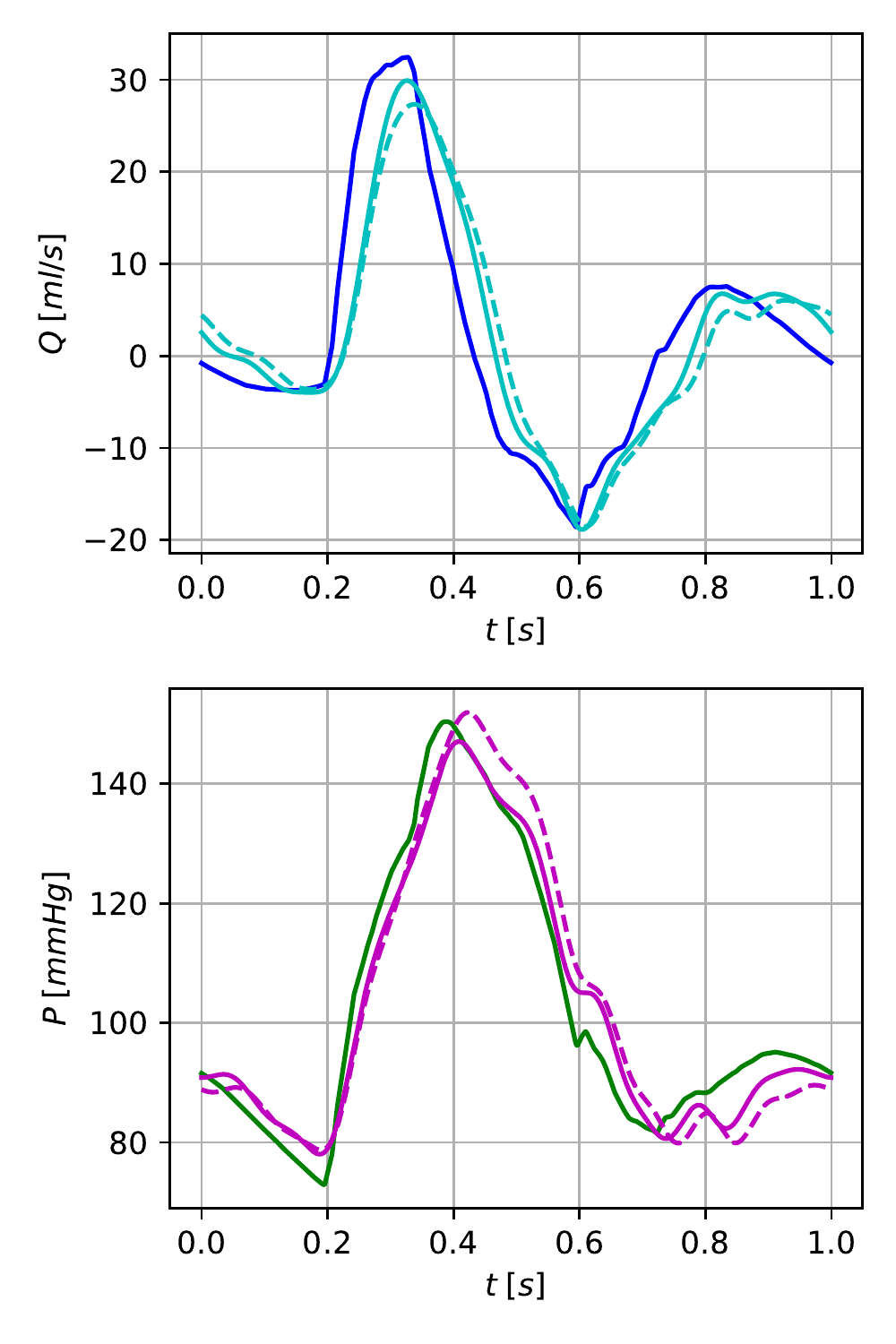}} 
		\subfloat[][Right anterior tibial] {\includegraphics[scale=0.55]{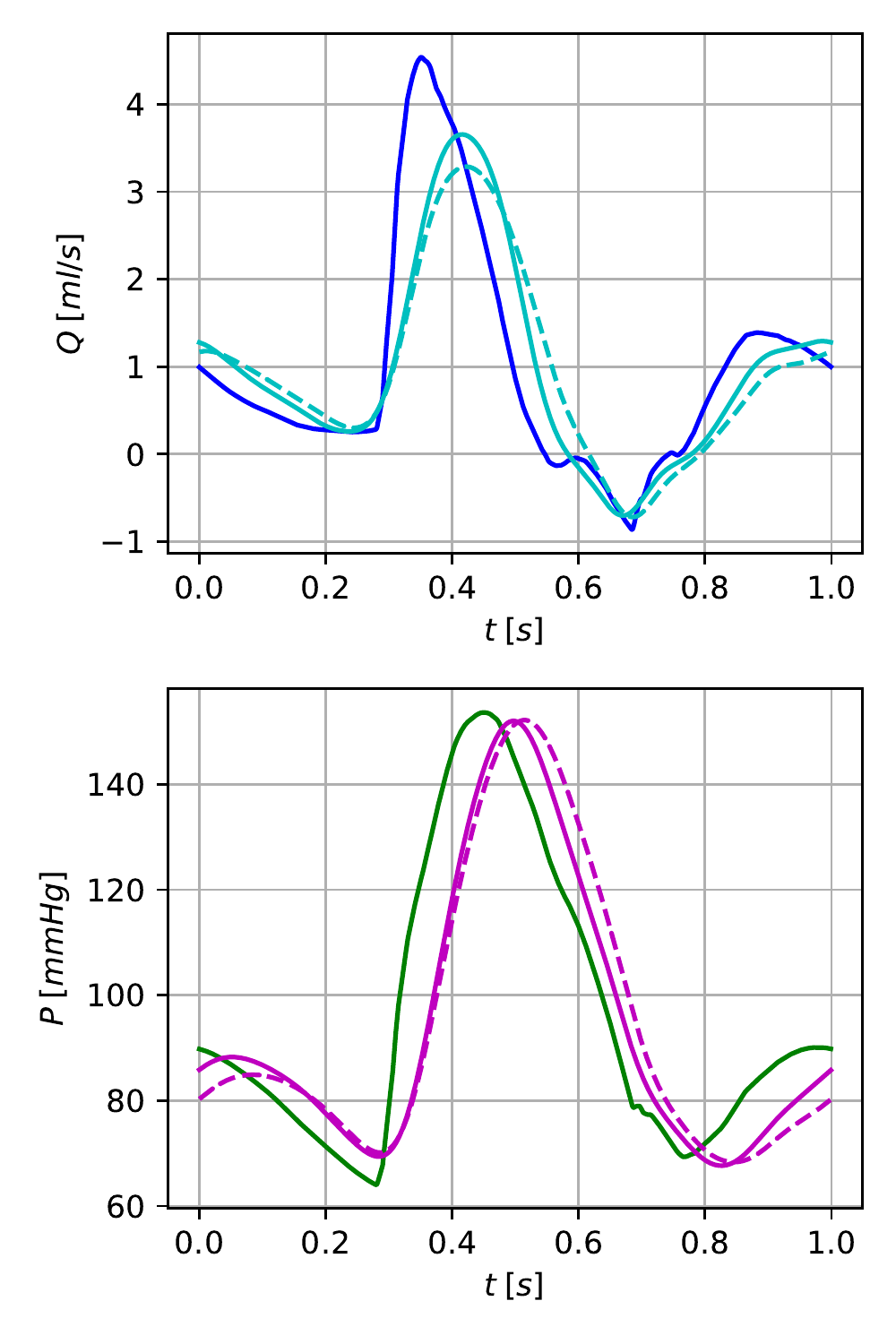}\label{fig:plot:test3:rightanteriortibial}}	
		\caption{ADAN56 model. Comparison between 1D (at vessel midpoint), nonlinear 0D and linear 0D results in three third- and fourth-generation vessels.}\label{fig:plot:test3:gen34}
	\end{center}
\end{figure}

The last benchmark model considered is a reduced version of the anatomically detailed arterial network model developed by Blanco \textit{et al.} \cite{Blanco:2014a,Blanco:2015a}. This model contains the largest 56 vessels of the human arterial system, as described in \cite{Blanco:2015a}. General properties of this arterial model are shown in Table \ref{table:test3}, while for the topology of the network and a complete set of parameters, we refer the reader to \cite{Boileau:2015a}. The inflow boundary condition $Q_{in}(t)$, inspired from one of the inflow signals reported in \cite{Murgo:1980a}, is taken from \cite{Boileau:2015a}. Terminal vessels are coupled to $RCR$ Windkessel terminal elements. In the tube law (\ref{tubelaw}), we set $P_0 = P_d = 10^5$ dyne/cm$^2$ and $A_0 = A_d$, from which the 1D/0D initial cross-sectional areas $A(x,0)/\widehat{A}(0)$ corresponding to the prescribed initial pressure are computed in all vessels. The vessel wall thickness is computed using the following empirical expression \cite{Blanco:2015a}
\begin{equation}
h = r_0\left[ \widetilde{a}\exp(\widetilde{b} r_0) + \widetilde{c}\exp(\widetilde{d} r_0) \right],
\end{equation}
where $r_0$ is the reference radius, related to $A_0$, $\widetilde{a} = 0.2802$, $\widetilde{b} = -5.053$ cm$^{-1}$, $\widetilde{c} = 0.1324$ and $\widetilde{d} = -0.1114$ cm$^{-1}$. As for the 37-artery network, also here we do not consider vessel tapering, but we assume all vessels to have a constant reference cross-sectional area $A_0$, related to the constant reference radius $r_0$ computed as the mean value between the proximal and distal radii given in \cite{Boileau:2015a}.\\
For the 0D simulations, the choice of the 0D vessel configurations is the same as that for the 37-artery network: the first vessel of the arterial model, that is the first portion of the aortic arch, is of $(Q_{in}, Q_{out})$-type; terminal vessels are of $(P_{in}, P_{out})$-type and coupled to $RCR$ Windkessel terminal elements; all other vessels of the network are modelled as two-split $(P_{in}, Q_{out})$-type 0D vessels. Results are shown for three aortic segments (aortic arch I, thoracic aorta III and abdominal aorta V), three first-generation vessels (right common carotid, right renal and right common iliac), three second-generation vessels (right internal carotid, right radial and right internal iliac) and three third-/fourth-generation vessels (right posterior interosseous, right femoral and right anterior tibial). Qualitative comparisons between nonlinear 0D, linear 0D and reference 1D solutions are illustrated in Figure \ref{fig:plot:test3:aorta} for the aortic segments and in Figures \ref{fig:plot:test3:gen1}-\ref{fig:plot:test3:gen34} for the vessels of first, second and third/fourth generations, respectively. Table \ref{table:errors:test3} summarizes the relative errors computed for both nonlinear and linear 0D results with respect to 1D results.\\
As for the 37-artery network considered in Section \ref{sec:numexp:test2}, for this arterial model we also observe that, even if some oscillations in the pressure and flow waveforms are amplified, there is globally a good agreement between the nonlinear 0D predicted results and the reference 1D results. Overall, pressure and flow profiles are reproduced with a reasonably good level of accuracy; the main features, shape and amplitude, of the waveforms are well-captured by the nonlinear 0D models, even if the differences between 1D and 0D results are not negligible. For instance, in Figure \ref{fig:plot:test3:gen34}, we note that for these vessels of third and fourth generations of bifurcations the 0D waveforms are slightly delayed to the right with respect to the 1D waves. However, once again, we note that generally, by preserving in the 0D models adopted the nonlinearity of the original 1D model, the 0D results are strongly improved with respect to the linear case also for this benchmark arterial network. This is quite evident in the pressure waves, where the systolic pressure is well-reproduced by the nonlinear 0D models. There are only few isolated cases where the linear 0D results seem to be better than the nonlinear ones, as, for example, in the right radial artery, where from Figure \ref{fig:plot:test3:rightradial} we can observe the systolic peak in pressure to be approximated better by the linear 0D solution, rather than by the nonlinear 0D solution.\\
The quantitative assessment presented in Table \ref{table:errors:test3} supports all these observations. From the linear to the nonlinear 0D results, the RMS relative errors computed with respect to the reference 1D results are reduced, suggesting that the nonlinear 0D models are able to better reproduce the overall dynamics of pressure and flow rate in all vessels. Furthermore, also in this case, the systolic relative errors for pressure are confirmed to be, in general, significantly smaller in the nonlinear 0D results with respect to the linear ones, illustrating the ability of the nonlinear 0D models to capture and reproduced the essential features of pressure and flow waveforms. For the nonlinear 0D results, RMS relative errors are all smaller than 12.0\% for the pressure and smaller than 18.0\% for the flow rate; systolic relative errors are all smaller than 4.0\% for the pressure and smaller than 10.0\% for the flow rate, with an exception in the right anterior tibial artery, where the peak in the flow rate waveform seems to be underestimated, as confirmed by Figure \ref{fig:plot:test3:rightanteriortibial}.\\
\begin{table}[h!]\footnotesize
	\centering
	\begin{tabular}{lc}
		\toprule
		\textbf{Property} & \textbf{Value} \\
		\midrule
		Blood density, $\rho$ & 1.040 g/cm$^3$ \\
		Blood viscosity, $\mu$ & 0.04 dyne$\cdot$s/cm$^2$ \\
		Velocity profile order, $\zeta$ & 2 \\
		Young's modulus, $E$ & 2.25$\cdot$10$^6$ dyne/cm$^2$ \\
		Diastolic pressure, $P_d$ & 10$^5$ dyne/cm$^2$ \\
		External pressure, $p_{ext}$ & 0 \\
		Outflow pressure, $P_{out}$ & 0 \\
		Initial velocity, $u(x, 0)/U(0)$ & 0\\
		Initial pressure, $p(x, 0)/P(0)$ & $P_d$\\
		\bottomrule
	\end{tabular}
	\caption{General parameters of the reduced ADAN56 model.}\label{table:test3}
\end{table}

\begin{table}[h!]\footnotesize
	\centering
	\renewcommand\arraystretch{1.2}
	\begin{tabular}{l | c | cccccc}
		\toprule
		\textbf{Vessel} & \textbf{0D model} & $\varepsilon_P^{\text{RMS}} (\%)$ & $\varepsilon_Q^{\text{RMS}} (\%)$ & $\varepsilon_P^{\text{SYS}} (\%)$ & $\varepsilon_Q^{\text{SYS}} (\%)$ & $\varepsilon_P^{\text{DIAS}} (\%)$ & $\varepsilon_Q^{\text{DIAS}} (\%)$ \\
		\midrule
		\multirow{2}{*}{Aortic arch I} & 0D-NL & 1.010 & 3.230 & 0.479 & 1.331 & -0.361 &  0.808 \\
		& 0D-L & 1.790 & 4.147 & 1.210 & 0.955 & -1.075 & -3.114\\
		\midrule
		\multirow{2}{*}{Thoracic aorta III} & 0D-NL & 1.200 & 4.131 & 0.649 & 2.303 & -1.332 & 5.122\\
		& 0D-L & 2.198 & 6.403 & 3.011 & 2.788 & -1.412 & 8.198\\
		\midrule
		\multirow{2}{*}{Abdominal aorta V} & 0D-NL & 1.735 & 9.719 & 1.740 & -7.392 & 0.732 & 13.722\\
		& 0D-L & 3.447 & 12.800 & 5.656 & -7.777 & -0.319 & 16.683\\
		\midrule
		\multirow{2}{*}{R common carotid} & 0D-NL & 1.111 & 6.857 & 1.272 & 6.378 &  -0.926 & 0.372 \\
		& 0D-L & 2.081 & 9.335 & 4.937 & 7.623 & -0.694 & -1.530 \\
		\midrule
		\multirow{2}{*}{R renal} & 0D-NL & 1.281 & 3.211 & 1.230 & 0.118 & 0.398 & -0.780\\
		& 0D-L & 2.765 & 5.940 & 4.733 & 3.339 & -0.429 & -2.143\\
		\midrule
		\multirow{2}{*}{R common iliac} & 0D-NL & 2.119 & 9.885 & 1.522 & -7.062 & 2.242 & 8.388 \\
		& 0D-L & 3.718 & 13.111 & 5.680 & -5.566 & 1.030 & 9.935\\
		\midrule
		\multirow{2}{*}{R internal carotid} & 0D-NL & 2.268 & 5.845 & 3.489 & 3.459 & -0.712 & 4.082 \\
		& 0D-L & 3.552 & 8.918 & 7.456 & 5.690 & -1.081 & 2.035 \\
		\midrule
		\multirow{2}{*}{R radial} & 0D-NL & 3.324 & 4.985 & 1.529 & -3.377 & 0.453 & 11.274\\
		& 0D-L & 3.586 & 5.514 & -0.319 & -7.678 & -0.358 & 13.079 \\
		\midrule
		\multirow{2}{*}{R internal iliac} & 0D-NL & 2.620 & 6.387 & -0.204 & -8.080 &  2.152 & 10.220\\
		& 0D-L & 4.012 & 8.659 & 3.104 & -7.054 & 1.711 & 7.333 \\
		\midrule
		\multirow{2}{*}{R post. inteross.} & 0D-NL & 3.275 & 4.132 & 0.134 & -1.734 & 0.998 & -0.9951 \\
		& 0D-L & 3.801 & 4.737 &  -2.797 & -4.547 & 0.469 & -0.155 \\
		\midrule
		\multirow{2}{*}{R femoral II} & 0D-NL & 4.888 & 16.954 & -2.224 & -7.812 & 6.997 & -0.868\\
		& 0D-L & 7.103 & 21.109 & 1.012 & -15.742 & 7.955 & 0.225 \\
		\midrule
		\multirow{2}{*}{R anterior tibial} & 0D-NL & 11.454 & 14.681 & -1.044 & -19.454 & 5.643 & 3.510 \\
		& 0D-L & 14.342 & 17.424 & -0.949 & -27.568 & 6.699 & 3.128 \\
		\bottomrule
	\end{tabular}
	\caption{ADAN56 model. Relative errors (in \%) for pressure and flow between for both nonlinear (0DNL) and linear (0DL) 0D results with respect to the 1D results (1D) at the midpoint of the vessel, computed according to the relative error metrics (\ref{metrics}).} \label{table:errors:test3}
\end{table}

\subsection{CPU times}\label{sec:numexp:CPU}

With this work, by constructing nonlinear 0D blood flow models to replace the original 1D model in not necessarily all, but in certain vessels of a network (according to some \textit{a-priori} model selection criteria), we aim at obtaining cheap simulations of blood flow in large and highly complex vascular networks, by drastically reducing their computational cost and execution time, while still preserving a good level of details in the results.\\
In Table \ref{table:CPUtimes} we compare the mean CPU times per cardiac cycle measured to perform the fully 1D and fully 0D (nonlinear and linear) simulations for the three arterial networks. 1D and 0D simulations were executed in Python, running on a Linux Ubuntu 18.04 machine with Gold Intel$^{\textregistered}$ Xeon$^{\textregistered}$ 6130 processor (2.1 GHz, 3.7 GHz Turbo, 16C/32T, 10,4 GT/s 2UPI, 22 MB Cache). From Table \ref{table:CPUtimes}, we first observe that the CPU time gradually increases with increasing network complexity. We note that, in the fully 1D simulations, the mean CPU times needed to solve the 37-artery network and ADAN56 model are approximately comparable, even if the second arterial network is larger than the first one. This observation can be justified by the fact that, because of the very fine mesh we are using to perform the 1D mesh-independent simulations, the time step $\Delta t$, computed according to the CFL condition, results to be smaller in the 37-artery network ($\Delta t\approxeq$ 6.05e-05 s) than in the ADAN56 model ($\Delta t\approxeq$ 1.3e-04 s), thus increasing the CPU time demanded to solve the first network. Furthermore, most importantly, the computational cost of the 0D simulations in all arterial networks is strongly reduced by orders of magnitude with respect to the corresponding 1D simulation, while the difference between the CPU times of nonlinear and linear 0D simulations is negligible. In Table \ref{table:CPUtimes}, in the columns of the CPU times corresponding to the 0D simulations, the speed-up gained when passing from the fully 1D simulation to the fully 0D simulation (either linear or nonlinear) is also indicated between brackets, measured as the ratio CPUtime1D/CPUtime0D. For the 37-artery network, the mean CPU time required to perform the fully 1D simulation with mesh-independent results is about 6245.1 s, while the CPU time demanded for the nonlinear, fully 0D simulation is significantly reduced to approximately 27.0 s, with a speed-up in the simulation of more than 230 times. Similarly, for the reduced ADAN56 model, the fully 1D simulation with mesh-independent results is performed with a mean CPU time of about 6492.7 s, while the CPU time demanded for the nonlinear, fully 0D simulation is decreased to approximately 84.9 s, with a speed-up in the simulation of more than 75 times. These results prove that the derived family of nonlinear 0D models may represent a powerful tool to improve the computational efficiency of blood flow simulations, while still preserving and well reproducing the main features of pressure and flow waveforms in networks of vessels. Of course, we expect this to be even more pronounced when this methodology will be applied to more complex networks, such as the global, closed-loop, multiscale cardiovascular model developed by M\"{u}ller and Toro \cite{Mueller:2014a,Mueller:2014b}, the ADAN model presented by Blanco \textit{et al.} \cite{Blanco:2014a,Blanco:2015a} and the comprehensive 1D model of the entire adult cardiovascular system reported by Mynard and Smolich \cite{MynardSmolich:2015a}. Indeed, in these networks, the number of vessels and spatial scales change drastically and, since we expect to always use 1D models for the larger vessels, the agreement in the results and the gain in terms of computational efficiency and execution time will be strongly improved by using hybrid 1D-0D networks.

We observe that, in many practical situations and applications, the mesh-independence of the 1D solution may not be necessary and cheaper 1D models may be employed. In this case, the computational cost of the fully 1D simulations would be considerably reduced with respect to that of the 1D mesh-independent model and thus, also the speed-up obtained by adopting the 0D models would be not as relevant. However, our main goal here is to investigate and assess the ability of the newly derived nonlinear 0D models to reproduce the physics and physical properties of blood flow through the vessel, when compared to the original 1D model. For this reason, it is important to have mesh-independent 1D results, to ensure that the physics of the problem is reproduced by the 1D model as accurately as possible. In addition, as pointed out above, we expect really impressive improvements in terms of computational efficiency when this methodology will be applied to much more complex networks and models, where the heterogeneity of spatial and/or temporal scales is relevant and also non mesh-independent simulations are really expensive.

Finally, we remark that the CPU time analysis presented in this section represents a preliminary study to get a first insight and estimate of the impact of using 0D models on the computational cost and efficiency of the simulations. More than the single CPU time values, the 1D/0D CPU time ratios are more relevant and informative about the speed-up obtained with respect to the reference 1D simulations. This analysis was performed using Python as programming language. We believe that what may drastically increment/decrement the speeds-up reported in this work is how efficiently 1D and 0D model discretizations are implemented, rather than the specific choice of the programming language. Indeed, in the case of serial code, we expect these ratios/speeds-up to be approximately maintained also when moving to another programming language (like C\ensuremath{++}), even if this would not be necessarily true in the case of a parallel code implementation.\\
\begin{table}[h!]\footnotesize
	\centering
	\renewcommand\arraystretch{1.2}
	\begin{tabular}{l|ccc}
		\toprule
		& \multicolumn{3}{|c}{\textbf{CPU time per card.cycle}}\\
		\midrule
		\multirow{2}{*}{\textbf{Test case}} & \textbf{1D sim.} & \textbf{0D-NL sim.} & \textbf{0D-L sim.}\\
		& & \textbf{(speed-up)} & \textbf{(speed-up)}\\
		\toprule
		Aortic bifurcation & 51.251 & 0.400 (128.096) & 0.332 (154.179) \\
		37-artery network & 6245.066 & 27.010 (231.215) & 19.566 (339.424) 19.923 (313.455) \\
		ADAN56 model & 6492.655 & 84.887 (76.486) & 78.486 (82.724) \\
		\bottomrule
	\end{tabular}
	\caption{Comparison between mean CPU times per cardiac cycle (in seconds) of the 1D, nonlinear 0D (0D-NL) and linear 0D (0D-L) simulations for the aortic bifurcation model, the 37-artery network and ADAN56 model. The speed-up gained in the 0D simulation (either linear or nonlinear) with respect to the 1D simulation is also reported between brackets, measured as the ratio CPUtime1D/CPUtime0D.} \label{table:CPUtimes}
\end{table}

\subsection{Analysis of relative contributions of single nonlinear terms}\label{sec:numexp:singlenonlinearities}

In this section, we perform the analysis of the relative contributions of considering as nonlinear single components of the 0D models (resistance, inductance and pressure-area relation) in determining pressure and flow waveforms and in producing more accurate 0D results. Indeed, in Sections \ref{sec:numexp:test2} and \ref{sec:numexp:test3}, we found that, for both the 37-artery network and ADAN56 model, overall considering nonlinear pressure-area relation and nonlinear parameters $L$ and $R$ in the 0D blood flow models strongly improves their ability of accurately predict pressure and flow waveforms, which are in good agreement with the reference 1D results. 
Here, we want to investigate and assess the contribution to the 0D model accuracy obtained by including only a single component as nonlinear and to compare these different possible scenarios with the previously considered ones where all components are either linear or nonlinear.\\
To this aim, we performed, for both the 37-artery network and ADAN56 model, a new set of fully 0D simulations, where in each simulation a single nonlinear term is activated in the 0D models (either the resistance $R$, the inductance $L$ or the pressure-area relation), while all the other components are kept linear. The obtained 0D results are then compared to the available fully-linear and fully-nonlinear 0D results.

\begin{figure}[h!]
	\begin{center}
		\includegraphics[scale=0.38]{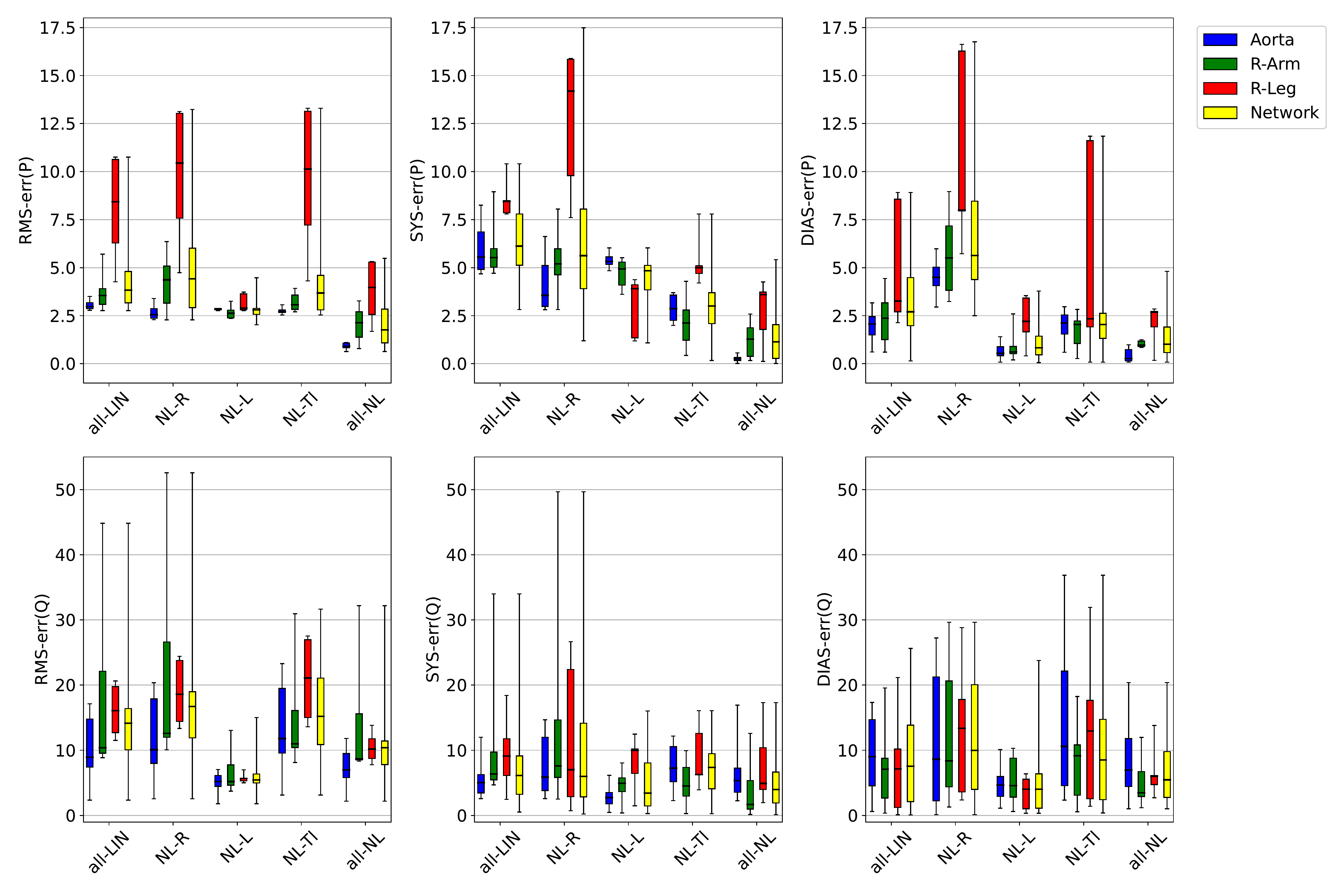}
		\caption{37-artery network. Boxplots comparing the errors computed in the different relative error metrics (\ref{metrics}) with respect to the reference 1D results for (in each boxplot, from left to right): fully-linear 0D results (all-LIN), 0D results with nonlinear resistance (NL-R), 0D results with nonlinear inductance (NL-L), 0D results with nonlinear pressure-area relation (NL-Tl) and fully-nonlinear 0D results (all-NL). In each boxplot, errors are compared over the entire network (yellow) and in several specific vascular regions: aorta (blue), right arm (green) and right leg (red). Relative errors for pressure $P$ are shown in the top row boxplots (from left to right boxplots: RMS, systolic and diastolic errors), while relative errors for flow rate $Q$ are reported in the bottom row boxplots. In each box, the box mid-line represents the median value, the box itself corresponds to the interquartile range (IQR), while the lines (whiskers) extending from the box indicate the minimum and maximum values.}\label{fig:singlenonlinear:network37}
	\end{center}
\end{figure}

\begin{figure}[h!]
	\begin{center}
		\includegraphics[scale=0.38]{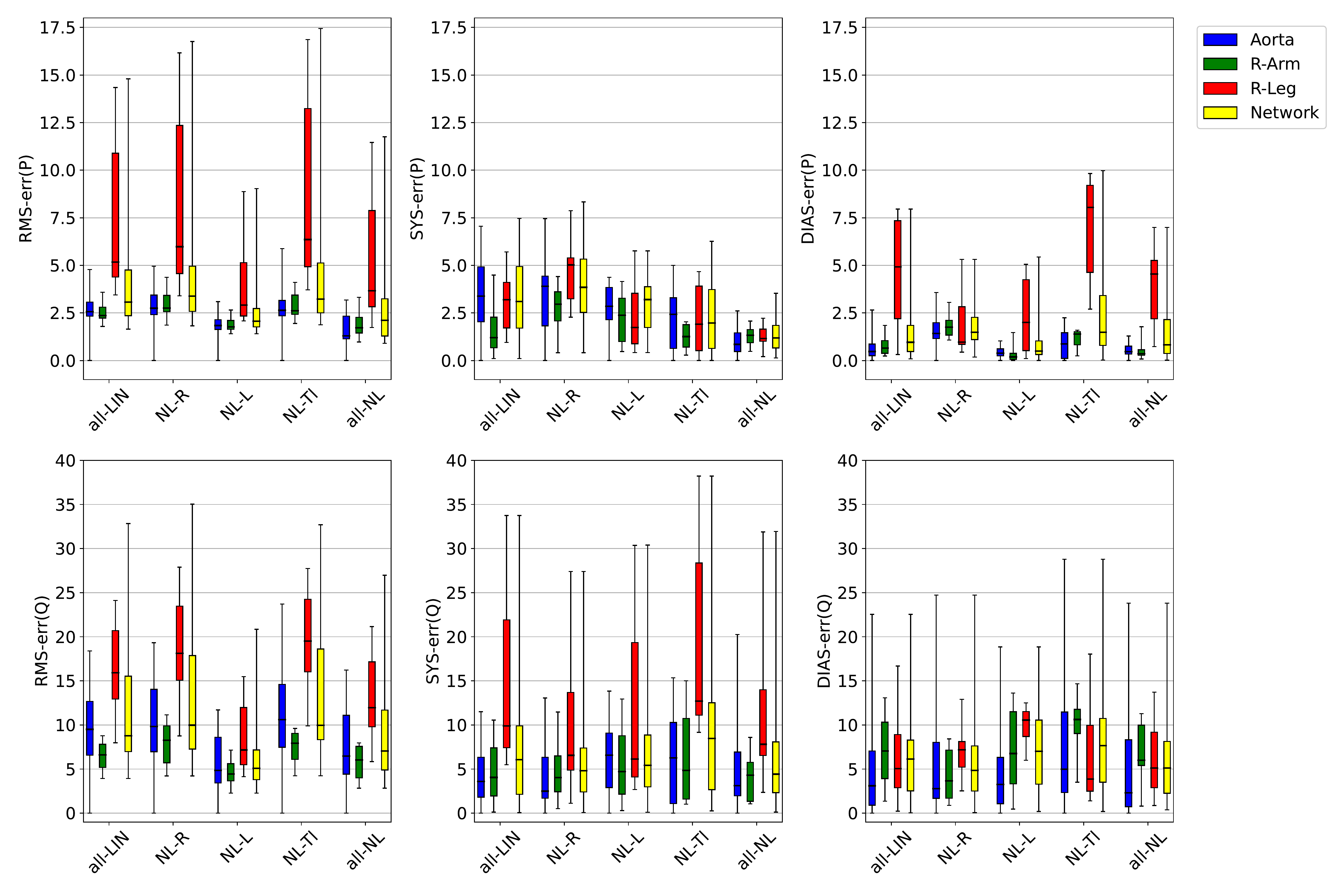}
		\caption{ADAN56 model. Boxplots comparing the errors computed in the different relative error metrics (\ref{metrics}) with respect to the reference 1D results for (in each boxplot, from left to right): fully-linear 0D results (all-LIN), 0D results with nonlinear resistance (NL-R), 0D results with nonlinear inductance (NL-L), 0D results with nonlinear pressure-area relation (NL-Tl) and fully-nonlinear 0D results (all-NL). In each boxplot, errors are compared in the entire network (yellow) and in several specific vascular regions: aorta (blue), right arm (green) and right leg (red). Relative errors for pressure $P$ are shown in the top row boxplots (from left to right boxplots: RMS, systolic and diastolic errors), while relative errors for flow rate $Q$ are reported in the bottom row boxplots. In each box, the box mid-line represents the median value, the box itself corresponds to the interquartile range (IQR), while the lines (whiskers) extending from the box indicate the minimum and maximum values.}\label{fig:singlenonlinear:adan56}
	\end{center}
\end{figure}

Figures \ref{fig:singlenonlinear:network37} and \ref{fig:singlenonlinear:adan56} display the error analysis carried out for the 37-artery network and ADAN56 model, respectively. For each type of 0D model (fully-linear, with nonlinear resistance, with nonlinear inductance, with nonlinear pressure-area relation and fully-nonlinear), the errors in the predicted 0D results are computed with respect to the reference 1D results in the different relative error metrics (\ref{metrics}).
In each boxplot, errors are compared over the entire network and in several specific vascular regions, namely the aorta, right arm and right leg, to examine whether specific features of the performance of the different 0D models emerge in any of the vascular districts considered. 
In all boxplots, for all error metrics, 0D model types and vascular regions, we report the median value, the interquartile range (IQR) and the minimum and maximum values of the error data. The IQR is computed as the difference between the 75th (upper) and 25th (lower) quartiles of each error dataset.

From this error analysis, several observations and outcomes are worthy to be pointed out.
Overall, the fully-nonlinear scenario improves the performance of the corresponding 0D blood flow models and the accuracy of the predicted results with respect to the fully-linear scenario; in few isolated cases, the fully-nonlinear and fully-linear results are comparable, with no significant gain obtained when all the nonlinearities are included in the 0D models, but also with no relevant worsening. 
Furthermore, when activating single nonlinear terms in the 0D models, the nonlinear inductance $L$ results to be the most determining nonlinear component which brings the major benefits to these models in terms of reduction of errors with respect to the reference 1D solution. Indeed, the contribution of this nonlinear parameter alone is significant in accurately defining both pressure and flow waveforms, with errors having small median values and small variability (small-sized IQR).
The nonlinear pressure-area relation also plays a considerable role in defining more accurate pressure results. This fact is evident when looking at the boxplots for the relative systolic errors, for both arterial networks considered.
Finally, it is worth mentioning that the main motivation for the introduction of nonlinear components, especially the pressure-area relation, regards the application of these 0D models to situations in which large deviations from a baseline state are to be considered, such as hypertension, an haemorrhage, a collapsed state or postural changes. Indeed, in these non physiological conditions, the combined contribution of the nonlinear pressure-area relation and of the nonlinear pressure-dependent parameters is expected to play an important role in accurately predicting the corresponding haemodynamic states.\\

\section{Concluding remarks and future work}\label{sec:conclusion}

In this paper we have extended the well-known family of lumped-parameter models for blood flow in a vascular segment, in order to construct 0D models that preserve important properties of the original 1D model. In particular, the main features characterizing these models are: (i) the nonlinearity in the pressure-area relationship and (ii) the nonlinearity in the model parameters $R$ and $L$, which depend on the average time-varying cross-sectional area. The different 0D vessel configurations depending on the data prescribed at the inlet and outlet of the vessel have been described; 0D junctions have been introduced to couple 0D vessels in a network.\\
In order to validate the newly derived family of nonlinear 0D models, fully 1D and fully 0D simulations have been performed and compared for three different benchmark arterial networks, showing that, even if there are discrepancies between 0D and 1D results, these nonlinear 0D models are able to capture and reproduce the essential features of pressure and flow waveforms. Furthermore, remarkable differences between linear and nonlinear 0D results have been observed in these physiologically normal cases and we expect that these differences will be even more relevant when addressing large deviations from the baseline haemodynamic states, like, for example, when simulating hypertension or haemorrhagic events.\\
Our ultimate goal is to construct hybrid 1D-0D networks, where certain vessels are treated 1D, while others are modelled as 0D vessels. Indeed, our work is primarily motivated by the attempt of facing the issues of computational efficiency and execution time, related to the modelling of blood flow in highly complex networks. Potential applications of the proposed framework will regard for instance the ADAN model developed by Blanco  \textit{et al.} \cite{Blanco:2014a,Blanco:2015a} and the global, closed-loop, multiscale mathematical model for the human cardiovascular system of M\"{u}ller and Toro \cite{Mueller:2014a,Mueller:2014b}. By approximating not necessarily all, but some vessels of such complex networks with these nonlinear 0D models, defining the different 0D vessel configurations, we expect that the computational cost and execution time of the resulting hybrid 1D-0D vessel networks will be dramatically reduced with respect to the fully 1D networks, even for long term simulations, and that the results will be significantly improved with respect to the fully 0D networks, still preserving the topology of the original 1D network and without losing any essential information with respect to the full 1D model. For this purpose, future work will include to find relevant \textit{a-priori} model selection criteria to define an adaptive model selection strategy that will allow to determine, given a network, the most suitable model to be used for each vessel of the network, either the 1D model or the newly derived nonlinear 0D models.\\
Moreover, we will consider purely 0D and hybrid 1D-0D couplings in order to find appropriate coupling conditions between 1D/0D vessels that do not require any restrictions on the admissible 0D vessel types converging at a junction.\\
Finally, in the proposed nonlinear lumped-parameter models we will incorporate additional aspects of the pressure-area relation, such as the viscoelastic properties of vessel walls, in order to derive more realistic and reliable 0D models.

\section*{Acknowledgements}
Beatrice Ghitti acknowledges the University of Trento for financing her Ph.D. studentship.\\
Lucas O. M\"{u}ller acknowledges funding from the Italian Ministry of Education, University and Research (MIUR) in the frame of the Departments of Excellence Initiative 2018–2022 attributed to the Department of Mathematics of the University of Trento (grant L. 232/2016) and in the frame of the PRIN 2017 project Innovative numerical methods for evolutionary partial differential equations and applications.\\
All authors warmly thank Prof. Marco Sabatini (Department of Mathematics, University of Trento, Italy) for very helpful discussions that contributed to the completion of the ODE stability analysis presented in this work.

\bibliographystyle{plain}

\appendix
\section{Stability analysis}\label{appendix:stabilityanalysis}

\subsection{First-order ordinary differential equations with periodic forcing function.}\label{appendix:stabilityanalysis:part1}

It is a well-known result that, given a linear homogeneous system of ODEs with constant coefficients, that is
\begin{equation}\label{homogeneous}
	\frac{d \bm{x}(t)}{dt} =  \bm{A} \bm{x}(t),
\end{equation}
the stability of the exact solution of such ODE system is determined by the real part of the eigenvalues of the coefficient matrix $\bm{A}$. Precisely, a necessary and sufficient condition for the ODE system (\ref{homogeneous}) to be \textit{asymptotically stable} is that all eigenvalues of $\bm{A}$ have strictly negative real part. If this is the case, then there exist positive constants $\alpha$, $\beta$ such that 
\begin{equation}\label{asympstability}
	\lVert e^{\bm{A} t} \rVert \leq \beta e^{-\alpha t}, \quad t \geq 0. 
\end{equation}
Proof of this result can be found in classical books on ordinary differential equations, see for instance \cite{Sanchez:1968a,Hale:1969a}.\\
Here, we want to prove that, when a periodic forcing function $\bm{b}(t)$ is added to system (\ref{homogeneous}) to obtain the inhomogeneous ODE system (\ref{matrixForm}), if the homogeneous part is asymptotically stable, then any solution of the complete non-homogeneous system will converge to a periodic solution as $t$ increases. Namely, the periodic solution of the inhomogeneous linear problem (\ref{matrixForm}) is an attracting stable solution for the system.\\
The exact solution of system (\ref{matrixForm}) obtained with the method of variation of parameters reads
\begin{equation}\label{exactsolution}
	\bm{x}(t) = e^{\bm{A} t} \bm{x}_0 + e^{\bm{A} t} \int_{0}^{t} e^{- \bm{A} s} \bm{b}(s)\, ds = e^{\bm{A} t} \bm{x}_0 + \int_{0}^{t} e^{\bm{A} (t-s)} \bm{b}(s)\, ds,
\end{equation}
where $\bm{x}(0) = \bm{x}_0$ and $e^{\bm{A} t}$ is the matrix exponential. Under the assumption of periodic forcing function $\bm{b}(t)$, \textit{i.e.} $\bm{b}(T_0)=\bm{b}(0)$ for some period $T_0>0$, we are able to obtain a particular solution of (\ref{matrixForm}), that is its periodic solution satisfying $\bm{x}(T_0)=\bm{x}(0)$, by choosing appropriately the initial condition $\bm{x}_0$. By enforcing this last condition, after straightforward calculations, we find that the initial condition $\bm{x}_0^P$ which provides a periodic solution of system (\ref{matrixForm}) is given by
\begin{equation}
	\bm{x}_0^P = \left( \bm{I} - e^{\bm{A} T_0} \right)^{-1} e^{\bm{A} T_0} \int_{0}^{T_0} e^{-\bm{A} s} \bm{b}(s)\, ds.
\end{equation}
At this point, we want to show that, if we start from any initial position different from the one which produces the periodic solution, the values of the corresponding solution will converge to the values of the periodic solution as $t$ increases. Let $\bm{x}^P$ be the periodic solution corresponding to the initial condition $\bm{x}_0^P$ and let $\bm{x}^{NP}$ be any other solution of system (\ref{matrixForm}) associated to some initial condition $\bm{x}_0^{NP}$. Then, by computing the norm of the difference between the two solutions and applying (\ref{asympstability}), the following inequality holds
\begin{equation}
	\lVert \bm{x}^P(t) - \bm{x}^{NP}(t) \rVert = \lVert e^{\bm{A} t}\left( \bm{x}_0^P - \bm{x}_0^{NP} \right) \rVert \leq \beta e^{-\alpha t} \lVert \bm{x}_0^P - \bm{x}_0^{NP} \rVert, 
\end{equation}
from which it clearly follows that $\beta e^{-\alpha t} \rightarrow 0$ as $t\rightarrow +\infty$, that is the sought result. In conclusion, if the forcing function $\bm{b}(t)$ is periodic and if the homogeneous part of system (\ref{matrixForm}) is asymptotically stable, then the exact solution of the original inhomogeneous problem will converge to the periodic solution of system (\ref{matrixForm}) as $t$ increases, for any admissible choice of the initial condition $\bm{x}(0) = \bm{x}_0 = \left[ V(0), Q(0) \right]^T$. As a consequence, the stability of the exact solution of the complete ODE system (\ref{matrixForm}) is determined by the real part of the eigenvalues of the coefficient matrix $\bm{A}$.

\subsection{First-order system of ordinary differential equations with null eigenvalue.}\label{appendix:stabilityanalysis:part2}

Here, as discussed is Section \ref{sec:stability}, we derive the additional assumption on the periodic forcing function $\bm{b}(t)$, which is needed to preserve the stability of the solution of an inhomogeneous ODE system of the form (\ref{matrixForm}) whose coefficient matrix $\bm{A}$ has a null eigenvalue (and all other eigenvalues with strictly negative real part). This assumption will ensure that the stability of the homogeneous system (which is stable, but not asymptotically stable) is inherited also by the corresponding inhomogeneous system, namely that any solution of (\ref{matrixForm}) for any admissible choice of the initial condition will converge to the periodic one as $t \to +\infty$.\\
This condition on the periodic forcing function is first derived in the scalar case for a single equation; then, the obtained result is extended to the case of a system of ODEs, specifically to system (\ref{QinQoutLin}) governing the linear $(Q_{in}, Q_{out})$-type 0D vessel.

\subsubsection{Scalar case.}
Let us start from the scalar case of a single ODE and consider the following mass conservation equation
\begin{equation}\label{singleEquationQinQout}
	\dfrac{d V}{d t} = Q_{in} - Q_{out}.
\end{equation}
For a scalar ODE, the coefficient matrix $\bm{A}$ in (\ref{matrixForm}) is reduced to a single constant coefficient $a$, which in the case of the above equation is equal to zero and plays the role of the null eigenvalue of matrix $\bm{A}$. As a consequence, the right-hand side of equation (\ref{singleEquationQinQout}) only depends on the scalar, periodic forcing function $b(t) = Q_{in}(t) - Q_{out}(t)$, where $Q_{in}(t)$ and $Q_{out}(t)$ represent the inflow entering the vessel and the outflow leaving the vessel, respectively, and are both periodic signals with the same period $T_0>0$. The exact solution of equation (\ref{singleEquationQinQout}) is
\begin{equation}
	V(t) = V(0) + \int_{0}^{t} b(s)\, ds,
\end{equation}
with initial condition $V(0)$. At this point, we are interested in studying the asymptotic properties of such a solution, namely the behaviour of $V(t)$ as $t$ approaches $+\infty$.\\
Let $t = nT_0 + \tau$, with $n\in \mathbb{N}$ sufficiently large and $0\leq \tau < T_0$. Then, we have
\begin{equation}
	\begin{split}
		V(t) = V(nT_0 + \tau) & = V(0) + \int_{0}^{nT_0+\tau} b(s)\, ds\\
		& = V(0) + \int_{0}^{nT_0} b(s)\, ds + \int_{nT_0}^{nT_0+\tau} b(s)\, ds\\
		& = V(0) + n \int_{0}^{T_0} b(s)\, ds + \int_{0}^{\tau} b(s)\, ds,
	\end{split}
\end{equation}
where, in the last equality, we have used the periodicity of the forcing term $b(t)$. In the above expression of the exact solution $V(t)$, the initial volume $V(0)$ is a fixed finite quantity, the second integral, $\int_{0}^{\tau} b(s)\, ds,$ is also finite, since it is the integral of a continuous periodic function over a fraction of its period $[0, T_0]$. Therefore, studying the asymptotic behaviour of the solution $V(t)$, that is
\begin{equation}
	\lim_{t\to +\infty} V(t)
\end{equation}
is completely equivalent to study
\begin{equation}
	\lim_{n \to +\infty} n \int_{0}^{T_0} b(s)\, ds.
\end{equation}
We distinguish the following three cases:
\begin{equation}
	\lim_{n \to +\infty} n \int_{0}^{T_0} b(s)\, ds = \begin{cases}
		+\infty & \text{if } \overline{b}>0,\\
		-\infty & \text{if } \overline{b}<0,\\
		0 & \text{if } \overline{b}=0,\\
	\end{cases}
\end{equation}
where we have set
\begin{equation}
	\overline{b} := \int_{0}^{T_0} b(s)\, ds.
\end{equation}
In other words, if the integral average of the periodic forcing function $b(t)$ over its period $[0,T_0]$ is different from zero, either positive or negative, then the solution of equation (\ref{singleEquationQinQout}) will diverge either to $+\infty$ or to $-\infty$, depending on the sign of $\overline{b}$, as $t\to +\infty$. On the other hand, if the integral average of $b(t)$ over its period is zero, then any solution $V(t)$ of (\ref{singleEquationQinQout}) will be bounded, namely
\begin{equation}\label{conditionperiodicity}
	\lim_{t\to +\infty} \big| V(t) \big| < +\infty \quad \Longleftrightarrow \quad \overline{b} = 0 \quad \Longleftrightarrow \quad \dfrac{1}{T_0} \int_{0}^{T_0} b(s)\, ds = 0,
\end{equation}
and will converge to the periodic one, that is
\begin{equation}
	V(nT_0 + \tau)  \xrightarrow[]{n\to +\infty}  V(0) + \int_{0}^{\tau} b(s)\, ds.
\end{equation}
We point out that this result is also physically consistent. Indeed, for equation (\ref{singleEquationQinQout}) the condition stated in (\ref{conditionperiodicity}) simply becomes
\begin{equation}
	\dfrac{1}{T_0} \int_{0}^{T_0} \left[ Q_{in}(s) - Q_{out}(s) \right]\, ds = 0 \quad \Longleftrightarrow \quad \int_{0}^{T_0} Q_{in}(s)\, ds = \int_{0}^{T_0} Q_{out}(s) \, ds,
\end{equation}
which means that, in order for the volume $V(t)$ not to constantly increase/decrease and asymptotically explode, the inflow entering the vessel and the outflow leaving the vessel must be balanced over each period.

\subsubsection{ODE system case.}
Let us now extend the previous result to the case of a system of ODEs, specifically to system (\ref{QinQoutLin}). For any admissible choice of the initial conditions $\bm{x}_0 = \left[ V(0), Q(0), V_d (0) \right]^T$, the general solution of (\ref{QinQoutLin}) is given by the formula in (\ref{exactsolution}), where the constant coefficient matrix $\bm{A}$ and the periodic forcing function $\bm{b}(t)$ are
\begin{equation}\label{matrixAvectorb}
	\bm{A} = \left[ \begin{matrix}
		0 & -1 & 0\\[2ex] 
		\dfrac{2}{C_0 L_0} & -\dfrac{R_0}{L_0} & -\dfrac{2}{C_0 L_0}\\[2ex] 
		0 & 1 & 0
	\end{matrix} \right],
	\quad
	\bm{b}(t) =  \left[ \begin{matrix}
		Q_{in}(t)\\[2ex] 
		0\\[2ex] 
		-Q_{out}(t)
	\end{matrix}\right].
\end{equation}
In general, if an $n\times n$ matrix $\bm{A}$ is diagonalizable, with linearly independent eigenvectors $\bm{v}_1, \bm{v}_2, \ldots, \bm{v}_n$ and corresponding eigenvalues $\lambda_1, \lambda_2, \ldots, \lambda_n$, then, denoted by $\bm{S}$ the matrix whose columns are the eigenvectors of $\bm{A}$, we have
\begin{equation}
	\bm{S}^{-1} \bm{A} \bm{S} = \left[ \begin{matrix}
		\lambda_1 & 0 & \cdots & 0\\
		0 & \lambda_2 & \cdots & 0\\
		\vdots & \vdots & \ddots & \vdots \\
		0 & 0 & \cdots & \lambda_n
	\end{matrix} \right] = \bm{D}.
\end{equation}
The exponential of the diagonal matrix $\bm{D}$ is easily written as
\begin{equation}
	e^{\bm{D}} = \left[ \begin{matrix}
		e^{\lambda_1} & 0 & \cdots & 0\\
		0 & e^{\lambda_2} & \cdots & 0\\
		\vdots & \vdots & \ddots & \vdots \\
		0 & 0 & \cdots & e^{\lambda_n}
	\end{matrix} \right] 
\end{equation}
and thus, it is straightforward to check that, under the assumption that $\bm{A}$ is diagonalizable, the corresponding matrix exponential $e^{\bm{A}t}$ can be computed as follows
\begin{equation}
	e^{\bm{A}t} = \bm{S} \left[ \begin{matrix}
		e^{\lambda_1 t} & 0 & \cdots & 0\\
		0 & e^{\lambda_2 t} & \cdots & 0\\
		\vdots & \vdots & \ddots & \vdots \\
		0 & 0 & \cdots & e^{\lambda_n t}
	\end{matrix} \right]  \bm{S}^{-1} = \bm{S} e^{\bm{D}t} \bm{S}^{-1}.
\end{equation}
Then, by introducing the above result in (\ref{exactsolution}), the exact solution of an inhomogeneous ODE system of the form (\ref{matrixForm}) with diagonalizable coefficient matrix $\bm{A}$ can be reformulated as 
\begin{equation}\label{exactsolutionNew}
	\bm{x}(t) = \left( \bm{S} e^{\bm{D}t} \bm{S}^{-1} \right) \bm{x}_0 + \int_{0}^{t} \left( \bm{S} e^{\bm{D} (t-s)} \bm{S}^{-1} \right) \bm{b}(s)\, ds,
\end{equation}
As observed at the beginning of this paragraph, if the coefficient matrix $\bm{A}$ has a zero eigenvalue and all other eigenvalues with strictly negative real part, then the homogeneous part of system (\ref{matrixForm}) is stable, but not asymptotically stable, with an attractive line of equilibrium points. Therefore, to study the asymptotic properties of the exact solution of the inhomogeneous system (\ref{matrixForm}) obtained by introducing a periodic forcing function $\bm{b}(t)$, we have to focus on the integral term in (\ref{exactsolutionNew}), to study its behaviour as $t$ tends to $+\infty$. We set
\begin{equation}
	\widetilde{\bm{S}}(t) := \bm{S} e^{\bm{D} t} \bm{S}^{-1}.
\end{equation}
By restricting to the order $n=3$ of the coefficient matrix of system (\ref{QinQoutLin}), after straightforward calculations, we get that the $i$th component of the integrand vector function in (\ref{exactsolutionNew}), obtained by multiplying the $i$th row of matrix $\widetilde{\bm{S}}$ by the vector forcing function $\bm{b}(t)$ and using the fact that $\lambda_1 = 0$, is given by
\begin{equation}
	\begin{split}
		( \widetilde{\bm{S}}(t-s) )_i \bm{b}(s) & = \sum_{j=1}^{3 } ( \widetilde{\bm{S}}(t-s) )_{ij} \bm{b}_j(s) \\
		& = \dfrac{1}{|\bm{S}|} \left[ \left( \bm{v}_{22} \bm{v}_{33} - \bm{v}_{23} \bm{v}_{32} \right) \bm{v}_{1i} + \left( \bm{v}_{23} \bm{v}_{31} - \bm{v}_{21} \bm{v}_{33} \right) \bm{v}_{2i} e^{\lambda_2 (t-s)} + \left( \bm{v}_{21} \bm{v}_{32} - \bm{v}_{22} \bm{v}_{31} \right) \bm{v}_{3i} e^{\lambda_3 (t-s)} \right] \bm{b}_1(s)\\
		& + \dfrac{1}{|\bm{S}|} \left[ \left( \bm{v}_{13} \bm{v}_{32} - \bm{v}_{12} \bm{v}_{33} \right) \bm{v}_{1i} + \left( \bm{v}_{11} \bm{v}_{33} - \bm{v}_{13} \bm{v}_{31} \right) \bm{v}_{2i} e^{\lambda_2 (t-s)} + \left( \bm{v}_{12} \bm{v}_{31} - \bm{v}_{11} \bm{v}_{32} \right) \bm{v}_{3i} e^{\lambda_3 (t-s)} \right] \bm{b}_2(s)\\
		& + \dfrac{1}{|\bm{S}|} \left[ \left( \bm{v}_{12} \bm{v}_{23} - \bm{v}_{13} \bm{v}_{22} \right) \bm{v}_{1i} + \left( \bm{v}_{13} \bm{v}_{21} - \bm{v}_{11} \bm{v}_{23} \right) \bm{v}_{2i} e^{\lambda_2 (t-s)} + \left( \bm{v}_{11} \bm{v}_{22} - \bm{v}_{12} \bm{v}_{21} \right) \bm{v}_{3i} e^{\lambda_3 (t-s)} \right] \bm{b}_3(s)
	\end{split}
\end{equation}
with $|\bm{S}|= \text{det}(\bm{S})$, for $i=1,\ldots,3$. As both eigenvalues $\lambda_{2,3}$ have strictly negative real part, the presence of the exponential factors $e^{\lambda_2 t}$ and $e^{\lambda_3 t}$ ensures the convergence of the corresponding terms as $t \to+\infty$, independently of any assumption on the periodic forcing function $\bm{b}(t)$. On the other hand, since $e^{\lambda_1 t} = e^{0} = 1$ for any $t\geq 0$, when integrating the corresponding terms without decaying exponential factor, namely
\begin{equation}\label{integralb}
	\int_{0}^{t} \dfrac{\bm{v}_{1i}}{|\bm{S}|} \left[ \left( \bm{v}_{22} \bm{v}_{33} - \bm{v}_{23} \bm{v}_{32} \right)  \bm{b}_1(s) + \left( \bm{v}_{13} \bm{v}_{32} - \bm{v}_{12} \bm{v}_{33} \right) \bm{b}_2(s) + \left( \bm{v}_{12} \bm{v}_{23} - \bm{v}_{13} \bm{v}_{22} \right) \bm{b}_3(s) \right]\, ds
\end{equation}
for $i=1,\ldots,3$, we need to introduce additional assumptions on the periodic forcing function $\bm{b}(t)$ in order for the above integral terms, and hence the solution of system (\ref{QinQoutLin}), not to explode as $t\to +\infty$.\\
As done in the scalar case, let $t = nT_0 + \tau$, with $n\in \mathbb{N}$ sufficiently large and $0\leq \tau < T_0$. Then, using explicit expressions of the components of $\bm{b}(t)$ given in (\ref{matrixAvectorb}), we have that the integral in (\ref{integralb}) can be rewritten as follows
\begin{equation}
	\dfrac{\bm{v}_{1i}}{|\bm{S}|} \int_{0}^{nT_0+\tau}  \left[ \alpha_1 Q_{in}(s) - \alpha_2 Q_{out} (s) \right]\, ds, \quad \text{for } i=1,\dots,3,
\end{equation}
where we have set
\begin{equation}\label{alpha1alpha2}
	\alpha_1 := \left( \bm{v}_{22} \bm{v}_{33} - \bm{v}_{23} \bm{v}_{32} \right), \quad \alpha_2 := \left( \bm{v}_{12} \bm{v}_{23} - \bm{v}_{13} \bm{v}_{22} \right).
\end{equation}
If we discard the multiplicative factor $\frac{\bm{v}_{1i}}{|\bm{S}|}$, where $\bm{v}_{1i}$ is the $i$th element of the eigenvector $\bm{v}_1$ associated to the null eigenvalue $\lambda_1$, for all the three components of the solution of system (\ref{QinQoutLin}) we obtain the same integral, which, by exploiting the periodicity of $Q_{in}(t)$ and $Q_{out}(t)$, becomes
\begin{equation}
	\begin{split}
		\int_{0}^{nT_0+\tau}  \left[ \alpha_1 Q_{in}(s) - \alpha_2 Q_{out} (s) \right]\, ds&  =
		n \int_{0}^{T_0}  \left[ \alpha_1 Q_{in}(s) - \alpha_2 Q_{out} (s) \right]\, ds\\ 
		& + \int_{0}^{\tau}  \left[ \alpha_1 Q_{in}(s) - \alpha_2 Q_{out} (s) \right]\, ds.
	\end{split}
\end{equation}
The second integral in the above expression is finite, since it is the integral of the linear combination of two continuous, periodic functions, with the same period $[0, T_0]$, over a fraction of their period. Then, by means of the same procedure used in the previous scalar case, also here we can observe that, in order for each component of the solution of system (\ref{QinQoutLin}) not to diverge either to $+\infty$ or to $-\infty$ as $t\to +\infty$, namely
\begin{equation}
	\lim_{t\to +\infty} |\bm{x}_i (t)| < +\infty \quad \text{for } i=1,\ldots,3,
\end{equation}
and to converge to the periodic solution, a necessary and sufficient condition is
\begin{equation}\label{integralQinQout}
	\int_{0}^{T_0}  \left[ \alpha_1 Q_{in}(s) - \alpha_2 Q_{out} (s) \right]\, ds = 0
	\quad \Longleftrightarrow \quad 
	\alpha_1 \int_{0}^{T_0}  Q_{in}(s)\, ds = \alpha_2 \int_{0}^{T_0}  Q_{out} (s)\, ds.
\end{equation}
We recall that this result is strictly valid only in the case where the coefficient matrix $\bm{A}$ is diagonalizable and, from linear algebra, it is well-known that an $n\times n$ matrix is diagonalizable if it has $n$ independent eigenvectors, which is true, for example, when the matrix has $n$ distinct eigenvalues.\\
It is easy to check that the matrix $\bm{A}$ of system (\ref{QinQoutLin}) given in (\ref{matrixAvectorb}) is diagonalizable in both cases $\Delta>0$ and $\Delta<0$, described in (\ref{eigenvalQinQout}). If $\Delta>0$, matrix $\bm{A}$ has three real and distinct eigenvalues
\begin{equation}
	\lambda_1 = 0, \quad \lambda_{2,3} = -\dfrac{R_0}{2L_0} \pm \dfrac{\sqrt{\Delta}}{2}
\end{equation}
and corresponding set of linearly independent eigenvectors
\begin{equation}
	\bm{v}_1 = \left[ \begin{matrix}
		1\\[1.5ex] 
		0\\[1.5ex] 
		1
	\end{matrix}\right], \quad 
	\bm{v}_2 =\left[ \begin{matrix}
		1\\
		\left( \dfrac{R_0}{2L_0} + \dfrac{\sqrt{\Delta}}{2}\right)\\
		-1
	\end{matrix}\right], \quad 
	\bm{v}_3 =\left[ \begin{matrix}
		1\\
		\left( \dfrac{R_0}{2L_0} - \dfrac{\sqrt{\Delta}}{2}\right)\\
		-1
	\end{matrix}\right];
\end{equation}
while, if $\Delta<0$, matrix $\bm{A}$ has three distinct eigenvalues, one equal to zero and the other two complex and conjugate, given by
\begin{equation}
	\lambda_1 = 0, \quad \lambda_{2,3} = -\dfrac{R_0}{2L_0} \pm i \dfrac{\sqrt{-\Delta}}{2}
\end{equation}
and corresponding set of linearly independent eigenvectors
\begin{equation}
	\bm{v}_1 = \left[ \begin{matrix}
		1\\[1.5ex] 
		0\\[1.5ex] 
		1
	\end{matrix}\right], \quad 
	\bm{v}_2 =\left[ \begin{matrix}
		1\\
		\left( \dfrac{R_0}{2L_0} + i \dfrac{\sqrt{-\Delta}}{2}\right)\\
		-1
	\end{matrix}\right], \quad 
	\bm{v}_3 =\left[ \begin{matrix}
		1\\
		\left( \dfrac{R_0}{2L_0} - i \dfrac{\sqrt{-\Delta}}{2}\right)\\
		-1
	\end{matrix}\right].
\end{equation}
Equipped with explicit expressions of the eigenvectors of the coefficient matrix $\bm{A}$, we can then compute factors $\alpha_1$ and $\alpha_2$ given in (\ref{alpha1alpha2}) and defining the linear combination between the inflow function $Q_{in}(t)$ and the outflow function $Q_{out}(t)$ in the integral (\ref{integralQinQout}). We obtain
\begin{equation}
	\alpha_1 = \alpha_2 = \begin{cases}
		-\sqrt{\Delta}, & \text{if } \Delta > 0,\\
		-i \sqrt{-\Delta}, & \text{if } \Delta < 0,\\
	\end{cases}
\end{equation}
so that condition (\ref{integralQinQout}) on the periodic forcing function $\bm{b}(t)$ simply becomes
\begin{equation}\label{condtionQinQout}
	\int_{0}^{T_0}  \left[ Q_{in}(s) - Q_{out} (s) \right]\, ds = 0
	\quad \Longleftrightarrow \quad 
	\int_{0}^{T_0}  Q_{in}(s)\, ds = \int_{0}^{T_0}  Q_{out} (s)\, ds,
\end{equation}
which is exactly the same physically consistent condition found in the scalar case to get asymptotically the periodic solution.

\end{document}